\documentclass[11pt]{article}
\usepackage{amsmath,amssymb,amsthm,bm}
\usepackage{dsfont}
\usepackage{verbatim} 
\usepackage{graphicx,rotating}
\usepackage{xcolor}
\usepackage{natbib}
\usepackage{placeins}
\usepackage{ulem} 
\usepackage{enumitem} 
\usepackage{pstricks,pst-node,pst-plot}
\usepackage{xr}
\usepackage[colorlinks,citecolor=blue,urlcolor=blue]{hyperref}

\newcommand{\figurepath}{}
\newcommand{\auxpath}{}

\newcommand{\R}{\mathbb{R}}
\newcommand{\dH}{ {\rm d}H}

\newcommand{\N}{\mathbb{N}}

\newcommand{\cL}{\mathcal{L}}

\newcommand{\iid}{ {\stackrel{i.i.d.}{\sim}} }
\newcommand{\Prb}{\mathbb{P}}
\newcommand{\E}{\mathbb{E}}

\newcommand{ \cB }{ \mathcal{ B } }

\newcommand{\var}{\mbox{\rm Var}}
\newcommand{\cov}{\mbox{\rm Cov}}
\newcommand{\frkc}{\mathfrak{c}}

\newcommand{\bmu}{ \boldsymbol{ \mu } }
\newcommand{\br}{ \boldsymbol{r} }
\newcommand{\bX}{ \boldsymbol{X} }
\newcommand{\bR}{ \boldsymbol{R} }
\newcommand{\bG}{ \boldsymbol{G} }
\newcommand{\btheta}{ \boldsymbol{\theta} }
\newcommand{\bfrkc}{ \boldsymbol{\mathfrak{c}} }

\newtheorem{theorem}{Theorem}
\newtheorem{lemma}{Lemma}
\newtheorem{proposition}{Proposition}
\newtheorem{corollary}{Corollary}

\theoremstyle{definition}
\newtheorem{definition}{Definition}

\theoremstyle{plain}

\theoremstyle{remark}

\newtheorem{remark}{Remark}

\begin{document}

\title{Functional delta residuals and applications to simultaneous confidence bands of moment based statistics}
\author{Fabian J.E. Telschow$^{1,}$\footnote{Corresponding author: fabian.telschow@hu-berlin.de, Rudower Chaussee 25, 12489 Berlin, Germany}\,, Samuel Davenport$^2$ and Armin Schwartzman$^{2,3}$\\[5mm]
	$^1$Institute of Mathematics, Humboldt Universit\"at zu Berlin \\
	$^2$Division of Biostatistics, University of California, San Diego\\
	$^3$Hal{\i}c{\i}o\u{g}lu Data Science Institute, University of California, San Diego
}
\date{\today}
\maketitle

\begin{abstract}
		Given a functional central limit (fCLT) for an estimator and a parameter transformation,
	we construct random processes, called functional delta residuals, which
	asymptotically have the same covariance structure as the limit process of the
	functional delta method.
	An explicit construction of these residuals for transformations of moment-based estimators
	and a multiplier bootstrap fCLT for the resulting functional delta residuals are proven.
	The latter is used to consistently estimate the quantiles of the maximum of the limit
	process of the functional delta method in order to construct asymptotically valid
	simultaneous confidence bands for the transformed functional parameters.
	Performance of the coverage rate of the developed construction,
	applied to functional versions of Cohen's d, skewness and kurtosis, is illustrated in simulations
	and their application to test Gaussianity is discussed.

\end{abstract}

\section{Introduction}
Recently there has been an increased interest in modeling functional data in the Banach space $ C(S) $ of continuous functions endowed with the supremum norm, as, unlike Hilbert space methods based on integral norms such as the $L^2$-norm, it allows for relevant differences to be localized and tested. This approach was pioneered in  \citet{Dette2020functional, Dette2020bio} for functional time series data. Results regarding functional data in the subspace of $C^3$-processes of the Banach space $ C(S) $ have notably been applied, primarily via the framework of random field theory, in order to control the familywise error rate for smooth test statistics over continuous compact domains (see e.g., \citet{AdlerBook, Worsley:2004, Taylor:2007}). Furthermore, analysis of functional data in the context of  $ C(S) $ is key in the construction of simultaneous confidence bands for the mean function or its derivatives \citep{Degras2011, Cao2012, Cao2014, Chang2017, Wang2019, Telschow:2019, Liebl2019}. Other emerging applications include simultaneous confidence bands for covariance functions \citep{Cao2016, Wang2020}, testing for equality of covariance functions using the supremum norm \cite{Guo2018}, Confidence Probability Excursion (CoPE) sets \citep{Sommerfeld:2018} and techniques to detect relevant differences in the mean and covariance functions of two samples \cite{Dette2020bio, Dette2020cov}.

Non-linear functional parameters, with the notable exception of the covariance function, have not yet received much attention in the context of $ C(S) $ random variables. This article will provide insight into this topic by introducing a construction of residuals -- called functional delta residuals -- which can be used to simulate from the limiting Gaussian field of a non-linear statistic over an arbitrary compact domain. We will use this to construct asymptotic simultaneous confidence bands for moment-based statistics, such as Cohen's d, skewness and kurtosis, which are differentiable functions of pointwise moments of the data.

\begin{figure}[h]
	\begin{center}
		\includegraphics[trim=0 0 0 0,clip,width=1.9in]{\figurepath 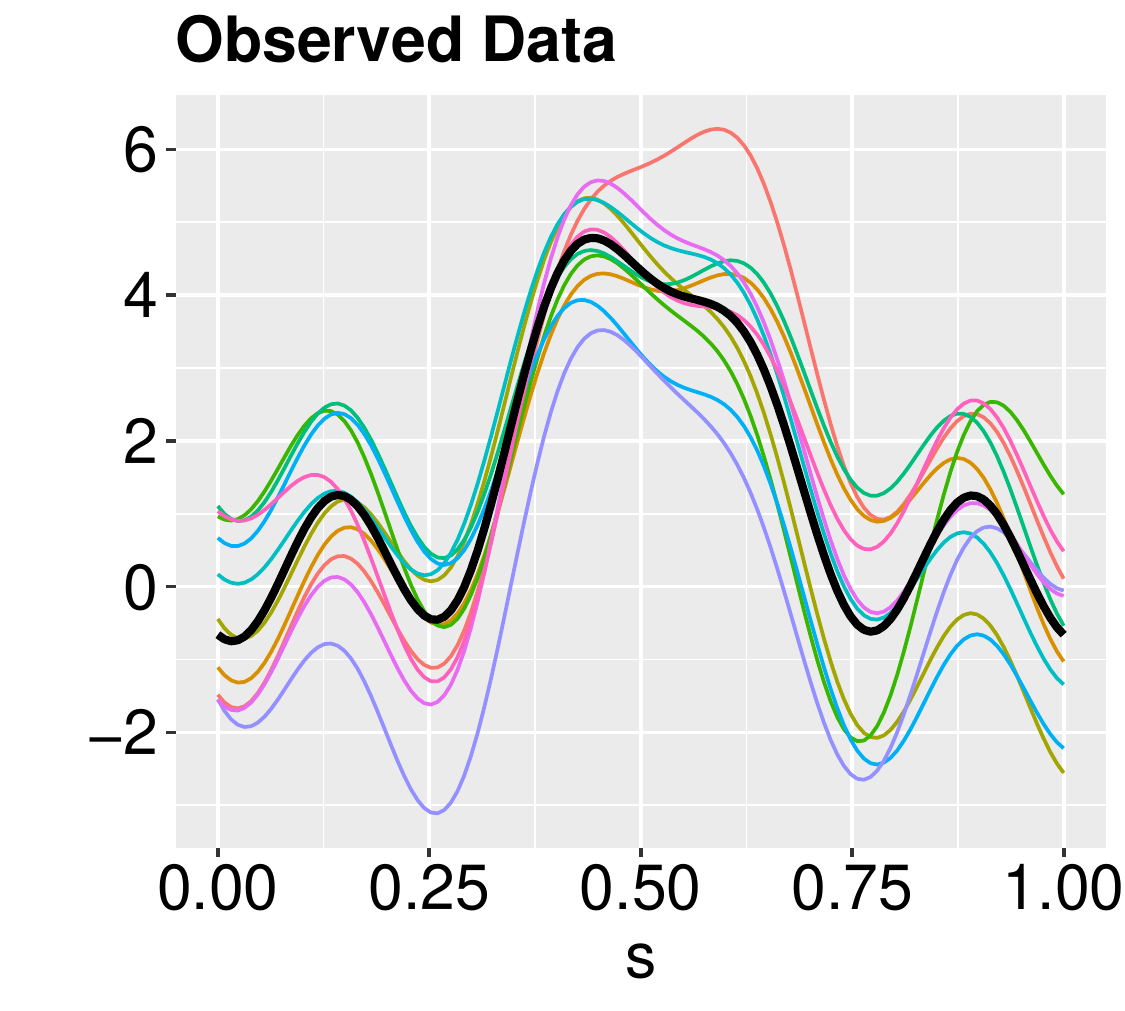}
		\includegraphics[trim=0 0 0 0,clip,width=1.9in]{\figurepath 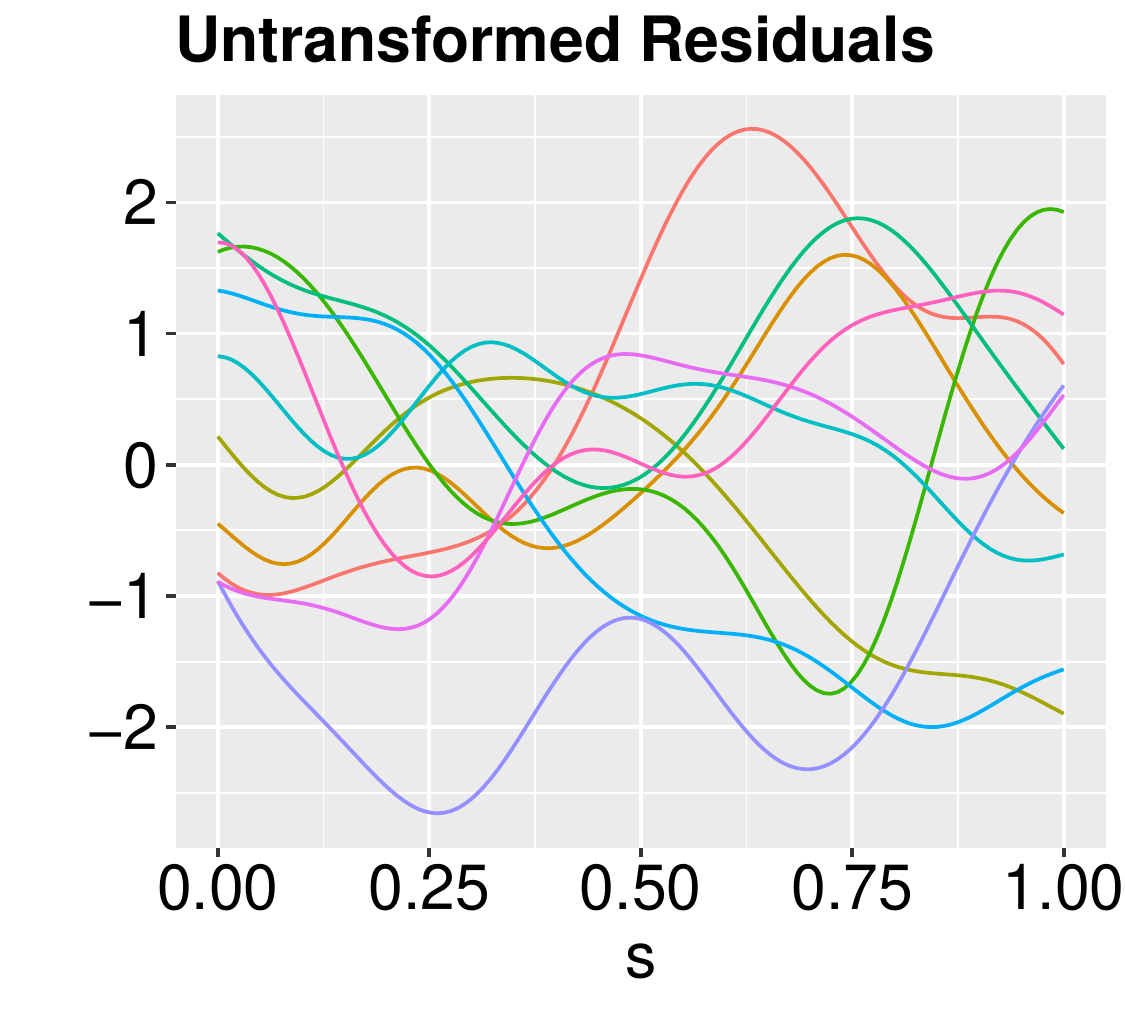}
		\includegraphics[trim=0 0 0 0,clip,width=1.9in]{\figurepath 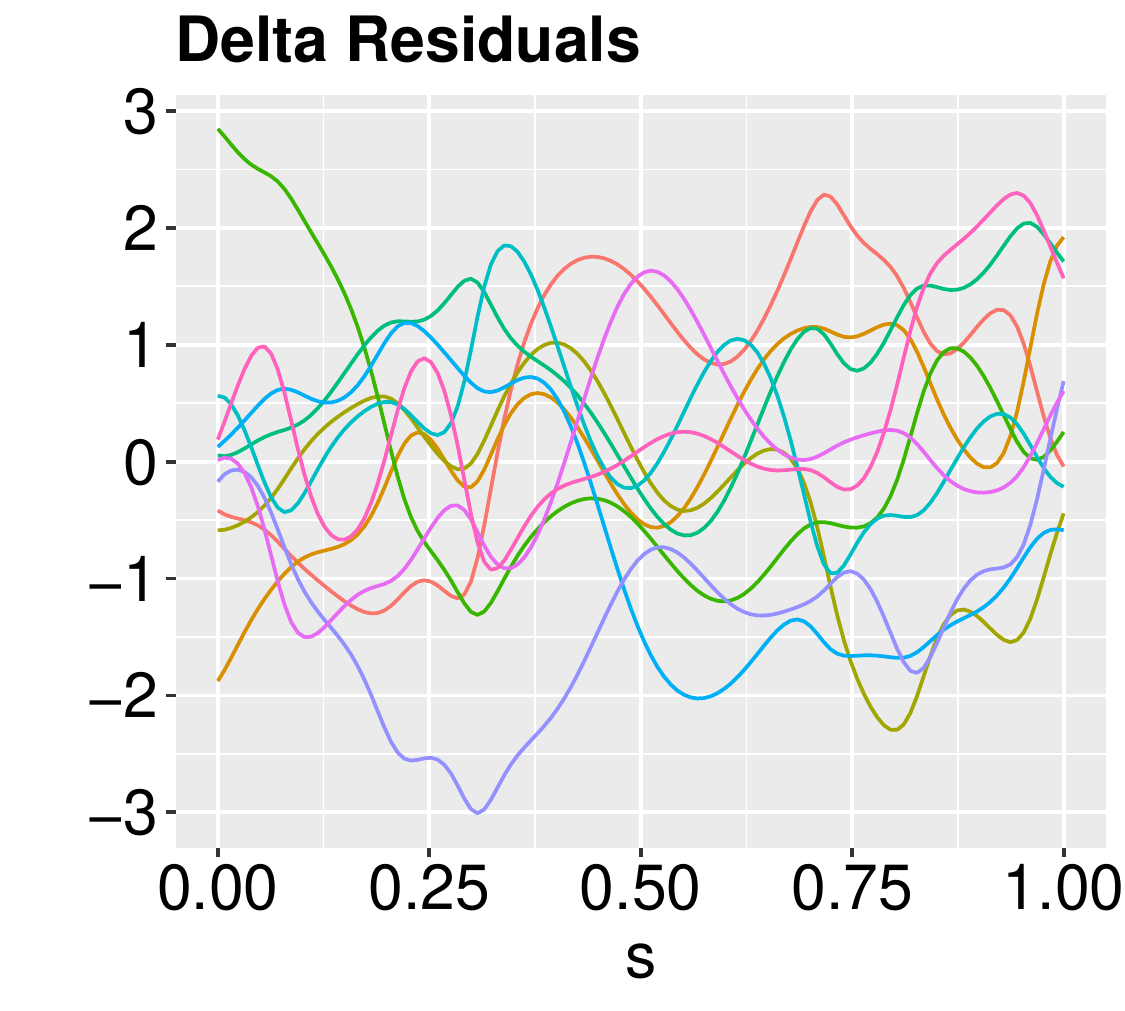}	
	\end{center}
		\vspace{-0.6cm}
	\caption{\textbf{Left:} samples of a Gaussian process with square exponential covariance function with a linear combination of Gaussian densities as  mean (bold black line). \textbf{Middle:} the untransformed residuals of this process. \textbf{Right:} the functional delta residuals of Cohen's $d$ of this process. \label{fig:SNR-samples}}
\end{figure}

\begin{figure}[h]
	\begin{center}
		\includegraphics[trim=5 5 0 0,clip,width=1.9in]{\figurepath 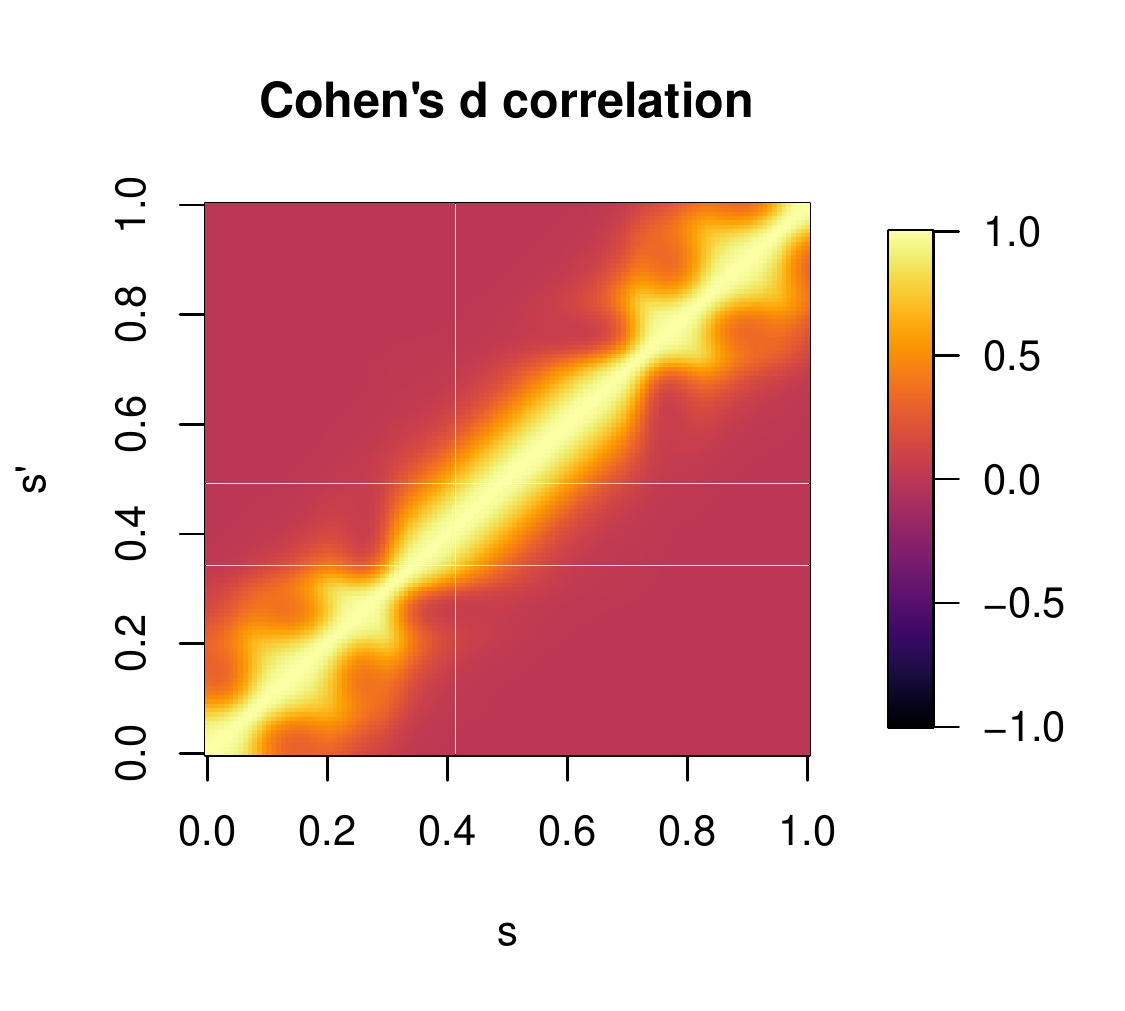}
		\includegraphics[trim=5 5 0 0,clip,width=1.9in]{\figurepath 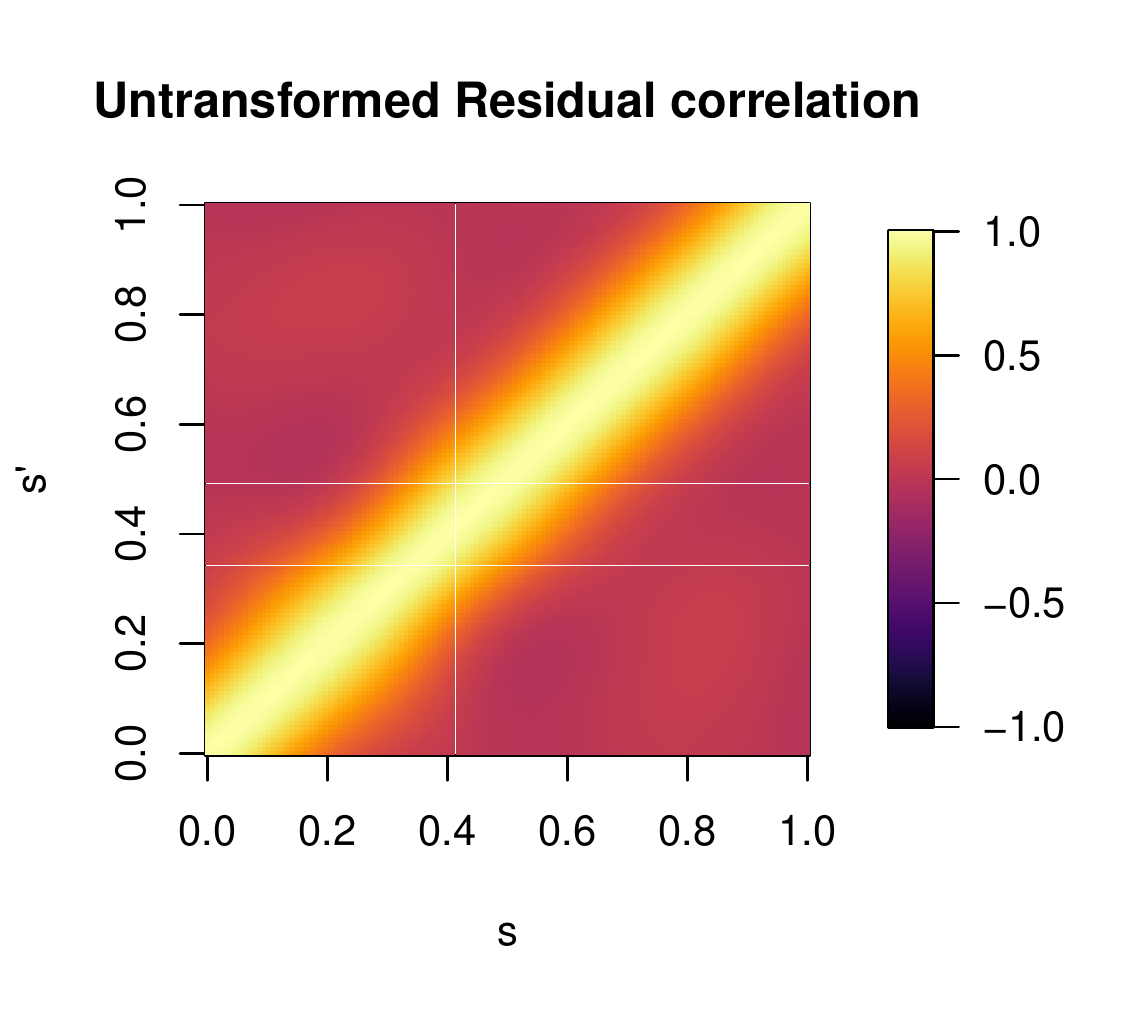}
		\includegraphics[trim=5 5 0 0,clip,width=1.9in]{\figurepath 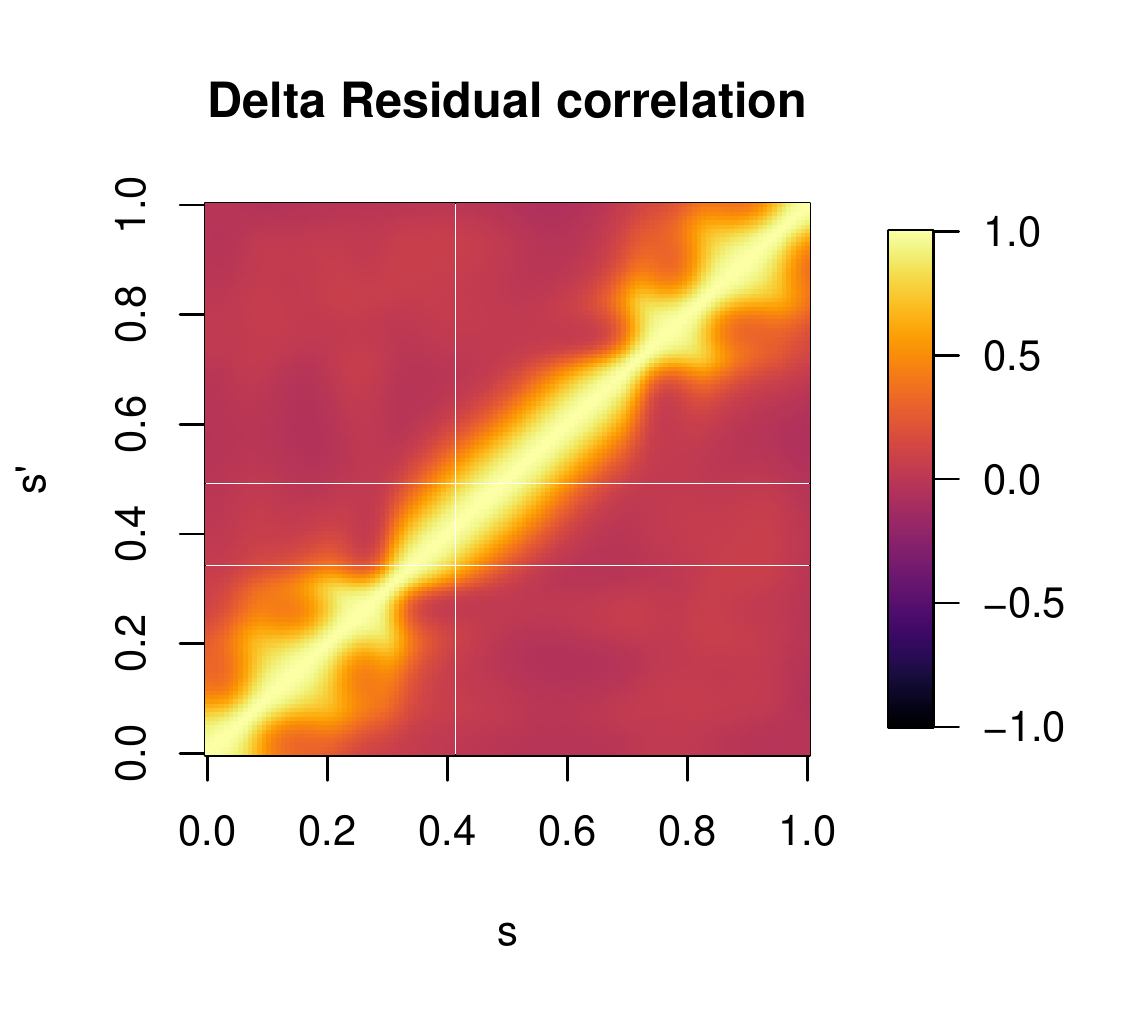}
	\end{center}
	\vspace{-0.6cm}
	\caption{\textbf{Left:}  the true asymptotic correlation function of Cohen's $d$ as given in Corollary \ref{cor:SNRfCLT} of the same process. \textbf{Middle:} correlation structure estimated from the sample correlation of a sample of untransformed residuals of size $ 100 $.
	\textbf{Right:} correlation structure of Cohen's $d$ estimated from the sample correlation of a sample of delta residuals of size $ 100 $.
	\label{fig:SNR-cov}}
\end{figure}

The motivation for our work comes from the following problem in spatial functional data analysis. \citet{Sommerfeld:2018}, in the context of climate data, and \citet{Bowring:2019a}, in the context
 of functional magnetic resonance imaging, study confidence statements about excursion sets of the mean function $ \mu(s) $ from a sample $ Y_1(s),\ldots,Y_N(s) $ observed from a signal plus noise model $Y(s) = \mu(s)+\varepsilon(s)$. Here $\varepsilon(s)$ is a stochastic error process with variance function $\sigma^2(s)$ and $s$ is a spatial index in a compact set $S\subset\R^D$. Their method requires estimation of the quantiles of the maximum of a limiting Gaussian process. These quantiles are estimated from the residuals $Y_n - \bar Y$ through a multiplier bootstrap \citep{Chang2009, Chang2017} or the Gaussian kinematic formula \citep{Worsley:2004,Adler:2009}. These methods successfully approximate the quantiles since the empirical covariance structure estimated from these residuals, asymptotically, has the same covariance structure as the limiting Gaussian process of $\sqrt{N}(\, \bar Y - \mu \,)$. Unsurprisingly, this approach no longer works if the object of interest is a non-linear transformation $H$ of the parameters, since the limiting Gaussian process of  $\sqrt{N}\big(\, H(\bar Y) - H(\mu) \,\big)$ in general has a different covariance structure than the residuals $Y_n - \bar Y$, compare Figure \ref{fig:SNR-samples} and \ref{fig:SNR-cov}. However, it is not immediately clear how to obtain residuals with the correct correlation structure, because applying the non-linear transformation to the residuals directly does not provide the correct correlation structure. In order to solve this problem we introduce the concept of functional delta residuals.

Our main result, Theorem \ref{thm:DeltaResiduals}, shows that functional delta residuals have the same asymptotic covariance structure as the limiting process from the fCLT for the transformed estimator. As an application, we derive the functional delta residuals for moment-based statistics such as the functional Cohen's $d$, skewness and kurtosis. 
Therefore we need to show that fCLTs in the Banach space of continuous functions hold for vectors of sample moments (Theorem \ref{thm:2DfCLT}) which in particular implies that they hold for moment-based statistics Corollary \ref{cor:momentfCLT}(a). The resulting functional delta residuals for moment-based statistics will be derived in Section \ref{scn:delta_res_cohensd}. In Theorem \ref{thm:DeltaResidualsBoots} we provide a conditional multiplier functional limit theorem using these residuals, which we prove using similar arguments to those used in \citet{Kosorok2003} and \citet{Chang2009}. In Theorem \ref{thm:SCB} we show that combining the previous results it is possible to construct asymptotic simultaneous confidence bands for moment-based statistics.
In particular the theory developed in our article provides a rigorous justification for the methods that we developed in \citet{Bowring:2020}. In that paper we extended our previous work on CoPE sets of signals in fMRI experiments to spatial inference using CoPE sets for Cohen's $d$, a statistic which important for measuring the power of a test, see \cite{Davenport:2020}. Moreover, we demonstrate how a transformation of pointwise skewness and kurtosis can be used to test Gaussianity of $C(S)$-valued samples. This uses transformations which transform the skewness and
kurtosis estimators to have approximately standard normal distributions \cite{DAgostino:1990}. In order to incorporate these results it is necessary to extend the theory to
work also for transformations depending on $H_N$. This mainly requires an extension of the Delta method to transformations
depending on $N$ which is proven in \ref{app:auxlemmas}.


The paper is organized as follows: Section \ref{scn:GeneralDeltaResiduals} introduces the general concept of functional delta residuals. Section \ref{scn:MomentDeltaResiduals} shows how to apply the general concept to moment-based statistics and includes the statements and proofs of our main results. Section \ref{scn:Simulations} studies the performance of simultaneous confidence bands for Cohen's $d$, skewness, kurtosis and certain transformations of skewness and kurtosis.

The proposed methods for simultaneous confidence bands are implemented in the \textsf{R}-package SIRF (Spatial Inference for Random Fields) available at \url{https://github.com/ftelschow/SIRF}. Code reproducing the simulations is available at \url{https://github.com/ftelschow/SIRF/DeltaResiduals}.

\section{Functional Delta Residuals}\label{scn:GeneralDeltaResiduals}

In this section we introduce the construction of functional delta residuals.
Throughout the article $ C^k\big(\, S, \R^P \,\big) $, $ k, P \in \mathbb{ N } $, denotes the
space of $ k $-times continuously differentiable functions with values in $ \R^P $
and domain $ S $.
For ease of readability $ C^k\big(\, S, \R \,\big) $ will be denoted by $ C^k(S) $.
We develop delta residuals in the framework of the Banach space $ C\big(\, S, \R^P \,\big) $
of continuous functions with values in $ \R^P $ over a compact domain $ S \subseteq \R^D $, $ D \in \mathbb{N} $. However the concept can be extended to other Banach spaces of functions and more general domains.
The norm on $ C\big(\, S, \R^P \,\big) $ is the maximum norm
$ \Vert f \Vert_\infty = \max_{ s \in S} \vert\, f( s ) \,\vert $, where $ \vert \cdot \vert $
denotes the standard norm on $ \R^P $.
The notation "$ \rightsquigarrow $" will denote weak convergence in $ C\big(\, S, \R^P \,\big) $, while
a bold symbol denotes a vector. Moreover, $ \boldsymbol{v}^T \in \R^{ 1 \times P } $ denotes the transpose of a column vector $ \boldsymbol{v} \in \R^P $. Given a function $f$ from
$ S \times S $ to $ \mathbb{R} $ we will write $f(\, s, s' \,)$ to refer to $ f $ evaluated at $s,s'\in S$ without explicitly defining $s,s'$. 
under this more general setup
Since a purely formal treatment hides the basic idea behind delta residuals, we
motivate them with a special case.
Let $ \lbrace X_{n} \rbrace_{n \in \mathbb{N}} $ be a sequence of
random processes in $ C(\,S) $ such that all elements are independent and identically
distributed as $ X $ with $ \frkc(\, s, s' \,) = \cov[\, X(s), X(s') \,] < \infty $ for all
$ s,s' \in S $ and $ \mu = \E[ X ] $. Assume that this array satisfies a functional CLT, i.e.,
\begin{equation}\label{eq:eq1}
	\sqrt{ N } \left(\, \bar X_N - \mu \,\right)
		= \frac{ 1 }{ \sqrt{ N } } \left(\, \sum_{ n = 1 }^N X_{n} - \mu \,\right) \rightsquigarrow G\,,
\end{equation} 
where $ \bar X_N = N^{-1} \sum_{ n = 1 }^N X_{ n } $ is the sample mean and $ G $ is a
tight zero-mean Gaussian process in $ C(\,S) $ with covariance function $ \frkc $.
Assume further that the residuals $ R_{ N, n } = X_{ n } - \bar X_N $
satisfy
\begin{equation}
\begin{aligned}
	\lim_{ N \rightarrow \infty } \frac{ 1 }{ N } \sum_{ n = 1 }^N R_{ N, n }( s )\, R_{ N, n }( s' ) 
 			&= \frkc(\, s, s' \,)\,.\label{eq:covGen}
\end{aligned}
\end{equation}
almost surely for all $ s, s' \in S $. For reasons which become clear in the next step we call
the residuals $ R_{ N, 1 }, \ldots , R_{ N, N } $ the \emph{untransformed residuals}.

Let $ H \in C^1( \R ) $ and denote with $ \dH_{x} $ the derivative of $ H $ at $ x \in \R $.
Suppose we are interested in inferring on the function $ s \mapsto H\big( \mu( s ) \big) $. Then
Equations \eqref{eq:eq1} and \eqref{eq:covGen} imply that the transformed processes
$
\tilde R_{ N, n }( s ) = \dH_{ \bar X_N( s ) }R_{ N, n }( s )\,,
$
which we call \emph{functional delta residuals} (\emph{delta residuals} for short), satisfy
\begin{equation*}
\begin{aligned}
	\lim_{ N \rightarrow \infty } \frac{ 1 }{ N } \sum_{ n = 1 }^N \tilde R_{ N, n }( s )\, \tilde R_{ N, n }( s' ) 
 &=  \lim_{ N \rightarrow \infty }\dH_{ \bar X_N( s ) }\dH_{ \bar X_N( s' ) } \frac{ 1 }{ N } \sum_{ n = 1 }^N R_{ N, n }( s )\, R_{ N, n }( s' ) \\
 &=\dH_{ \mu(s) }\frkc(\, s, s' \,)\dH_{ \mu(s') }
\end{aligned}
\end{equation*}
almost surely for all $ s, s' \in S $.
They can thus be used to approximate the covariance structure of $\tilde G$, which is
the Gaussian limiting process appearing in the fCLT obtained from the delta residuals
method, since
\begin{align*}
	\sqrt{ N } \Big(\, H \big( \bar X_N \big) - H \big( \mu \big) \,\Big)
								 \rightsquigarrow \tilde G = \dH_{ \mu } G\,.
\end{align*}

The next result generalizes the outlined concept of functional delta residuals to
arbitrary functional estimators $\hat\theta_N$ and clarifies some of the underlying
necessary conditions.

\begin{theorem}\label{thm:DeltaResiduals}
 Let $ N \in \N $ and $ \hat\btheta_N \in C\big(\, S, \R^P \,\big) $ be an estimator of a parameter
 $ \btheta \in  C\big(\, S, \R^P \,\big) $ such that as $ N \rightarrow \infty $
 \begin{equation}\label{eq:fclt}
  \sqrt{ N } \left(\, \hat\btheta_{N} - \btheta \,\right) \rightsquigarrow  \bG\,,
 \end{equation}
weakly in $ C\big(\, S, \R^P \,\big) $, where $ \bG $ denotes a zero-mean Gaussian process on
$ C\big(\, S, \R^P \,\big) $ with covariance function $ \bfrkc $. Let
$ H \in C^1\big(\, \R^P, \R^{P'} \,\big) $ then
\begin{enumerate}
	\item[(a)] the functional delta method implies that
		\begin{equation*}
	  		\sqrt{ N }\Big(\, H( \hat\btheta_{N} ) - H( \btheta ) \,\Big)
			  		\rightsquigarrow \dH_{ \btheta } \bG = \tilde \bG\,,~ N \rightarrow \infty
		\end{equation*}
		where $ \dH_{x} $ is the derivative of $ H $ at $ x \in \R^{P} $ and $\tilde \bG$ is a zero-mean Gaussian process with covariance	$ \tilde{ \frkc }
		(\, s, s' \,) = \dH_{ \btheta( s ) } \bfrkc\big(\, s, s' \,\big)\dH_{ \btheta( s' ) }^T $,
	\item[(b)] if $ \lbrace \bR_{ N, n }: N \in \mathbb{ N },
		  1 \leq n \leq N \rbrace $ is a triangular array of random processes in $ C\big(\, S, \R^P \,\big) $ such that $ \sum_{ n = 1 }^N  \bR_{ N, n } = 0 $ and
 		  \begin{equation}\label{eq:conv}
   				\lim_{ N \rightarrow \infty } N^{-1} \sum_{ n = 1 }^N \bR_{ N, n }(s)\, \bR_{ N, n }^T( s' )
   		 				= \bfrkc\big(\, s, s' \,\big) \in \R^{ P \times P } \,.
 			\end{equation}
			uniformly in probability, then the
			\emph{functional delta residuals}
			 $ \tilde \bR_{ N, n }( s ) =
					\dH_{ \hat \btheta_{ N }( s ) } \bR_{ N, n }( s ) $, $ n = 1, \ldots, N $, satisfy
			\begin{align*}
				\lim_{ N \rightarrow \infty } N^{ -1 } \sum_{ n = 1 }^{ N } \tilde \bR_{ N, n }( s )
			\,\tilde \bR_{ N, n }^T( s' ) = \tilde{ \bfrkc }\big(\, s, s' \,\big)
			\end{align*}
 			uniformly in probability and $ \sum_{ n = 1 }^N \tilde \bR_{ N, n } = 0 $.
			\end{enumerate}
\end{theorem}
\begin{proof}
	Part (a) is a simple Taylor expansion argument showing that $ H $,
	considered as a function of $ C\big(\,S,\R^P\,\big) \rightarrow C\big(\,S,\R^{P'}\,\big) $,
	is Hadamard differentiable tangential to $ C\big(\, S, \R^P \,\big) $.
	Thus \cite[Theorem 2.8]{Kosorok:2008} implies that the
	delta method is applicable, which proves (a).
	
	To prove (b) note that, by linearity of the differential,
    \begin{equation}\label{eq:samp-cov}
		\begin{split}
			N^{ -1 } \sum_{ k = 1 }^N& \tilde \bR_{ N, n }( s )\, \tilde \bR_{ N, n }^T( s' )
         		= \dH_{ \hat \btheta_{N}(s) } \left(\, N^{-1} \sum_{ k = 1 }^N \bR_{N,n}(s)\, \bR_{ N, n }^T(s')\,\right) \dH_{\hat \btheta_{N}(s')}^T\,.
		\end{split}
	\end{equation}
    The fCLT \eqref{eq:fclt} implies that $\hat \btheta_{ N } \rightarrow \btheta$ uniformly
    in probability. Additionally,
    $ \dH_{ \hat \btheta_{ N }( s ) } \rightarrow \dH_{ \btheta( s ) } $ uniformly  in probability
    as $ N \rightarrow \infty $, by the continuous mapping theorem, and so the claim follows from 
    \eqref{eq:samp-cov}.
\end{proof}

\begin{remark}
	Two observations are noteworthy. Firstly, the factors $ \sqrt{ N } $ in equation
	\eqref{eq:fclt} and $ N^{ -1 } $ in equation \eqref{eq:conv} can be replaced by
	general factors tending to infinity and zero respectively provided they are also
	changed in the subsequent equations (we stick to these factors here in order
	to keep the notation simple).
	Secondly, if $ \dH_{ \btheta( s ) } = 0 $ for all $ s \in S $, then the delta
	residuals can be identically equal to zero. In that case an assumption of higher
	differentiability of $ H $ can be used to establish a similar result using
	a second-order delta method.
\end{remark}

\section{Delta Residuals for Moment-based Statistics} \label{scn:MomentDeltaResiduals}

In this section we illustrate how Theorem \ref{thm:DeltaResiduals} can be applied to
statistics based on differentiable functions of pointwise sample moments which we call moment-based statistics.
Hereafter, unless otherwise stated, we assume that $ \{ X_n \}_{n\in \mathbb{N}} $ is a sequence of random
processes in $ C(\,S) $ such that all elements are independent and identically
distributed as $ X $. For $r \in \mathbb{N}$, the $ r $-th pointwise population moment of $ X $ is defined
by
$
 \mu^{(r)}( s ) = \E\big[ X^r( s ) \big]
$
and the (non-centered) sample moment as
\begin{equation*}
 \hat\mu^{(r)}_N( s ) = \frac{ 1 }{ N } \sum_{ n = 1 }^N X_{ n }^r( s )\,.
\end{equation*}
Statistics such as Cohen's $ d $, skewness or kurtosis can be expressed as continuously
differentiable functions of the sample moments and therefore functional delta residuals for these
statistics can be constructed from the general framework described below. Specific examples
will be discussed in Section \ref{scn:ExamplesDeltaResiduals}.

\subsection{A Functional Central Limit Theorem for Moments}

In order to apply Theorem \ref{thm:DeltaResiduals} we need to establish a fCLT for vectors
of different sample moments $ \hat\mu^{(r)}_{N} $.
We base the proof of our fCLT on the following sample path property for the process $ X $.
However, other properties allowing for fCLTs could be used. 
\begin{definition}\label{def:L4Lipshitz}
    Let $ Y $ be a process in $C(\,S) $. Given $ p \in \mathbb{ N } $, we say that $ Y $ has
    \emph{$ \mathcal{ L }^p $-H\"older continuous paths} of order $ \alpha \in (0,1] $,
    if
    \begin{equation}
    	\big\vert\, Y(s) - Y(s') \,\big\vert \leq L\, \vert\, s-s' \,\vert^\alpha
    \end{equation}
    almost surely for all $ s, s' \in S$ with $ L $ a positive random variable satisfying $ \E\big[ L^p \big] < \infty $.
\end{definition}
\begin{remark}
	$ \cL^2 $-H\"older continuous paths ensure that $ Y $ satisfies a fCLT, i.e.,
	for $ Y, Y_1, \ldots, Y_N, \ldots $ an i.i.d. sequence in $C(\,S)$, the sum $ N^{-1/2} \sum Y_n $
	converges weakly to a tight zero-mean Gaussian process which has the same covariance
	structure as $Y$, see \citet[Theorem 1]{Jain:1975}. Similar fCLTs for dependent functional
	arrays requiring a mixing condition have been recently shown in \citet[Theorem 2.1]{Dette2020functional}.	
\end{remark}
\begin{remark}
	It is obvious that $ Y $ satisfies a fCLT, eq. \eqref{eq:eq1}, if $ Y $ is a Gaussian process.
	However, it is still interesting to note that if $ Y $ has almost surely $ C^1 $-sample paths,
	then applying \citep[Theorem 4]{Landau:1970} to the components of the gradient process implies
	that $ Y $ has $ \cL^p $-H\"older continuous paths for any $ p \in \N $.
\end{remark}
The following Lemma states useful properties of processes with $ \cL^p $-H\"older
continuous paths and is an adaptation Lemma 5 of \citet{Telschow:2019}.
\begin{lemma}\label{Lem:UniformConvergence}
	Let $ \lbrace Y_{ n } \rbrace_{n \in \mathbb{N}} \sim Y $ be i.i.d.
	processes in $ C(\,S) $ and
	$ \lbrace Z_{ n } \rbrace_{n \in \mathbb{N}} \sim Z $ be i.i.d.
	processes in $ C(\,S) $ such that $ Z $ and $ Y $
	both have $ \cL^p $-H\"older continuous paths with $ p \geq 1 $ and the domain $S$ is compact.
	Assume that there exist $ s',s'' \in S $	such that $ \E\big[\,\vert Y( s' ) \vert^p\,\big] $ and
	$ \E\big[\, \vert Z( s'' ) \vert^p\,\big] $ are finite.
	Then
	\begin{enumerate}
		\item[(a)] $ \E\big[\, \Vert Y \Vert^q_\infty\,\big] < \infty $ for all $ q \leq p $.
		\item[(b)] $ \big\Vert\, \bar Y_N - \E[Y] \,\big\Vert_\infty \rightarrow 0 $ almost surely as $N$ tends to infinity.
		\item[(c)] If $ p \geq 2 $, then
			$ \big\Vert\, N^{-1} \sum_{ n = 1 }^N \left(\, Y_{ n } - \bar Y_N \,\right)
	\left(\, Z_{ n } - \bar Z_N \,\right) - \cov\!\left[\, Y,  Z\, \right] \,\big\Vert_\infty  \rightarrow 0$
			almost surely as $N$ tends to infinity. Here $\Vert \cdot \Vert_\infty $
			denotes the maximum norm on $ C(\, S\times S ) $.
	\end{enumerate}
\end{lemma}
\begin{proof}
	\emph{Proof of (a):} Using the convexity of $ \vert \cdot \vert^p $,
	$ \E\big[\,\vert Y( s' ) \vert^p\,\big] < \infty $ and $ \delta(\,s,s'\,) = \vert\, s - s' \,\vert^\alpha$  we have
    \begin{align*}
       \E\big[\,\Vert Y \Vert^p_\infty \,\big] &\leq 2^{p-1}\Bigg(\, \E\Big[\,\Vert Y-Y(s') \Vert^p_\infty\,\Big] + \E\Big[\,\vert Y(s') \vert^p\,\Big] \,\Bigg) \\
       &\leq
2^{p-1}\left(\, \E\big[ L^p\big] \max_{s\in S} \delta(\, s, s' \,) + \E\Big[\,\vert Y(s') \vert^p\,\Big]\,\right)\,,
     \end{align*}
      where $L$ is the random variable from the $ \cL^p $-H\"older property. This yields
      $ \E\big[\, \Vert Y \Vert^q_\infty\,\big] < \infty $ for all $ q \leq p $. 
      
      \emph{Proof of (b):} We apply the generic uniform convergence result in
      \citet[Theorem 21.8]{Davidson1994}.
      Since pointwise convergence holds by the strong law of large numbers, it is sufficient
      to establish strong stochastic equicontinuity of the random function $\bar Y_N - \E[Y]$.
      This is established using \citet[Theorem 21.10 (ii)]{Davidson1994}, since
		  \begin{equation*}
			\left\vert\, \bar Y_N(s)- \bar Y_N(s') - \E\big[\,Y(s)-Y(s')\,\big] \,\right\vert \leq \left(\, \sum_{ n = 1 }^N \frac{L_{ n }}{N} + \E[L]\,\right) \,\delta(\, s, s' \,) = C_N\, \delta(\, s, s' \,)
		  \end{equation*}
		  for all $s,s'\in S$. Here $L_{1},\ldots,L_{ N }\sim L$ i.i.d.
		  denote the random variables from the $\cL^p$-H\"older paths of the
		  $Y_{ n }$'s and $Y$. Hence the random variable $ C_N $ converges
		  almost surely to the constant $2 \E[L]$ by the strong law of
		  large numbers.
		  
		  \emph{Proof of (c):} First we have that, for all $s, s' \in S$,
		  \begin{equation*}
		  \begin{split}
		  	\big\vert\, Y_{ n }(s) - \bar Y_N(s)
		  			- Y_{ n }(s') + \bar Y_N(s') \,\big\vert
		  	&\leq
		  		\big\vert\, Y_{ n }(s) - Y_{ n }(s') \,\big\vert
		  			+
		  			\big\vert\, \bar Y_N(s') - \bar Y_N(s) \,\big\vert \\
		  	&\leq
		  	\left(\, L_{n}
		  		+ \frac{1}{N} \sum_{n=1}^N L_{n} \,\right)\, \delta(\, s, s' \,)\\
		  		&\leq
		  	K_{n}\, \delta(\, s, s' \,)\,,
		  \end{split}
		  \end{equation*}
		  where $ L_{n} $, $ n \in \{ 1, ..., N \} $, is the
		  random $\cL^p$-H\"older
		  constant of $Y_{ n }$.
		  Therefore $ A_n = Y_{ n } - \bar Y_N $ has $\cL^2$-H\"older paths since
		  $p\geq 2$. The same holds for $ B_n = Z_{ n } - \bar Z_N $.
		  Hence we compute
		  \begin{align*}
		   \Bigg\vert\, \frac{1}{N} \sum_{ n = 1 }^N A_n(s)&\,B_n(t) - A_n(s')\,B_n(t') \,\Bigg\vert\\
		  		&\leq \frac{1}{N} \sum_{ n = 1 }^N \Vert\, A_n\,\Vert_\infty \,\tilde K_{ n }\,\delta(\, t, t' \,) + \Vert\, B_n\,\Vert_\infty \,K_{ n } \,\delta(\, s, s' \,)  \\
		  		&\leq \sqrt{ \sum_{ n = 1 }^N \tfrac{\Vert\, A_n \,\Vert_\infty^2}{N} \sum_{ n = 1 }^N \tfrac{\tilde K_{ n }^2}{N} }\,\delta(\, t, t' \,) + \sqrt{ \sum_{ n = 1 }^N \tfrac{\Vert\, B_n \,\Vert_\infty^2}{N} \,\sum_{ n = 1 }^N \tfrac{K_{ n }^2}{N} }\, \delta(\, s, s' \,)\,,
		  \end{align*}
		  for each $ s,s',t,t' \in S $ and $\tilde K_{1},\ldots,\tilde K_{N}$ i.i.d. denote the random variables from the $\cL^2$-H\"older property of the $B_{ n }$'s. By the strong law of large numbers the random H\"older constant converges almost surely and is finite. Thus, again the generic uniform convergence result in \citet[Theorem 21.8]{Davidson1994} together with \citet[Theorem 21.10]{Davidson1994} yield the claim.
\end{proof}

To handle  vector-valued random processes, we require the following Lemma in order to prove Theorem \ref{thm:2DfCLT}. This states simple conditions for obtaining weak convergence of a vector-valued process from its components. Its generalization to arbitrary dimensional vector-valued processes is immediate.
\begin{lemma}\label{lem:2dtightness}
	Let $ X_1, X_2, \dots, X, Y_1, Y_2, \dots, Y$ be $ C(\,S) $-valued random variables on
	the probability space $ (\, \Omega, \mathcal{F}, \Prb \,) $ such that $
	X_N \rightsquigarrow X $ and $ Y_N \rightsquigarrow Y $.
	If the finite dimensional distributions of $ (\, X_N, Y_N \,) $ converge to those of $ (\, X, Y \,) $, we have
	$
	(\, X_N, Y_N \,) \rightsquigarrow (\, X, Y \,)
	$
	in $ C(\,S) \times C(\,S) $.
\end{lemma}
\begin{proof}
	Tightness of the pair $(\, X_N, Y_N \,)$ is implied by Lemma 1.4.3 and Problem 9 in Section 1.3 from \cite{Vaart1996} (hereafter VW).
	
	Moreover, the finite dimensional distributions converge and form a separating class in $ C(\,S) \times C(\,S)$ (the proof is along the lines of Example 1.3 in \cite[p. 12]{Billingsley:1999}) so in particular the joint distribution converges (arguing as in Example 5.1 in \cite[p. 57]{Billingsley:1999}).
\end{proof}  

With these preparatory results we can now prove the main theorem of this section.
\begin{theorem}\label{thm:2DfCLT}
    Let $ S $ be a compact space and $ \lbrace X_{ n } \rbrace_{n\in\mathbb{N}} $
	be i.i.d. processes in $ C(\,S) $ distributed as $X$.
    Let $ r_1, \ldots, r_K \in \N $ such that $ 1 < r_1 < \ldots < r_K $,
    for some $K \in \mathbb{N}$.
    Denote the corresponding vector as $ \br = ( r_1, ..., r_K ) $.
    Assume that $ X $ has $ \cL^{ 2r_K } $-H\"older continuous paths with $ L$ the random
    H\"older bound and there
    exists an $ s \in S $ such that $ \E\big[ X^{ 2r_K }(s) \big] < \infty $. Then
    \begin{equation}\label{eq:CLT2D}
            \sqrt{N}\,\Big(\, \hat\bmu^{ (\br) }_N 
            				- \bmu^{ (\br) } \,\Big)
             \rightsquigarrow \bG\,.
    \end{equation}
    Here $ \hat\bmu^{ (\br) }_N = \Big(\, \hat\mu^{(r_1)}_N, \ldots, \hat\mu^{(r_k)}_N \,\Big) $
    and $ \bmu^{ (\br) } = \Big(\, \mu^{(r_1)}, \ldots, \mu^{(r_k)} \,\Big) $.
    Moreover, $ \bG $ is a zero-mean Gaussian process with paths in $ C\big(\, S, \R^K \,\big) $
    and covariance matrix function $ \frkc $ having entries
    \begin{equation*}
            \frkc_{kl}(\, s, s' \,) = \E\big[\, X^{r_k}(s)\,X^{r_l}(s') \,\big] - \mu^{(r_k)}(s)\,\mu^{(r_l)}(s')~\text{ for }k,l=1,\ldots,K, \text{ and } s, s' \in S\,.
    \end{equation*} 
\end{theorem}
\begin{proof}
	First we need to establish that for all $ k = 1, \ldots, K $ the sequence
	$ \big\lbrace X_{ n }^{r_k} \big\rbrace_{n \in \mathbb{N}}$
	satisfy the CLT in $ C(\,S) $.
	To do so observe that for all $ x, y $ we have that $ \vert\, x^{r_k} - y^{r_k} \,\vert
	\leq (\, r_k - 1 \,)\, \max_{ \xi \in \{\, x,y \,\} } \xi^{ r_k - 1 } \,\vert\, x - y \,\vert $, since 
	we can factor as
	\begin{align*}
		 x^{r_k} - y^{r_k}
		 	=    (\, x-y \,)\, \sum_{ l = 0 }^{ r_k - 1 } x^{ r_k - 1 -l }\,y^{ l }
			\leq (\, x-y \,)\, (\, r_k - 1 \,)\, x^{ r_k - 1 }\,, ~ ~ ~\text{ for }  x\geq  y\\
		 y^{r_k} - x^{r_k}
		 	=    (\, y-x \,)\, \sum_{ l = 0 }^{ r_k - 1 } y^{ r_k - 1 -l }\, x^{ l }
			\leq (\, y-x \,)\, (\, r_k - 1 \,)\, y^{ r_k - 1 }\,, ~ ~ ~\text{ for }  y\geq x\,.
	\end{align*}		
	Hence
	\begin{equation}
	\begin{split}\label{eq:XkisLipschitz}
		\big\vert\, X^{r_k}(s) - X^{r_k}(s') \,\big\vert 
			&\leq (\, r_k - 1 \,)\, \max_{ s\in S } \vert\, X(s) \,\vert^{ r_k - 1 }
																\vert\, X(s) - X(s') \,\vert\\
			&\leq (\, r_k - 1 \,)\, L\, \Vert X \Vert_\infty^{ r_k - 1 } \vert\, s - s' \,\vert^\alpha\,.
	\end{split}
	\end{equation}
	Using H\"older's inequality for $ p = r_k $ and $ q = r_k \,/\, (\, r_k - 1 \,) $ we obtain for
	all $ k = 1, \ldots, K $ that
	\begin{align*}
		\E\left[\, L^2\, \Vert X \Vert_\infty^{ 2( r_k - 1 ) } \,\right]
			&\leq \E\left[ L^{ 2r_k } \right]^\frac{ 1 }{ r_k }
			      \E\left[\, \Vert X \Vert_\infty^{ 2r_k } \,\right]^\frac{ r_k - 1 }{ r_k }
			 < \infty
	\end{align*}
	by Lemma \ref{Lem:UniformConvergence}(a) applied to $ X $. As such the positive random variable
	$ M = (\, r_k - 1 \,)\, L \,\Vert X \Vert_\infty^{ r_k - 1 }(s) $ satisfies
	$ \E\big[ M^2 \big] < \infty $
	and the fCLT for each $ k = 1, \ldots, K $ follows from \citet[Theorem 1]{Jain:1975}.
	
	In order to apply Lemma \ref{lem:2dtightness} it remains to show that the
	finite dimensional distributions converge to the finite dimensional
	distributions of $ \bG $ from the statement of the theorem.
	To see this, for $s \in S$, define $ \bX^{(\br)}_{ n }(s) =
	\big(\, X^{r_1}_{ n }(s), \ldots, X^{r_K}_{ n }(s) \,\big) $ and $ \bmu^{(\br)}(s) =
	\big(\, \mu^{r_1}(s), \ldots, \mu^{r_K}(s) \,\big) $. For $ I \in \mathbb{N} $ and any
	$ s_1, \ldots, s_I \in  S $, we apply the multivariate CLT to the sequence of random vectors
	\begin{align*}
		\Big(\, \bX^{(\br)}_{ n }(s_1) - \bmu^{(\br)}(s_1), \ldots,
				 	\bX^{(\br)}_{ n }(s_I) - \bmu^{(\br)}(s_I) \,\Big)\,,
	\end{align*} 
	which yields convergence to the finite dimensional distributions of a Gaussian
	random vector with covariance matrix
	$ \Sigma \in \mathbb{R}^{K I \times K I} $ given by
	\begin{align*}
		\Sigma_{ i \cdot K + k, j \cdot K + l }  &= \cov\big[\, X^{r_k}(s_i) - \mu^{(r_k)}(s_i),\, X^{r_l}(s_j) - \mu^{(r_l)}(s_j) \,\big]\\
						   &=  \E\big[\, X^{r_k}(s_i)\,X^{r_l}(s_j) \,\big] - \mu^{(r_k)}(s_i)\,\mu^{(r_l)}(s_j)\,,
	\end{align*}
	for $ i,j = 1,\ldots, I-1 $ and $ k,l = 1,\ldots, I $.
	Hence these finite dimensional distributions converge to those of $ \bG $, which finishes the proof.
\end{proof}
\begin{remark}\label{rmk:momentsLipshitz}
	During the proof of the previous theorem we showed in eq. \eqref{eq:XkisLipschitz} that
	if $ X $ has $ \cL^{ 2r } $-H\"older continuous paths and there
	exist an $ s \in S $ such that $ \E\big[ X^{ 2r }(s) \big] < \infty $
	then $ X^r $ has $ \cL^{2} $-H\"older continuous paths.
\end{remark}

\subsection{ Delta Residuals }\label{scn:delta_res_cohensd}

In the previous section we established a multivariate functional CLT for sample moments.
In order to apply Theorem 1 we must construct untransformed residuals having the covariance structure given in Theorem \ref{thm:2DfCLT}.
\begin{proposition}[Moment Residuals]\label{prop:MomentResiduals}
	We call the processes $ R_{ N, n }^{(r)} = X_{ n }^r - \hat\mu^{(r)}_N $ the $ r $-th
	moment residuals, which satisfy $ \sum_{ n = 1 }^N R_{ N, n }^{(r)} = 0 $.
	Let $ r_1, r_2 \in \N $ and assume that $ X $ has
	$ \cL^{ 2 \max(\, r_1, r_2 \,) } $-H\"older continuous paths, then
	\begin{equation}
		\lim_{ N \rightarrow \infty } N^{-1} \sum_{ n = 1 }^N
													R_{ N, n }^{(r_1)}(s)\,R_{ N, n }^{(r_2)}(s') 
		= \E\left[\, X^{r_1}(s)\, X^{r_2}(s') \,\right] - \mu^{(r_1)}(s)\,\mu^{(r_2)}(s')\,,~ ~ s,s'\in S
	\end{equation}
	almost surely uniformly in $C(\, S \times S )$.
\end{proposition} 
\begin{proof}
	The $r$-th moment residuals sum to zero by construction.
	By Lemma \ref{Lem:UniformConvergence}(c) and Remark \ref{rmk:momentsLipshitz} it follows that,
	if $ s,s'\in S $,
	\begin{align*}
		N^{-1}\sum_{ n = 1 }^N R_{ N, n }^{(r_1)}(s)\, R_{ N, n }^{(r_2)}(s') 
			&= N^{-1}\sum_{ n = 1 }^N \big(\, X_{ n }^{r_1}(s) - \hat\mu^{(r_1)}_N(s) \,\big)
									   \,\big(\, X_{ n }^{r_2}(s') - \hat\mu^{(r_2)}_N(s') \,\big)\\
			&\xrightarrow{ N \rightarrow \infty } \cov\big[\,X^{r_1}(s), X^{r_2}(s')\,\big]\\
			&= \E\left[\, X^{r_1}(s)\,X^{r_2}(s') \,\right] - \mu^{(r_1)}(s)\,\mu^{(r_2)}(s'))
	\end{align*}
	uniformly almost surely in $ C(\, S \times S ) $.
\end{proof}
Combining the fCLT for vectors of moments, Theorem \ref{thm:2DfCLT},
and the above Proposition \ref{prop:MomentResiduals}, which shows that moment residuals can play the
role of the untransformed residuals in Theorem \ref{thm:DeltaResiduals}, we obtain the following
corollary.
\begin{corollary}\label{cor:momentfCLT}
	Assume the assumptions of Theorem \ref{thm:2DfCLT}
	and assume that $ H \in C^1\big(\, \R^K \big) $. Then
	\begin{enumerate}
		\item[(a)] Let $ \tilde G $ be a zero-mean Gaussian
				process with covariance
				\begin{equation*}
					\tilde{ \frkc }(\, s, s' \,)
						= \dH_{ \bmu^{ (\br) }(s) }\, \frkc(\, s, s' \,)\,
			  				\dH_{ \bmu^{ (\br) }(s') }^T\,,
				\end{equation*}
				where $ \frkc(\, s, s' \,) $ is given in Theorem \ref{thm:2DfCLT}.
				Then
		\begin{equation}\label{eq:MomentfCLT}
		\sqrt{ N } \Bigg(\, H\Big( \hat\bmu^{ (\br) }_N \Big)
							- H\Big( \bmu^{ (\br) } \Big) \,\Bigg)
				\rightsquigarrow \tilde G\,,~ N \rightarrow \infty\,.
				\end{equation}
	\item[(b)] The delta residuals
	\begin{equation}\label{eq:DeltaMomentResiduals}
		\tilde R_{ N, n }^{ ( \br ) } 
			= \dH_{ \hat\bmu^{ (\br) }_N }\, \bR^{ ( \br ) }_{ N, n }  \,,
		~ n=1,\ldots,N\,,
	\end{equation}
 	where $ \bR^{ ( \br ) }_{ N, n } =
 	\Big(\, R^{ (r_1) }_{ N, n },\ldots, R^{ (r_K) }_{ N, n } \,\Big)^T $, satisfy
	\begin{align}\label{eq:covrates_estim}
		\lim_{ N \rightarrow \infty } N^{-1}\sum_{ n = 1 }^N
						\tilde R_{ N, n }(s)\, \tilde R_{ N, n }^T(s') = \tilde{\frkc}(\, s, s' \,)
 	\end{align}
 	uniformly almost surely in $C(\, S \times S )$.
	\end{enumerate}
\end{corollary}

\subsection{Examples of Delta Residuals Based on Moments}\label{scn:ExamplesDeltaResiduals}

We will now discuss a series of examples in which our results can be applied.
Throughout this section we will assume that we have a sequence
$ \left\{ X_{n} \right\}_{n \in \N} $ satisfying the assumptions from Theorem
\ref{thm:2DfCLT}.\vspace{0.2cm}
 
\textbf{Sample variance.\,} A first simple example of delta residuals can be constructed for the sample variance
\begin{equation}
	\hat\sigma^2_N(s)		= \frac{ 1 }{ N } \sum_{ n = 1 }^N \Big(\, X_{ n }(s) - \hat\mu^{(1)}_N(s) \,\Big)^2
		= \hat\mu^{(2)}_N(s) - \Big(\, \hat\mu^{(1)}_N(s) \,\Big)^2\,,
\end{equation}
which is a uniform almost surely consistent estimator for the pointwise population variance
$ s \mapsto \var[X(s)] $.
The transformation of $ (\, \mu^{(1)}, \mu^{(2)} \,) $ is given by $ H(\, x, y \,) = y - x^2 $ and the
resulting delta residuals are
\begin{align*}
	R_{ N, n }^{\hat\sigma^2_N}
		&= \left(\, \frac{\partial H}{\partial x}\Big(\, \hat\mu^{(1)}_N, \hat\mu^{(2)}_N \,\Big),\, \frac{\partial H}{\partial y}\Big(\, \hat\mu^{(1)}_N, \hat\mu^{(2)}_N \,\Big) \,\right) 
		\,\Big(\, X_{ n } -\hat\mu^{(1)}_N, X_{ n }^2 -\hat\mu^{(2)}_N \,\Big)^T\\
		&= \big(\, -2\,\hat\mu^{(1)}_N, 1 \,\big)\,\Big(\, X_{ n } -\hat\mu^{(1)}_N,\, X_{ n }^2 -\hat\mu^{(2)}_N \,\Big)^T\\
		&= \left(\, X_{ n } - \hat\mu^{(1)}_N \,\right)^2 - \hat\sigma^2_N\,.
\end{align*}\vspace{0.2cm}

\textbf{Cohen's $d$.\,} Recently, effect size measures gained popularity in the analysis of fMRI data \citep{Bowring:2020,Vandekar2021improving}. \citet{Bowring:2020} used the pointwise Cohen's $ d $ statistic defined by
\begin{equation*}
	\hat d_N(s) = \frac{\hat\mu^{(1)}_N(s)}{\sqrt{\, \hat\mu^{(2)}_N(s) - \left(\, \hat\mu^{(1)}_N(s) \,\right)^2 }\,}\,,~ s \in S
\end{equation*}
which is a uniform almost surely consistent estimator for the pointwise population Cohen's $d$
\begin{equation*}
	d(s) = \frac{\mu^{(1)}(s) }{ \sqrt{ \,\mu^{(2)}(s) - \Big(\, \mu^{(1)}(s) \,\Big)^2 \,} }\,,~ s \in S\,.
\end{equation*}
Note that the denominator will be non-zero, with probability $1$ for all $s\in S$, if
$ N \geq D+1 $ \citep[Lemma 11.2.10]{Adler:2009}. The residuals for Cohen's $d$ can be derived
from the transformation $H(\,x,y\,) = x\, /\, \sqrt{\, y-x^2 \,}$, i.e.,
\begin{align*}
	R_{ N, n }^{ \hat d_N }
		&= \left(\,  \frac{ \hat\mu^{(2)}_N }{ \left(\, \hat\mu^{(2)}_N - \left(\, \hat\mu^{(1)}_N \,\right)^2 \,\right)^{3/2} },
				  -\frac{ \hat\mu^{(1)}_N }{ 2\, \left(\, \hat\mu^{(2)}_N - \left(\, \hat\mu^{(1)}_N \,\right)^2 \,\right)^{3/2} } \,\right) \Big(\, X_{ n } - \hat\mu^{(1)}_N, X_{ n }^2 -\hat\mu^{(2)}_N \,\Big)^T\\
		&= \frac{ X_{ n } - \hat\mu^{(1)}_N }{ \hat\sigma_N } - \frac{ \hat d_N }{ 2\,\hat\sigma_N^2 } \left(\, \left(\, X_{ n } - \hat\mu^{(1)}_N \,\right)^2 - \hat\sigma_N^2 \,\right)\,.
\end{align*}
Deriving this identity is a little tedious.
An elegant shortcut can be taken using the observation that we can treat the
residuals $ R_n^{\hat\sigma^2_N } $ as the untransformed residuals, and use the equivalent definition for
Cohen's $ d $ of $ \hat d_N = \hat\mu^{(1)}_N\, / \sqrt{\, \hat\sigma^2_N \,}$.
Since Theorem \ref{thm:2DfCLT} together with the delta residuals method applied to the
continuously differentiable function $F(\,x, y \,) = (\,x,y-x^2\,)$ yield a fCLT for the vector
valued process $\big(\, \hat\mu^{(1)}_N, \hat\sigma^2_N \,\big)$, we can use the transformation
$\tilde H(\,x, y \,)=x\,/\sqrt{\,y\,}$ to obtain
\begin{align*}
	R_{ N, n }^{\hat d_N}
		&= \left(\,  \frac{ 1 }{ \hat\sigma_N }, -\frac{\hat\mu^{(1)}_N }{ 2\, \hat\sigma^{3}_N } \,\right)\, \Big(\, X_{ n } - \hat\mu^{(1)}_N, \Bigg(\, X_{ n } - \hat\mu^{(1)}_N \,\Big)^2 - \hat\sigma^2_N \,\Bigg)^T\\
		&= \frac{ X_{ n } - \hat\mu^{(1)}_N }{ \hat\sigma_N } - \frac{ \hat d_N }{ 2\,\hat\sigma_N^2 } \left(\, \left(\, X_{ n } - \hat\mu^{(1)}_N \,\right)^2 - \hat\sigma_N^2 \,\right)\,.
\end{align*}
This can also be understood as an application of the chain rule, i.e., $ dH_{(\,x, y \,)} =d( \tilde H \circ F ){(\,x, y \,)} = d\tilde H_{F(\,x, y \,)}dF_{(\,x, y \,)} $. As an illustration of Corollary 1, we derive the asymptotic covariance structure for Cohen's $ d $.
\begin{corollary}\label{cor:SNRfCLT}
 Under the assumptions of Theorem \ref{thm:2DfCLT} we have that
 $
 \sqrt{N} \Big(\, \hat d_{N}(s) - d(s) \,\Big) \rightsquigarrow \tilde G
$
with covariance structure $\tilde{\frkc}$ given by
\begin{equation}\label{eq:CohendCov}
 \tilde{\frkc}(\, s, s' \,) = \left(\, \sigma(s)^{-1}, -\mu(s)\,\sigma(s)^{-\tfrac{3}{2}} \,\right)\,\frkc(\, s, s' \,)\,
\left(\, \sigma(s')^{-1}, -\mu(s')\,\sigma(s')^{-\tfrac{3}{2}}\,\right)^T\,.
\end{equation}
Moreover, if $X$ is a Gaussian process, $\tilde{\frkc}$ simplifies to
\begin{equation}\label{eq:gaussSNRcorrelation}
 \tilde{\frkc}(\, s, s' \,) = \frac{\frkc_{11}(\, s, s' \,)}{\sigma(s)\,\sigma(s')} + \frkc_{11}(\, s, s' \,)^2 \frac{\mu(s)\,\mu(s')}{2\,\sigma^3(s)\,\sigma^3(s')}\,.
\end{equation}
\end{corollary}
\begin{proof}
	The general form of the covariance structure, i.e., \eqref{eq:CohendCov}, follows by applying Corollary
	\ref{cor:SNRfCLT}.
	 When $X$ is a Gaussian process the asymptotic covariance structure simplifies significantly.
	 To show this, we define $\varepsilon(s) = X(s)-\mu(s)$ and use the fact from the moments
	 of multivariate normal distributions, better known as Isserlis' theorem, cf. Theorem 1 in
	 \citet{Vignat:2012},
 \begin{equation*}
    \E\Big[\,\varepsilon^2(s)\,\varepsilon^2(s')\,\Big] = \sigma^2(s)\,\sigma^2(s')+2\,\frkc_{11}(\, s, s' \,)^2\,,
 \end{equation*}
for all $s,s'\in  S$ to compute
 \begin{align*}
    \frkc_{22}(\, s, s' \,) &= \E\Big[\, \Big(\,\varepsilon(s)^2 - \sigma^2(s)\,\Big) \,\Big(\,\varepsilon(s')^2 - \sigma^2(s')\,\Big) \,\Big] \\
                            &= \E\Big[\, \varepsilon^2(s)\, \varepsilon^2(s') \,\Big] -  \E\Big[\, \varepsilon^2(s')\,\Big]\,\sigma(s)^2 -  \E\Big[\,\varepsilon^2(s)\,\Big]\, \sigma^2(s') + \sigma^2(s)\,\sigma^2(s')\\
                            &= 2\frkc_{11}(\, s, s' \,)^2\,.
 \end{align*}
 Finally, we note that
 \begin{equation*}
  \frkc_{12}(\, s, s' \,) = \E\Big[\, \varepsilon(s)\,\Big( \varepsilon^2(s')-\sigma^2(s') \,\Big) \,\Big]= \E\Big[\, \varepsilon(s)\,\varepsilon^2(s')\,\Big] = 0 =\frkc_{21}(\, s, s' \,)\,,
 \end{equation*}
  yielding the simplified version of the limiting covariance structure.
\end{proof}

\textbf{Skewness and excess kurtosis.\,} There exist several measures for the skewness and kurtosis of a
 distribution: a broad overview of these can be found in \citet{Joanes:1998}.
 Here we will use some of the most standard measures which date back to \citet{Fisher:1930}. The skewness estimator is given by
\begin{equation}\label{eq:g1}
	g_1^N = \frac{ N^{-1}\sum_{ n = 1 }^N \left(\, X_{ n } - \hat\mu^{(1)}_N \,\right)^3 }{ \hat\sigma^{3}_N  } = \frac{ \hat\mu^{(3)}_N -3\,\hat\mu^{(1)}_N\, \hat\mu^{(2)}_N + 2\,\left(\, \hat\mu^{(1)}_N \,\right)^3 }{ \left(\, \hat\mu^{(2)}_N - \left(\, \hat\mu^{(1)}_N \,\right)^2 \,\right)^{3/2}  }\,.
\end{equation}
Hence the transformation $H(\,x,y,z\,) = (\,z - 3\,x\,y\, + 2\,x^3)\,/\,(\,y-x^2\,)^{3/2}$ together with the first three moment residuals yield functional delta residuals with the correct covariance structure.
The sample excess kurtosis can be defined as
\begin{equation}\label{eq:g2}
	g_2^N = \frac{ N^{-1}\sum_{ n = 1 }^N \left(\, X_{ n } - \hat\mu^{(1)}_N \,\right)^4 }{ \hat\sigma^{4}_N  } - 3
				 = \frac{ \hat\mu^{(4)}_N - 4\,\hat\mu^{(1)}_N\,\hat\mu^{(3)}_N
				 	+ \left(\, \hat\mu^{(1)}_N \,\right)^2\hat\mu^{(2)}_N
				 	- 3\,\left(\, \hat\mu^{(1)}_N \,\right)^4  }{ \left(\, \hat\mu^{(2)}_N
				 	- \left(\, \hat\mu^{(1)}_N \,\right)^2 \,\right)^2  } - 3\,,
\end{equation}
which shows that the transformation $H(\,x,y,z,w\,) = (\,w - 4\,x\,z + x^2\,y - 3\,x^4)\,/\,(\,y-x^2\,)^{2}$ together with the first three moment residuals yield functional delta residuals with the correct covariance structure.

\subsection{A Multiplier Bootstrap Functional Limit Theorem}\label{sec:MultfCLT}

We have shown that functional delta residuals can be used to
approximate the covariance structure of the limiting
process $ \tilde G $ given in Corollary \ref{cor:momentfCLT}.
As such, as we will prove formally in this section, the multiplier
bootstrap using the delta
residuals can be used to approximate sample path properties of the limiting process.
This will, importantly, enable us to estimate the quantiles of the maximum of $ \tilde G $.
To establish this we prove both weak convergence and the stronger conditional weak
convergence VW Chapter 2.6 of the multiplier bootstrap process to $ \tilde G $.
The main results used in the proof are the Jain-Marcus theorem \citep[Theorem 1]{Chang2009} and
Theorem 2 from \citep{Chang2009}.

In the following we assume that $ \left\lbrace\, g_{ N, 1 }, \ldots, g_{ N, N } :
N \in \mathbb{N},~1 \leq n \leq N  \,\right\rbrace $ is a triangular array of i.i.d.
random variables defined on the probability space $ (\, \Omega_g, \mathfrak{P}_g, \Prb_g \,) $
satisfying $ \E[\, g_{ N, 1 } \,] = 0 $ and $ \E[\, g_{ N, 1 }^2 \,] = 1 $.
Since the random variables in the sequence $ \lbrace X_{ n } \rbrace_{
n \in \mathbb{N}} $ are defined on the probability space $ (\, \Omega, \mathfrak{P}, \Prb \,) $,
we define extensions of these random variables to the product space
$ (\, \Omega \times \Omega_g, \mathfrak{P} \times \mathfrak{P}_g, \Prb \otimes \Prb_g \,) $
by defining
$ g_{ N, n }(\, \omega_g, \omega \,) = g_{ N, n }( \omega_g ) $ and $ X_{ n }(\, \omega_g, \omega \,) = X_{ n }( \omega )
$ for all $(\,\omega_g, \omega\,) \in \Omega \times \Omega_g$.

For $ N \geq 1 $, the multiplier bootstrap process of the functional delta residuals is defined on the product probability space by
\begin{equation}
  	\tilde G^{(g, \br)}_N(\, \omega_g, \omega \,) = \frac{ 1 }{ \sqrt{N} } \sum_{ n = 1 }^N g_{ N, n }(\omega_g) \tilde R_{ N, n }^{(\br)}( \omega )\,,~ ~ ~(\,\omega_g, \omega\,) \in \Omega \times \Omega_g\,.
\end{equation}
The multipliers $ g_{ N, n } $ and the $ X_{ n } $'s are assumed to be independent on the product space.
In particular this means that the $ g_{ N, n }$ are independent of the functional delta residuals $ \tilde R^{ ( \br ) }_{ N, n } $
defined in equation \eqref{eq:DeltaMomentResiduals}.
To shorten the notation, given a random variable $ Y $ on $ (\, \Omega \times \Omega_g, \mathfrak{P} \times \mathfrak{P}_g, \Prb \otimes \Prb_g \,) $ we define the random variable $ Y_\omega = Y(\, \cdot, \omega \,) $ on
$ (\, \Omega, \mathfrak{P}, \Prb \,) $ and, in a slight abuse of notation, will also write $ X_\omega = X(\omega) $ for a random variable $X$ on $ (\, \Omega, \mathfrak{P}, \Prb \,) $.
  	
\begin{theorem}\label{thm:DeltaResidualsBoots}
 	Under the assumptions of Theorem \ref{thm:2DfCLT} the following statements hold
	\begin{equation*}
  	(i)~~	\tilde G^{(g, \br)}_N \rightsquigarrow \tilde G~~~~~~~~~~~~~~~~~~~~
  	(ii)~~ \sup_{ h \in \cB } \Big\vert\, \E_g\Big[\, h\Big( \tilde G^{(g, \br)}_{ N } \Big)\,\Big] - \E\Big[\, h\big( \tilde G \big)\,\Big] \,\Big\vert
  						 \rightarrow 0\,.
	\end{equation*}
   	Here convergence in $(ii)$ is in probability with respect to $ (\, \Omega, \mathfrak{P},\Prb \,) $, $\E_g$ is the expectation with respect to $ (\, \Omega_g, \mathfrak{P}_g,\Prb_g \,) $ and $ \cB $ is
   	the set of all $ h: C(\,S) \rightarrow \R $ such that
   	$ \sup_{ f \in C(\,S) } \vert\, h ( f ) \,\vert \leq 1 $
   	and $ \vert\, h ( f ) - h ( f' )\,\vert \leq \Vert\, f - f' \,\Vert_\infty $ for all $ s, s' \in S $.
\end{theorem}

\begin{proof}
 	In order to prove $(i)$ we use the decomposition
 	\begin{equation}
 	\begin{split}
 	\label{eq:Decomp}
		\tilde G_N^{(g, \br)}
		&= \frac{ 1 }{ \sqrt{N} } \sum_{ n = 1 }^N g_{ N, n } \,\dH_{ \hat\bmu^{ (\br) }_N }
		\left(\, \bX_{ n}^{(\br)} - \hat\bmu^{ (\br) }_N \,\right)\\
		&= \boldsymbol{ B }_N\,\bX^{(g, \br)}_N + \dH_{ \bmu^{ (\br) } }\,\bX^{(g, \br)}_N
					  + \sqrt{ N }\, \bar g_{ N } \,\boldsymbol{ C }_N\,,
	\end{split}
	\end{equation}
	where 
	$ \boldsymbol{ B }_N = \dH_{ \hat\bmu^{ (\br) }_N } - \dH_{ \bmu^{ (\br) } } $,
	$\bX^{(g, \br)}_N = \sum_{ n = 1 }^N \frac{ g_{ N, n } }{ \sqrt{N} }
			\left(\, \bX_{ n }^{(\br)} - \bmu^{ (\br) } \,\right)$ and
	$\boldsymbol{ C }_N = \dH_{ \hat\bmu^{ (\br) }_N }
	\left(\, \bmu^{ (\br) } - \hat\bmu^{ (\br) }_N \,\right). $
	We first establish that $\bX^{(g, \br)}_N$ converges weakly to $\bG$ from Theorem \ref{thm:2DfCLT} by using Lemma
	\ref{lem:2dtightness}. To do so we demonstrate convergence of the component processes
	and of the finite dimensional distributions.
	For weak convergence of the components, it is sufficient to verify conditions
	$ (A) $, $ (B) $, $ (C) $ and $ (D) $  from
	\citet{Chang2009} for the processes
 	$ Z_{ N, n } =  X_{ n }^{r} \,/\, \sqrt{N\,} $, $r\in \{r_1,...,r_K\}$, as the result then follows by applying their Lemma 1 and our Theorem \ref{thm:2DfCLT}. $ (A) $
 	and $ (C) $ hold (see the proof of Theorem  \ref{thm:2DfCLT} and the Remark thereafter.
	For each $r$, using the convention
 	from \citet{Chang2009} that $\{ \mathcal{A} \}$ denotes the indicator function of a set
 	$ \mathcal{A} $, we have for every $ \eta > 0 $ that
	\begin{equation*}
 	\begin{split}
 		\sum_{ n = 1 }^N \E\left[\, \Vert Z_{ N, n } \Vert^2_\infty
 									\,\big\{ \Vert Z_{ N, n } \Vert_\infty > \eta \,\big\} \,\right]
 						&= \E\left[\, \left\Vert\, X^{r} \,\right\Vert_\infty^2
 		 							 \left\{\, \left\Vert X^{r} \,\right\Vert_\infty >
 		 							 \sqrt{N} \eta \,\right\} \,\right]
 						\xrightarrow{ N\rightarrow\infty } 0\,.
 	\end{split}
 	\end{equation*}
 	This follows by the Dominated Convergence Theorem since by assumption on $X$ we have
 	$ \E\big[\,\Vert X^{r} \Vert_\infty^2\,\big] < \infty $ by Lemma \ref{Lem:UniformConvergence}(a).
 	Thus $ ( B ) $ holds. Condition $ (D) $ holds since
	$
 		\sum_{ n = 1 }^N \E \big[\,\Vert Z_{ N, n } \Vert^2_\infty\,\big]
 						= \E\big[\,\Vert X^{r} \Vert_\infty^2\,\big]
 						< \infty\,.
	$
  	(Note that a square is missing in condition $ (D) $ in \citet{Chang2009} as can be seen
 	by following the proof of Lemma $ 1 $ in \citet{Chang2009}.)
 	This shows that for each $k \in \{ 1, \ldots, K \} $ the process
 	$N^{-1/2}\sum X_{n}^{r_k} - \mu^{(r_k)} $ converges weakly in the
 	space $\ell^\infty(S)$ of bounded functions over $ S $ to the process
 	given by the $k$-th component of $ \bG $.
 	Since by Theorem \ref{thm:2DfCLT} the components of $ \bG $ are
 	$ C(\,S) $-valued and the sample paths of all $Z_{ N, n }$'s are
 	also $C(\,S)$-valued,
 	VW's Lemma 1.3.10 establishes weak convergence in $C(\,S)$ of
 	the component processes.  
 	Convergence of the finite dimensional distributions follows directly from the
 	independence of the multipliers and the $X_{n}$'s by the multivariate CLT as
 	in Theorem \ref{thm:2DfCLT}.
 	Therefore $ \bX^{(g, \br)}_N \rightsquigarrow \bG$ in $ C\big(\,S, \R^K \big) $ by Lemma
 	\ref{lem:2dtightness}. 	
 	An application of the
 	Continuous Mapping Theorem implies that
 	$ \dH_{ \bmu^{ (\br) } } \bX^{(g, \br)}_N \rightsquigarrow \dH_{ \bmu^{ (\br) } } \bG = \tilde G $.
 	Moreover, $ \boldsymbol{ B }_N $ converges uniformly almost surely to zero by Lemma
 	\ref{Lem:UniformConvergence} and the Continuous Mapping Theorem.
 	Hence Slutsky's Lemma (VW Example 1.4.7) yields that $ \boldsymbol{ B }_N \bX^{(g, \br)}_N $
 	converges weakly to zero.
 	The same holds true for $  \sqrt{ N }\, \bar g_{ N }\, \boldsymbol{ C }_N $. Combining these observations
 	with the decomposition from eq. \eqref{eq:Decomp} it follows that
 	$ \tilde G_N^{(g, \br)} \rightsquigarrow \tilde G $ in $C(\, S )$.

	We turn to the proof of $(ii)$. Given $ \omega \in \Omega $, consider the decomposition
	\begin{equation}
	\begin{split}\label{eq:splittingExpec}
 		\sup_{ h \in \cB } \Big\vert\, \E_g\Big[\, &h\Big(\, \tilde G^{(g, \br)}_{N\omega} \,\Big) \,\Big] - \E\Big[\, h\big( \tilde G \big)\,\Big]\, \Big\vert\\
 		&\leq \sup_{ h \in  \cB } \Big\vert  \E_g\Big[\, h\Big(\, \tilde G^{(g, \br)}_{N\omega} \,\Big) \,\Big] - \E_g\Big[\, h\Big(\dH_{ \hat\bmu^{ (\br) }_{ N \omega } }\bX^{(g, \br)}_{N\omega} \Big)\,\Big] \,\Big\vert\\
 		&~~~~~~~~~~~~~~~~
 		+ \sup_{h \in \cB } \Big\vert\, \E_g\Big[\, h\Big( \dH_{ \hat{\bmu}^{ (\br) }_{ N \omega } }\bX^{(g, \br)}_{N\omega} \Big)\,\Big] - \E\Big[\, h( \tilde G )\,\Big] \,\Big\vert   \\
 	&\leq \E_g\Big[\, \big\Vert\, \tilde G^{(g, \br)}_{N\omega} - \dH_{ \hat\bmu^{ (\br) }_{ N \omega } } \bX^{(g, \br)}_{N\omega} \,\big\Vert_\infty\,\Big]
 	\\&~~~~~~~~~~~~~~~~
 		  + \sup_{ h \in \cB } \Big\vert\, \E_g\Big[\, h\Big( \dH_{ \hat{\bmu}^{ (\br) }_{ N \omega } } \bX^{(g, \br)}_{N\omega} \Big)\,\Big] - \E\Big[\, h(\tilde G)\,\Big] \,\Big\vert\\
 	&= c\, \Vert\, \boldsymbol{ C }_{N\omega} \,\Vert_\infty
 		  + \sup_{ h \in \cB } \Big\vert\, \E_g\Big[\, h\Big( \dH_{ \hat{\bmu}^{ (\br) }_{ N \omega } } \bX^{(g, \br)}_{N\omega} \Big)\,\Big] - \E\Big[\, h(\tilde G)\,\Big] \,\Big\vert
 	\end{split}
 	\end{equation}
 	Here $\boldsymbol{ C }_{N}$ is defined in the decomposition \eqref{eq:Decomp} and $ c = \E_g\big[\, \vert \sqrt{N\,}\, \bar g_{N} \vert \,\big] < \infty $, which follows from
 	$ g_{ N, n }$'s being i.i.d. with unit variance. Moreover, since
 	$ \Vert \boldsymbol{ C }_{N\omega} \Vert_\infty \leq \big\Vert  \dH_{ \hat\bmu^{ (\br) }_{ N \omega }} \big\Vert_\infty
 	 \Vert\, \hat\bmu^{ (\br) }_{N\omega} - \bmu^{ (\br) } \,\Vert_\infty $ it follows from Lemma
 	 \ref{Lem:UniformConvergence}(b) and the Continuous Mapping Theorem that the first term converges to zero as $N$ tends to infinity for almost all $ \omega \in \Omega $.
 	
 	It remains to show that the second term converges to zero in probability.
 	The proof of the previous part also established that Theorem $2$ of \citet{Chang2009}
 	is applicable to $ \bX^{(g, \br)}_{N\omega} $.
 	Therefore, for any subsequence $ \bX^{(g, \br)}_{N'\omega} $, we can choose a
 	subsubsequence $ \bX^{(g, \br)}_{N''\omega} $ (VW Lemma 1.9.2(ii))
 	such that for almost all $\omega\in\Omega$, 
 	\begin{equation}\label{eq:blabla}
 		\sup_{ h \in \cB } \Big\vert\, \E_g\Big[\, h\Big(\bX^{(g, \br)}_{N''\omega} \Big)\,\Big] - \E\Big[\, h(\bG) \,\Big] \,\Big\vert \xrightarrow{N''\rightarrow\infty} 0. 	
 	\end{equation}
	Since $ \hat\mu_{N} $ converges to $ \mu $ $\Prb$-almost surely, there exists $\Omega'\subset\Omega$
 	with $\Prb(\Omega')=1$ such that for all $\omega\in\Omega'$ it holds that $ \hat\mu_{N''\omega}
	 \rightarrow \mu $ as $ N'' $ tends to infinity and \eqref{eq:blabla} holds.
 	VW's Theorem 1.12.2 implies that $ \bX^{(g, \br)}_{N'' \omega } $ converges in distribution
 	for all $ \omega \in \Omega' $. Hence Slutsky's Lemma implies that $\dH_{ \hat\bmu^{ (\br) }_{N''\omega} }\!\bX^{(g, \br)}_{N''\omega} \rightsquigarrow \tilde G$ for all $ \omega \in \Omega' $.
 	A further application of 	VW's Theorem 1.12.2 yields
 	$$
 	\sup_{ h \in \cB } \Big\vert\, \E_g\Big[\, h\Big(\dH_{ \hat\bmu^{ (\br) }_{N''\omega}} \bX^{(g, \br)}_{N''\omega} \Big) \,\Big] - \E\Big[\,h(\tilde G)\,\Big] \,\Big\vert \rightarrow 0\,
 	 $$
 	 for all $ \omega \in \Omega' $.
 	 Since the subsequence $N'$ was arbitrary the claim follows from eq. \eqref{eq:splittingExpec}
 	 and VW's Lemma 1.9.2(ii).
\end{proof}
The usefulness of the above theorem is mainly due to the following corollary.
\begin{corollary}\label{cor:mbootconsistency}
Given any continuous function $ F: C(\, S ) \rightarrow \R $, for every point $ a \in \R $
at which $ \Prb\,\big(\, F( \tilde G ) \leq a \,\big) $ is continuous, we have that
$ \Prb\,\Big(\, F\Big( \tilde G_{N\omega}^{(g, \br)} \Big) \leq a \,\Big) \rightarrow \Prb\,\Big(\, F( \tilde G )\leq a \,\Big)$ for almost all $ \omega \in \Omega$.
\end{corollary}
\begin{proof}
Suppose the claim is false, then there exists $ a \in \R $ at which $ \Prb\big(\, F( \tilde G ) \leq a \,\big) $ is continuous, $ \epsilon > 0 $, $\tilde \Omega \subset \Omega $ with $\Prb( \tilde \Omega ) > 0 $ and a subsequence
$ ( N_{j} )_{j \in \N } $ such
that for all $j$, $ \big\vert\, \Prb\,\Big(\, F\big(\,\tilde G_{N_{j}\omega}^{(g, \br)} \big) \leq a \,\Big) - \Prb\,\big( \,F\big( \tilde G \big) \leq a \,\big) \,\big\vert > \epsilon$ for all $\omega \in \tilde \Omega$.
Now, applying Theorem \ref{thm:DeltaResidualsBoots}(ii) and Lemma 1.9.2 (ii) from VW, it follows that there exists a subsubsequence $(N_{j_k})_{k\in \mathbb{N}}$ such that $ \sup_{ h \in \cB } \big\vert\, \E_g\big[\, h\big( \tilde G^{(g, \br)}_{N_{j_k}\omega} \big)\,\big] - \E_g\big[\, h( \tilde G )\,\big] \,\big\vert $ converges to $0$ for almost all $\omega \in \Omega$. In particular by Theorem 1.12.2 in VW and the Continuous Mapping Theorem, $ F\Big(\tilde G_{N_{j_k}\omega}^{(g, \br)} \Big)$ converges weakly to $F(\tilde G)$. Hence $\Prb( \tilde \Omega ) = 0 $. This gives a contradiction.
\end{proof}

\begin{remark}\label{remark:max}
 The above corollary applies when $F$ is the maximum norm $\Vert \cdot \Vert_\infty$. This means that the multiplier bootstrap consistently estimates the quantiles of the maximum, which we will use for the construction of simultaneous confidence bands in the next section.
\end{remark}

\begin{remark}\label{remark:NdependExtension}
	The construction of functional delta residuals and the above multiplier theorem can be extended to converging
	sequences of transformations $ \big(H_N\big)_{N\in \mathbb{N}} \subset C^1\big( \mathbb{R}^D, \mathbb{R} \big) $ to a
	transformation $H \in C^1\big( \mathbb{R}^D, \mathbb{R} \big)$.
	This requires that the gradients $\nabla H\vert_x$, $\nabla H_N\vert_x$ satisfy the following uniform convergence for some $\epsilon >0$:
	\begin{equation}
	\begin{split}
		&\lim_{N\rightarrow\infty} \max_{x \in D^{\epsilon}_{\bmu^{ (\br) }}} \big\vert\, \nabla H_N\vert_x - \nabla H\vert_x \,\big\vert = 0\,\\
		\text{ with }~ ~ ~ ~		&D^{\epsilon}_{\bmu^{ (\br) }} = \Big\{\, x\in \mathbb{R}^D:~x = \bmu^{ (\br) }(s) \pm \eta ~\text{ for some }s\in S\,,~\vert \eta\vert \leq \epsilon \,\Big\}\,.
	\end{split}
	\end{equation}	
	The proof that the functional delta method remains valid for such transformations can be found in \ref{app:auxlemmas}. Moreover, the only change in the proof of Theorem \ref{thm:DeltaResidualsBoots} is that the extended continuous mapping theorem
	\cite[Theorem 7.24]{Kosorok:2008} needs to be used.
\end{remark}

\subsection{Simultaneous Confidence Bands}\label{sec:SCBs} 
Throughout this section we require that the assumption of Theorem \ref{thm:2DfCLT} hold. 
\paragraph{Construction of the SCBs}
Functional delta residuals can be applied to construct simultaneous confidence bands (SCBs) for  $ H\big(\, \bmu^{(\br)} \big ) $. The easiest way to do so is based on $t$-statistics (see e.g., \citet{Telschow:2019} ). To do so, in the context of moment-based statistics for $s \in S$ and a quantile $q_{\alpha, N} \in \mathbb{R}$, we define the collection of intervals $ \text{SCB}( s, N, q_{\alpha,N} ) $ with endpoints
\begin{equation}\label{eq:SCBconstruction}
   H\Big(\, \hat{\bmu}^{(\br)}_N(s) \Big ) \pm q_{\alpha,N} \,se\Big[\, H\Big(\, \hat{\bmu}^{(\br)}_N(s) \Big ) \,\Big]\,.
\end{equation}
The sample estimator $ \hat{\bmu}^{(\br)}_N $ fulfills a fCLT, see Theorem \ref{thm:2DfCLT}. Moreover, since the paths of $X$ are $\mathcal{L}^2$-H\"older continuous, $ \hat{\bmu}^{(\br)}_N $ is also strongly consistent, i.e., $ \Vert\, \hat{\bmu}^{(\br)}_N - \bmu^{(\br)} \,\Vert_\infty \rightarrow 0 $ almost surely by Lemma \ref{Lem:UniformConvergence}(b).
The standard error of the estimator is given by
\begin{equation}
\begin{split}
	se\Big[\, H\Big(\, \hat{\bmu}^{(\br)}_N(s) \Big ) \,\Big]
		&= \sqrt{\, \E\Bigg[\, H\Big(\, \hat{\bmu}^{(\br)}_N(s) \,\Big) - \E\Big[\, H\Big(\, \hat{\bmu}^{(\br)}_N(s) \,\Big)\,\Big]  \,\Bigg]^2  \,}\,,\label{eq:se_form}
\end{split}
\end{equation}
which, asymptotically, can be consistently estimated, by Corollary \ref{cor:momentfCLT}, using the sample variance of the functional delta residuals.
The intervals given by eq. \eqref{eq:SCBconstruction} are $(1-\alpha)$-SCBs, if the quantile $ q_{\alpha,N} $
satisfies
\begin{equation}\label{eq:nonasymQuant}
	\Prb\left(\, \max_{s \in S } \left\vert\, \frac{ H\Big(\, \hat{\bmu}^{(\br)}_N(s) \Big ) - H\Big(\, {\bmu}^{(\br)}(s)\Big)  }{ se\Big[\,H\Big(\, \hat{\bmu}^{(\br)}_N(s) \Big )\,\Big] } \,\right\vert > q_{\alpha,N} \,\right) = \alpha\,.
\end{equation}
The desired $ q_{\alpha,N} $ cannot be calculated easily in the finite sample because, in general, the standard error \eqref{eq:se_form} for finite $ N $ is hard to estimate. However, asymptotically the same arguments as in \citet{Chang2017}, yield the following asymptotic $(1-\alpha)$-SCBs for $H\big(\, {\bmu}^{(\br)}\big)$ and the asymptotic quantile $q_\alpha$ can be estimated using the delta residuals (as we will discuss different approaches to estimating this later on in this section).
\begin{theorem}\label{thm:SCB}
Under the assumptions of Theorem \ref{thm:2DfCLT} we define $q_{\alpha}$ such that
\begin{equation}\label{eq:asymQuant}
	\Prb\left(\, \max_{s \in S } \left\vert\, \frac{\tilde G(s)}{\sqrt{\,\var\big[\tilde G(s)\big]}\,} \,\right\vert > q_{\alpha} \,\right) = \alpha\,,
\end{equation}
where $\tilde G$ is defined in Corollary \ref{cor:momentfCLT}.
Then
\begin{equation*}
  \lim_{ N \rightarrow \infty }\Prb\,\Bigg(\, H\Big(\, {\bmu}^{(\br)}(s) \,\Big) \in \text{SCB}(\,s,N,q_\alpha\,) ~\text{ for all }s\in S \,\Bigg) = 1-\alpha\,.
\end{equation*}
\end{theorem}
\begin{proof}
Apply Corollary \ref{cor:momentfCLT} together with Slutsky's Lemma.
\end{proof}

If $ H $ is non-linear then the estimator $ H\big(\, \hat{\bmu}^{(\br)}_N \,\big ) $ for $H\big(\,\bmu^{(\br)}\,\big)$ is biased for finite $N$. As such it is possible to define bias corrected SCBs (as discussed in  \citet{Liebl2019}) to improve the finite sample size coverage of the SCBs. To do so we define, for $s \in S$, their endpoints as
\begin{equation}\label{eq:SCBconstructionBias}
   H\Big(\, \hat{\bmu}^{(\br)}_N(s) \,\Big ) - bias\Big[\, H\Big(\, \hat{\bmu}^{(\br)}_N(s) \,\Big) \,\Big] \pm q_{\alpha,N} se\Big[H \Big(\, \hat{\bmu}^{(\br)}_N(s)\, \Big )\Big]\,.
\end{equation}
Using $\nabla^2H$ to denote the Hessian of $H$, the bias can be approximated as follows, using a Taylor expansion of $ H $ around $ {\bmu}^{(\br)} $,
\begin{equation}
\begin{split}\label{eq:biasApproxHeuristic}
	bias\Big[\, H\Big(\, \hat{\bmu}^{(\br)}_N(s) \,\Big ) \,\Big]
		&= \E\left[\, H\Big(\, \hat{\bmu}^{(\br)}_N(s) \,\Big) \right] - H\Big(\, {\bmu}^{(\br)}(s)\,\Big)\\
		&\approx \E\left[\, \nabla H_{\bmu^{(\br)}}\Big(\, \hat{\bmu}^{(\br)}_N(s) - \bmu(s) \,\Big) + \frac{1}{2}\nabla^2H_{\bmu^{(\br)}}\Big(\, \hat{\bmu}^{(\br)}_N(s) - \bmu(s), \hat{\bmu}^{(\br)}_N(s) - \bmu(s) \,\Big) \right] \\
		&= \frac{1}{2} \E\left[\, \nabla^2H_{\bmu^{(\br)}}\Big(\, \hat{\bmu}^{(\br)}_N(s) - \bmu(s), \hat{\bmu}^{(\br)}_N(s) - \bmu(s) \,\Big) \,\right] \\
		&= \frac{1}{2} \sum_{k=1}^K\sum_{k'=1}^K \frac{\partial^2H}{ \partial s_k\partial s_{k'}}\Big(\,{\bmu^{(\br)}}(s)\,\Big) \cdot \E\left[\, \big(\, \hat{\mu}^{(r_k)}_N(s) - \mu^{(r_k)}(s)\,\big)\,  \big(\, \hat{\mu}^{(r_{k'})}_N(s) - \mu^{(r_{k'}) }(s) \,\Big) \right] \\
		&= \frac{1}{2N} \sum_{k=1}^K\sum_{k'=1}^K \frac{\partial^2H}{ \partial s_k\partial s_{k'}}\Big(\,{\bmu^{(\br)}}(s) \,\Big) \cdot \Big(\,\mu^{(r_k+r_{k'})}(s) - \mu^{(r_{k'})}(s)\,\mu^{(r_k)}(s)\,\Big)\,.
\end{split}
\end{equation}
In our simulations we use a simple plugin estimator based on the strongly consistent estimator
$ \hat{\bmu}^{(\br)}_N $ from equation \eqref{eq:biasApproxHeuristic}
\begin{equation}
\begin{split}
	\widehat{bias}\Big[\, H\Big(\, \hat{\bmu}^{(\br)}_N \,\Big ) \,\Big]
		= \frac{1}{2N} \sum_{k=1}^K\sum_{k'=1}^K \frac{\partial^2H}{ \partial s_k\partial s_{k'}}\Big(\,{\hat{\bmu}_N^{(\br)}}\,\Big) \cdot \Big(\,\hat\mu^{(r_k+r_{k'})}_{N} - \hat\mu^{(r_{k'})}_N\,\hat\mu^{(r_k)}_N\,\Big)\,.\label{eq:bias_est}
\end{split}
\end{equation}
to estimate the bias. Strong consistency of $ \hat{\bmu}^{(\br)}_N $
implies that this bias estimate is strongly consistent, too. Plugging it into the SCBs \eqref{eq:SCBconstructionBias} yields our bias corrected SCBs.

\paragraph{Estimation of the quantile $q_\alpha$}
We have described one approach to calculate
 the asymptotic quantile $q_\alpha$, i.e., estimating it using the $\alpha$-quantile of the maximum of the multiplier bootstrap process based on the functional delta residuals from Section \ref{sec:MultfCLT}. This is an adaptation of \citet{Chang2017} and Corollary \ref{cor:mbootconsistency} implies that this estimate is consistent for $q_\alpha$.

 A second approach assumes that the residuals have $C^3$ sample paths and utilizes the Gaussian kinematic formula, compare \citep{Telschow:2019}. Here the quantile $q_\alpha$ is approximated by exploiting the fact that for large $ u \in \mathbb{R}$,
\begin{equation*}
 \Prb\Big(\, \max_{s\in  S}  \hat G(s) > u \,\Big) \approx \cL_0\, \Phi^+(u) + \sum_{d=1}^{D} \cL_{d}\, \rho_d(u)\,,
\end{equation*}
as shown in \citet{Taylor2005}. Here $\hat G(s) = \tilde G(s) \,/ \sqrt{\,\var\big[\,\tilde G(s)\,\big]\,} $. The functions $
\rho_d(u) = (2\pi)^{-(d+1)/2} \mathcal{H}_{d-1}(u) e^{-u^2/2}$, $d=1,\ldots,D$, are the so-called Euler characteristic densities, where $\mathcal{H}_d$ is the $d$-th Hermite polynomial and $\Phi^+(u) = \Prb\big(\, N(\,0,1\,) > u\,\big)$. The coefficients $\cL_0,\ldots,\cL_{D}$ are referred to as the Lipshitz-Killing curvatures of $ S$, which are intrinsic volumes of $ S$ considered as a Riemannian manifold endowed with a Riemannian metric induced by $\hat G$ \citep[Chapter 12]{Adler:2009}. In particular, $\cL_0=\chi( \,S)$ is the Euler characteristic of the set $ S$, which is usually known.
Given consistent estimators $\hat\cL_{1},\ldots,\hat\cL_{D}$ of the  Lipshitz-Killing curvatures an estimate $\hat q_\alpha$ of $q_\alpha$ can be found by finding the largest $u$ such that
\begin{equation*}
    \cL_0 \,\Phi^+(u) + \sum_{d=1}^{D} \hat\cL_{d}\, \rho_d(u) = \alpha\,.
\end{equation*}

Currently there are only a few works dealing with estimation of  Lipshitz-Killing curvatures for nonstationary processes and arbitrary dimensional domains $S$, see \citet{Taylor:2007, Schwartzman:2019}. The estimators from the last two sources require residuals which asymptotically have the covariance structure of the limiting process $\hat G$. The functional delta residuals satisfy this, once we normalize them to have empirical variance $1$.
 We refer the reader to \citet{Telschow:2019} for more details on how to use the Gaussian kinematic formula to
estimate quantiles for SCBs.

\paragraph{Testing Gaussianity using Skewness and Kurtosis}
It is well-known that confidence intervals can be inverted to tests, see for example \cite{Lehmann2005}. By the same reasoning, simultaneous confidence bands for skewness and kurtosis can be used to test whether a sample is Gaussian. In order to do so we compute the simultaneous confidence bands for the pointwise skewness and kurtosis given in  equations \eqref{eq:g1} and \eqref{eq:g2} under the assumption of Gaussianity. Assume that
$X_1,\ldots,X_N\iid X$ with $X$ being a Gaussian process on $S$. Using the formulas for the sample variance of $g_1^N$ and $g_2^N$ from \cite{Fisher:1930}, as well as $\E\big[g_1^N\big] = 0$ and $\E\big[g_2^N\big] = -6\,/(N+1)$ for Gaussian samples, Theorem \ref{thm:SCB} implies, if $X$ is Gaussian,  that for any $\alpha \in (0,1)$
\begin{equation}
\begin{split}
	&\lim_{N\rightarrow\infty}\Prb\left(\, \max_{s \in S} \vert\, g_1^N(s) \,\vert >  q_\alpha^{(g_1^N)}\,\sqrt{\frac{6\,(\,N-2\,)}{(\,N+1\,)\,(\,N+3\,)}\,} \,\right) = \alpha\,,\\
	&\lim_{N\rightarrow\infty}\Prb\left(\, \max_{s \in S} \Bigg\vert\, g_2^N(s) + \frac{6}{N+1} \,\Bigg\vert >  q_\alpha^{(g_2^N)}\,\sqrt{\frac{24\,N\,(\,N-2\,)\,(\,N-3\,)}{(\,N+1\,)^2\,(\,N+3\,)\,(\,N+5\,)\,}} \,\right) = \alpha\,.
	\end{split}
\end{equation}
Here $q_\alpha^{(g_1^N)}$ and $q_\alpha^{(g_2^N)}$ are the quantiles obtained from applying Theorem \ref{thm:SCB} to $g_1^N$
and $g_2^N$. Note that the only difference, with regards to the construction of the SCBs, is that we have replaced the
unknown standard error and the bias by the known Gaussian quantities. This shows that under the null hypothesis that
$X$ is Gaussian, the tests rejecting Gaussianity if
\begin{equation}
\begin{split}
	&\max_{s\in S} \vert g_1^N(s) \vert >  q_\alpha^{(g_1^N)}\,\sqrt{\frac{6\,(\,N-2\,)}{(\,N+1\,)\,(\,N+3\,)}\,}\,,\\
	&\max_{s \in S} \Bigg\vert\, g_2^N(s) + \frac{6}{N+1} \,\Bigg\vert >  q_\alpha^{(g_2^N)}\,\sqrt{\frac{24\,N\,(\,N-2\,)\,(\,N-3\,)}{(\,N+1\,)^2\,(\,N+3\,)\,(\,N+5\,)\,}}
\end{split}
\end{equation}
are asymptotically exact. Hence either of the SCBs can be used to test departure for Gaussianity. The skewness SCBs check whether there is departure from Gaussianity due to significant non-zero skewness for any $s \in S$, while the excess kurtosis SCBs do the same for significant non-zero excess kurtosis.

In the simulations from Section \ref{scn:Simulations} we will see that nominal coverage for the above bands requires very large $N$. This is partially because the rate of convergence of $g_1^N$ and $g_2^N$ to their asymptotic normal distributions is
very slow. To solve this problem \cite{DAgostino:1970} proposed a transformation making $g_1$ approximately standard normal
for any $N\geq 8$ and \cite{Anscombe:1983} proposed a transformation doing the same for $g_2$ for $N \geq 20$. For samples of real valued random variables \cite{DAgostino:1990} argued that these transformations give powerful tests for Gaussianity.
The transformation which is applied to $g_1$ is given by
\begin{equation}\label{eq:Z1}
\begin{split}
	Z_{1,N}(x) &= \frac{1}{\sqrt{\,\log(W_N)N}\,} {\rm asinh}\left(\, \frac{x}{\alpha_N\sqrt{\,c_{1,N}\,}} \,\right) \xrightarrow{N\rightarrow\infty} 0.335\cdot {\rm asinh}\big(\, 1.216\cdot x \,\big)\,,\\
	&~ ~ ~ ~ ~W_N^2      = \sqrt{\,2c_{2,N} - 2\,} - 1\,,~ ~ ~ ~
	\alpha_N^2 = \frac{2}{W_N^2-1}\,.
\end{split}
\end{equation}
where $c_{1,N}$ and $c_{2,N}$ are constants depending only on $N$, which can be found in \cite{DAgostino:1970}.
The transformation which is applied to $g_2$ is given by
\begin{equation}\label{eq:Z2}
	\begin{split}
	Z_{2,N}(x) &= \sqrt{\,\frac{9A_N}{ 2N} \,}\, \left(\, 1 - \frac{2}{9A_N} - \left(\, \frac{1 - 2/A_N}{ 1 + \frac{x + 3 - b_{1,N}}{ \sqrt{\,b_{2,N}\,} } \sqrt{\, 2 / (A_N - 4) \,} } \,\right)^{1/3} \,\right)
	  \xrightarrow{N\rightarrow\infty} 0.8165 \cdot \left(\, 1 - \left(\,  \frac{1}{1 + 0.75\cdot (x+3)}  \,\right)^{1/3} \,\right)\,,\\
	A_N &= 6 + \frac{8}{\sqrt{\,b_{3,N}\,}} \left(\, \frac{2}{\sqrt{\,b_{3,N}\,}} + \sqrt{1 + \frac{4}{b_{3,N}}} \,\right).
	  \end{split}
\end{equation}
Here $b_{1,N} $, $b_{2,N}$ and $b_{3,N}$ are again constants depending on $N$, which can be found in \cite{DAgostino:1990}.

The transformation $Z_{1,N}\big(g_1^N(s)\big)$ of the sample skewness and the transformation $Z_{1,N}\big(g_1^N(s)\big)$ of the sample excess kurtosis are both moment-based statistics and satisfy the assumptions from Remark \ref{remark:NdependExtension} as is proven in \ref{app:auxlemmas}. Since, under the assumption that $X$ is Gaussian, the transformed estimators are unit-variance and zero-mean, we obtain asymptotic (1-$\alpha$)-SCBs for the transformed skewness with endpoints
\begin{equation}\label{eq:TestGaussNormalSkew}
	Z_{1,N}\big(g_1^N(s)\big) \pm  q_\alpha^{Z_{1,N}(g_1^N)}
\end{equation}
and therefore reject Gaussianity of the sample due to non-zero skewness, if zero is not contained in all of these intervals. 
Similarly, we can test departure of Gaussianity due to non-zero excess-kurtosis with the SCBs with endpoints
\begin{equation}\label{eq:TestGaussNormalKurt}
	Z_{2,N}\big(g_2^N(s)\big) \pm  q_\alpha^{Z_{2,N}(g_2^N)}\,.
\end{equation}
Since these transformations are bijective for all $N\in \mathbb{N}\cup \{ \infty \}$ the SCBs for the transformed
variable can be turned into (1-$\alpha$)-SCBs for skewness and excess kurtosis by applying their inverse. These SCBs
might be easier to interpret.

%
%
\section{Simulations of Coverage for Simultaneous Confidence Bands}\label{scn:Simulations}

In this section we study the coverage rate of simultaneous confidence bands for different moment-based statistics. We use $20,000$ Monte Carlo simulations to assess the coverage and evaluate the processes on a grid of $[0,1]$ composed of $175$ equally spaced points. The bands are calculated as described in Section \ref{sec:SCBs}.
In particular, the quantile $q_\alpha$ is estimated using either the functional delta residuals through the multiplier bootstraps given in \cite{Chang2017, Telschow:2019} with Rademacher (rMult/rtMult) or Gaussian multipliers (Mult/tMult) using $5000$ bootstrap replicates or the Gaussian kinematic formula (GKF/tGKF), see Section \ref{sec:SCBs}. Here the $ "t" $ for multiplier bootstraps denotes the $t$-multiplier bootstrap and tGKF refers to the Gaussian kinematic formula for a $t$-process rather than the
Gaussian one, see e.g. \cite{Telschow:2019}.
We compare different constructions for the SCBs with and without the bias correction (eq. \eqref{eq:bias_est}). The results for bias correction (except for Cohen's $d$) are deferred to the appendix, since the coverage rates of the SCBs in general are closer to nominal without estimating the bias. In the Gaussian simulations the standard error (s.e.) for the moment-based statistics are known. Hence in such cases we compare SCBs using the known Gaussian s.e. and SCBs using the s.e. estimate obtained from the  empirical variance of the functional delta residuals \eqref{eq:covrates_estim}. If the known Gaussian s.e. is used, this is indicated in the titles of the plots by writing \textit{Gaussian S.E.}. For the estimated s.e. we use \textit{Estimated S.E.} in the titles.

\subsection{Functional Models of the Simulations}
\begin{figure}[h]
\begin{center}
		\includegraphics[trim=0 0 0 0,clip,width=1.5in]{\figurepath 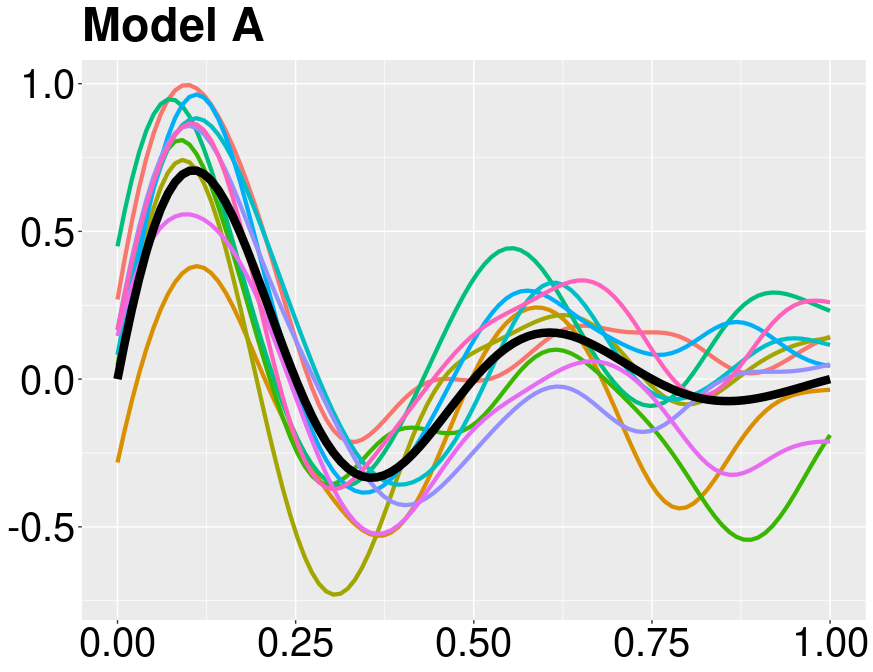}
		\includegraphics[trim=0 0 0 0,clip,width=1.5in]{\figurepath 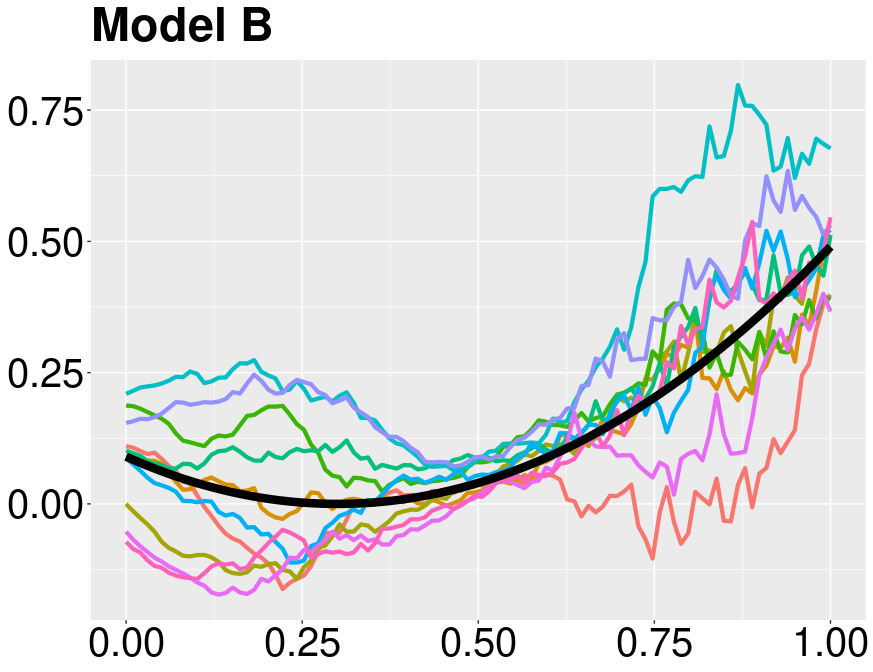}
		\includegraphics[trim=0 0 0 0,clip,width=1.5in]{\figurepath 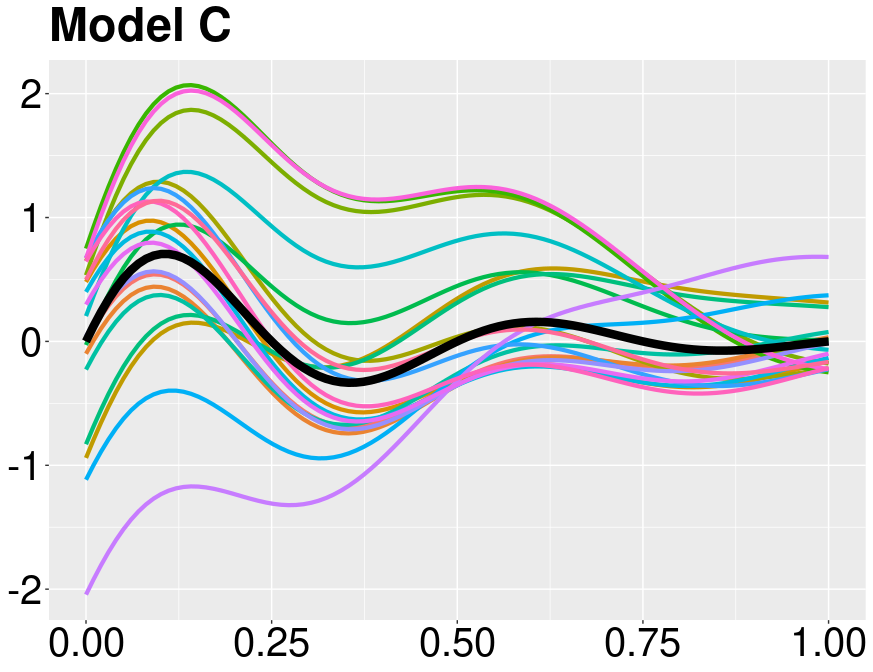}
	\caption{Examples of ten sample paths of the considered models. The bold black line is
		the true population mean. \label{fig:SamplePaths}}
\end{center}
\end{figure}

We compare the following three models defined on $S = [0,1]$:
    \begin{align*}
	      \mathbf{Model~A:}~~ Y^A(s)
    &= \sin(\,4 \pi s\,) \exp(\,-3 s\,) +  \frac{(\,1 - s - 0.4\,)^2 + 1 }{6}\cdot \frac{ \mathbf{a}^T\mathbf{K}^A(s)}{\Vert\, \mathbf{K}^A(s) \,\Vert} \\
	      \mathbf{Model~B:}~~Y^B(s)
	&= (\, s-0.3 \,)^2 +  \frac{\sin(\,3\,\pi\, s\,) + 1.5}{ 6 }\cdot Z_B(s) \\
	      \mathbf{Model~C:}~~Y^C(s)
	&= \sin(\,4 \pi s\,) \exp(\,-3 s\,) + (\,1.5 - s\,)\cdot Z_C(s)
    \end{align*}
    with $ \mathbf{K}^A(s) $ being a vector with entries
    $ K^A_i(s) = \exp\big(\, -\frac{ \left(\,s-x_i\,\right)^2 }{ 2h_i^2 } \,\big) $ for $x_i=i\,/\,21$ and $ \mathbf{a} \sim N(\,0, I_{21 \times 21}\,) $. Hence Model \textbf{A} is a smooth non-stationary Gaussian process.  The error process $ Z_B $
    is the zero-mean Gaussian process having the non-stationary modified Matern-type covariance $ c_{Z_B}(\,s,t\,) = c(\,s,t\,)\, / \sqrt{\,c(\,s,s\,)\, c(\,t,t\,)\,} $ with
    \begin{equation*}
    \begin{aligned}
    	c(\,s,t\,) = 0.4^2\, \Big(\, 2^{1-\nu_{ts}}/\,\Gamma(\nu_{ts}) \,\Big)
    	\,\Big(\, \sqrt{\,2\nu_{ts}}\,\vert\, t-s \,\vert \,\Big)^{\nu_{ts}}
    	K_{\nu_{ts}}\Big(\,\sqrt{\,2\nu_{ts}}\,\vert\, t-s \,\vert \,\Big)\,,
    \end{aligned}
    \end{equation*}
    where $\nu_{t,s} = 1 - 3\sqrt{\,\max(\,s,t\,)\,}\,/\,4$, compare the simulation section in \citet{Liebl2019}. This model has continuous but non-differentiable paths. Model \textbf{C} is a smooth non-Gaussian error process of the form $Z_C(s) = W(s) \,/ \sqrt{\,\var[\,W(s)\,]\,} $ with
    \begin{equation}
    	W(s) = \frac{\sqrt{2}}{6}(\, \eta_1-1 \,)\sin(\,\pi s\,) + \frac{2}{3}(\, \eta_2-1 \,)\,(\,s-0.5\,)\,,~ \eta_1 \sim \chi_1^2\,,~ \eta_2 \sim \text{Exponential}(1)\,.
    \end{equation}
Examples of the sample paths of these processes can be found in Figure
    \ref{fig:SamplePaths}.
\vspace{-0.7cm}
\begin{figure}[h]
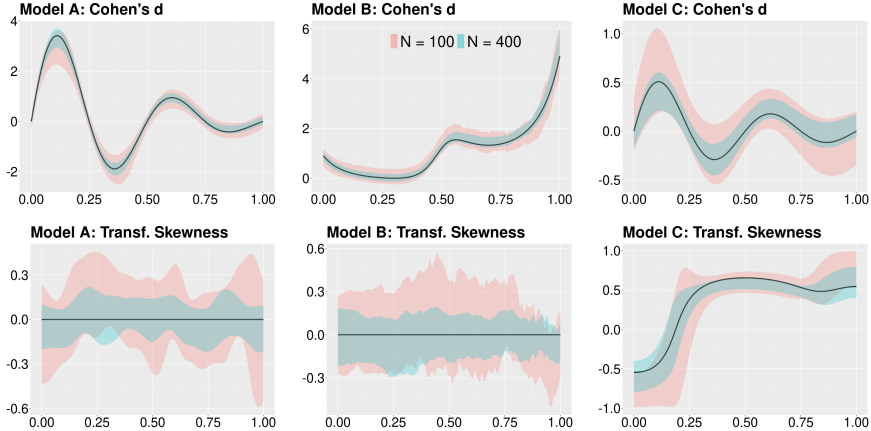

\begin{center}
		\includegraphics[trim=0 0 0 0,clip,width=1.5in]{\figurepath ModelA_SCBs_cohensd_example.pdf}
		\includegraphics[trim=0 0 0 0,clip,width=1.5in]{\figurepath ModelB_SCBs_cohensd_example.pdf}
		\includegraphics[trim=0 0 0 0,clip,width=1.5in]{\figurepath ModelC_SCBs_cohensd_example.pdf}\\
		\includegraphics[trim=0 0 0 0,clip,width=1.5in]{\figurepath ModelA_SCBs_skewness_normality_example.pdf}
		\includegraphics[trim=0 0 0 0,clip,width=1.5in]{\figurepath ModelB_SCBs_skewness_normality_example.pdf}
		\includegraphics[trim=0 0 0 0,clip,width=1.5in]{\figurepath ModelC_SCBs_skewness_normality_example.pdf}
	\caption{Examples of SCBs for moment-based statistics for two sample sizes. The bold black line is the true population parameter. \label{fig:SampleSCBs}}
\end{center}
\end{figure}

\FloatBarrier
\subsection{Coverage Rates for Cohen's $d$}
\begin{figure}[h]
\begin{center}
	\includegraphics[trim=0 0 0 0,clip,width=1.5in]{\figurepath 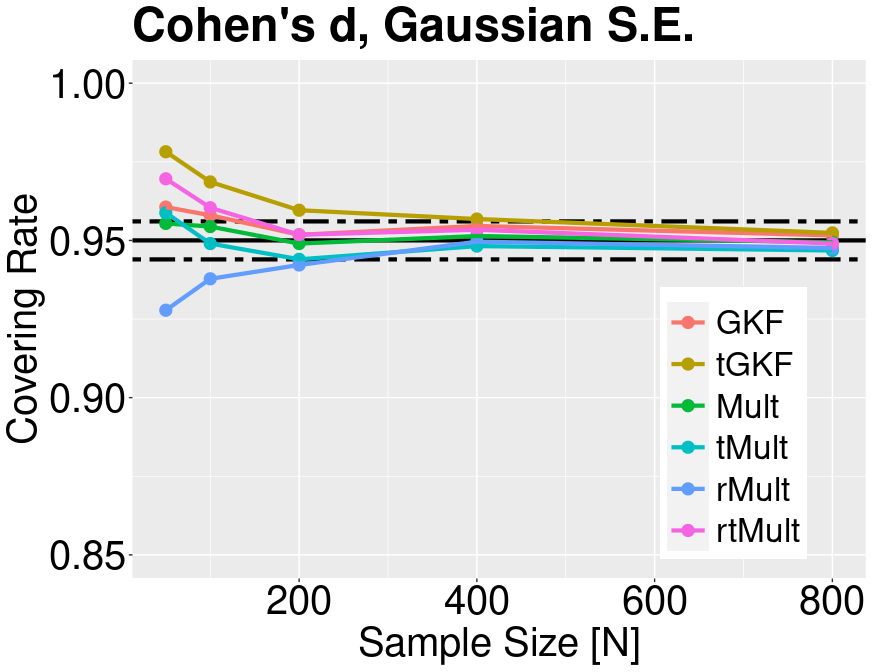}
	\includegraphics[trim=0 0 0 0,clip,width=1.5in]{\figurepath 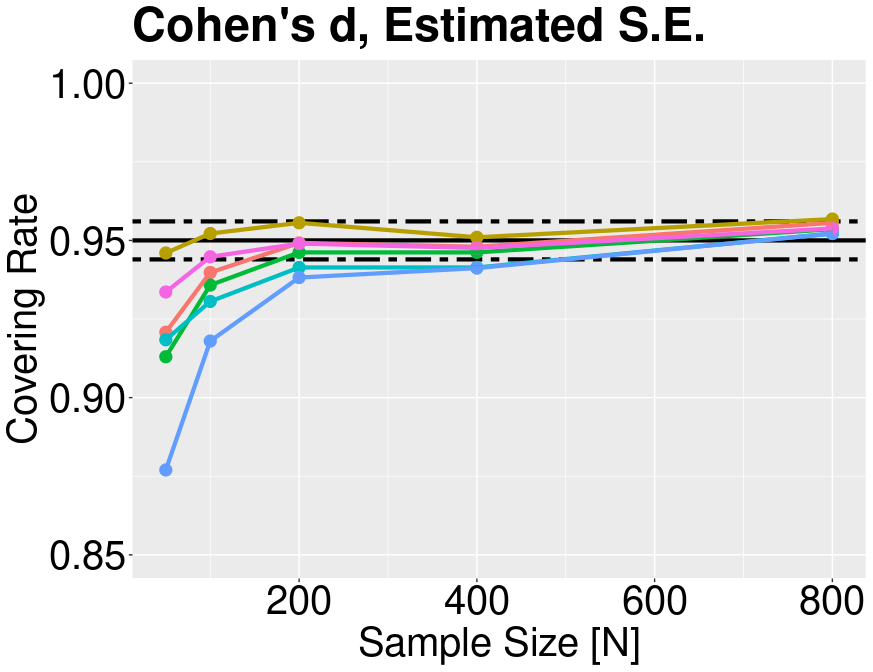}
	\includegraphics[trim=0 0 0 0,clip,width=1.5in]{\figurepath 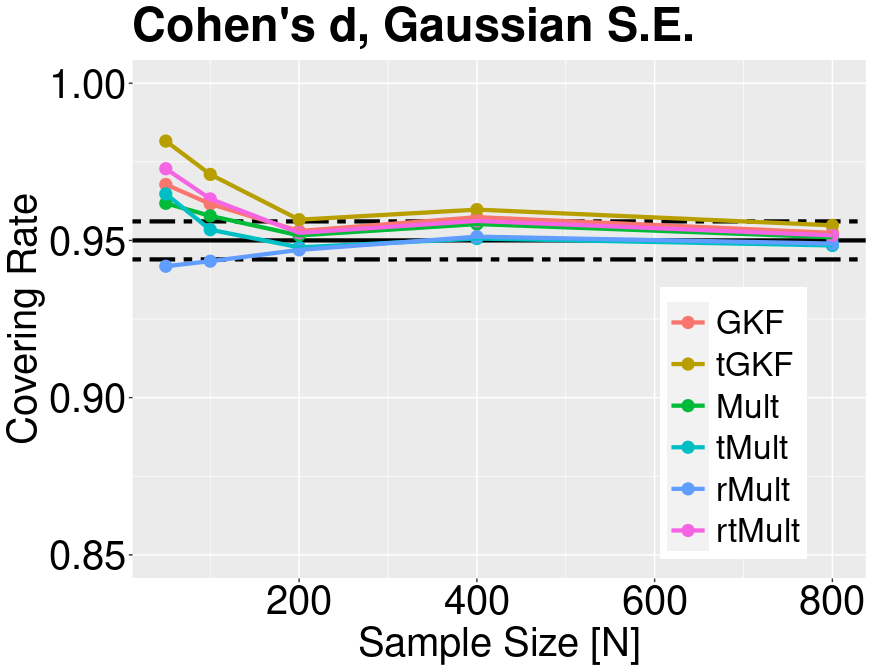}
	\includegraphics[trim=0 0 0 0,clip,width=1.5in]{\figurepath 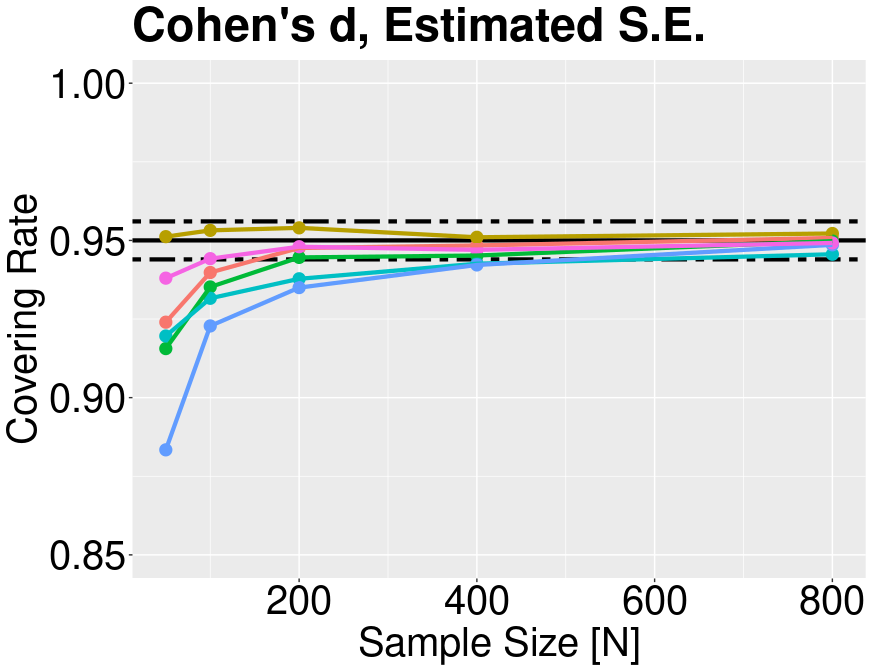}
	\caption{Simulations of coverage rates of SCBs for Model A. In the first two panels no bias correction is used in the construction of the SCBs, while in the third and fourth the coverage rate of bias-corrected SCBs is reported. The black dashed lines are
		95\% confidence intervals for the nominal level $0.95$. \label{fig:SCBsModelA}}
\end{center}
\end{figure}

\begin{figure}[h]
\begin{center}
	\includegraphics[trim=0 0 0 0,clip,width=1.5in]{\figurepath 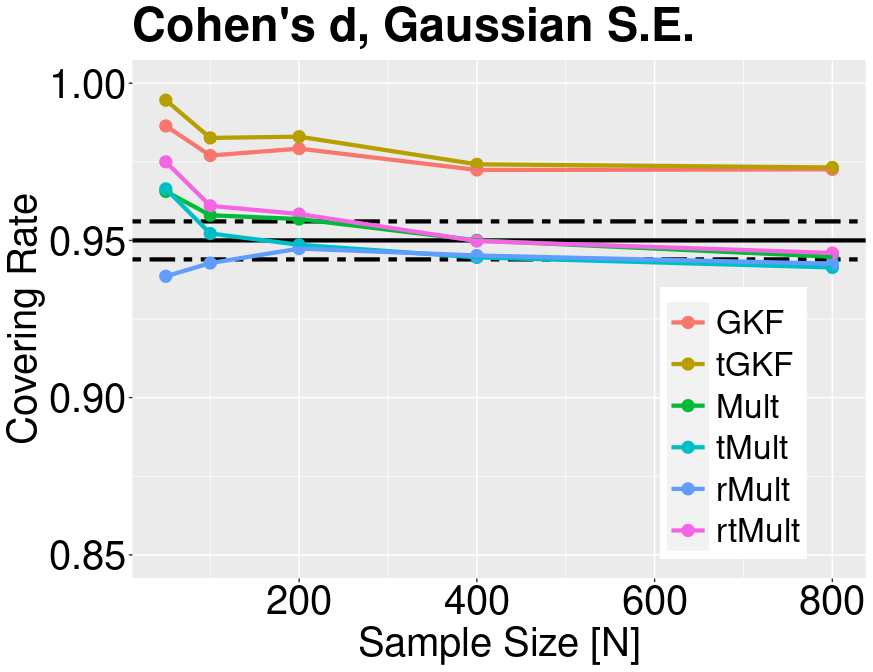}
	\includegraphics[trim=0 0 0 0,clip,width=1.5in]{\figurepath 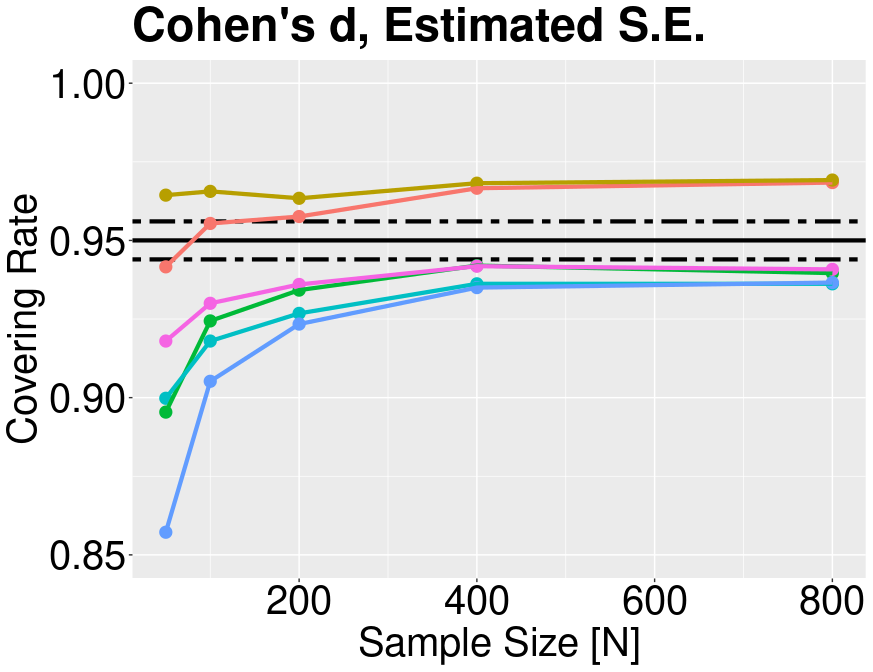}
	\includegraphics[trim=0 0 0 0,clip,width=1.5in]{\figurepath 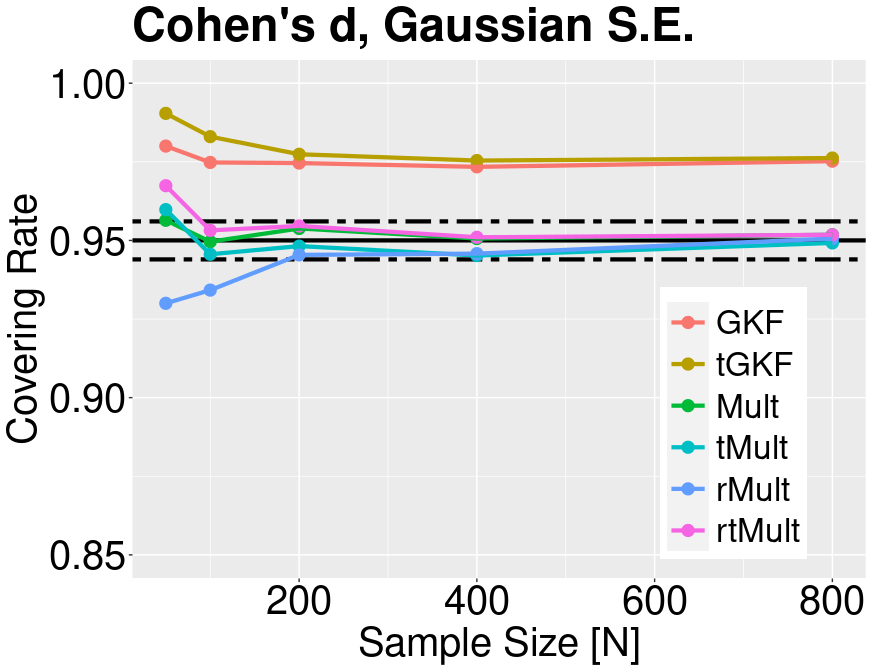}
	\includegraphics[trim=0 0 0 0,clip,width=1.5in]{\figurepath 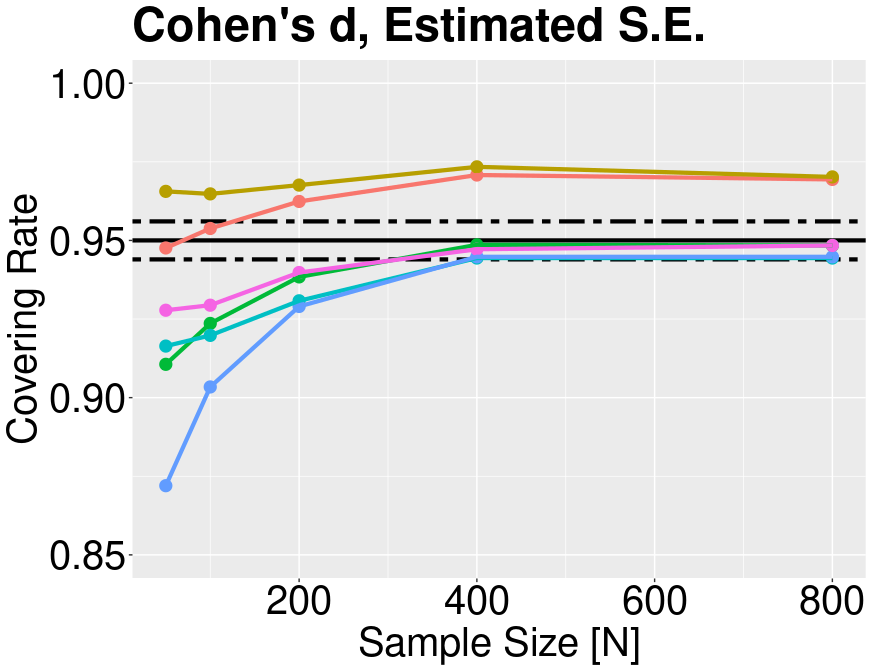}
	\caption{Simulations of coverage rates of SCBs for Model B. In the first two panels no bias correction
			is used in the construction of the SCBs, while in the third and fourth the coverage rate of
			bias-corrected SCBs is reported. The black dashed lines are 95\% confidence intervals
			 for the nominal level $0.95$.\label{fig:SCBsModelB}}
\end{center}
\end{figure}

\begin{figure}[h]
\begin{center}
\includegraphics[trim=0 0 0 0,clip,width=1.5in]{\figurepath 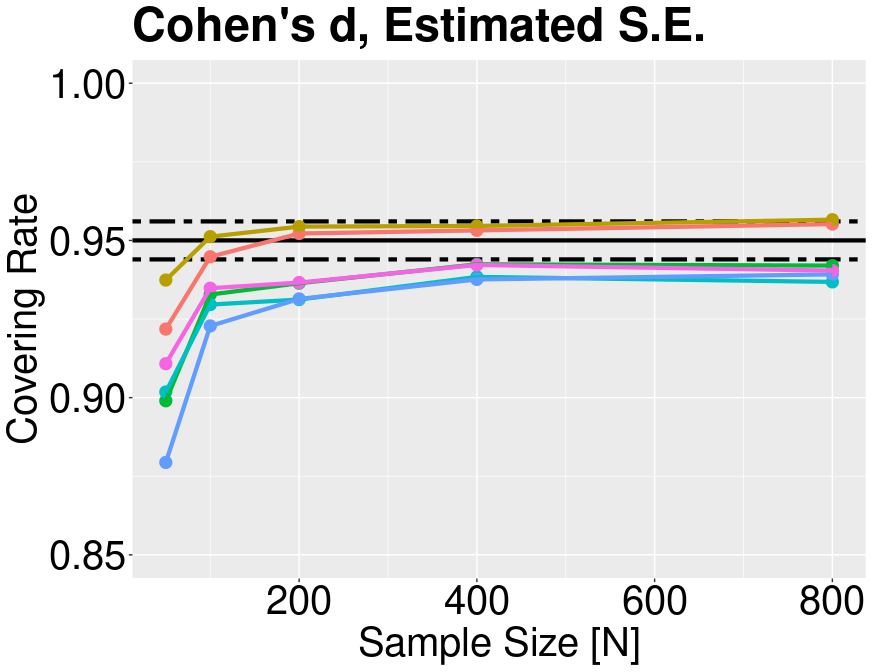}
\includegraphics[trim=0 0 0 0,clip,width=1.5in]{\figurepath 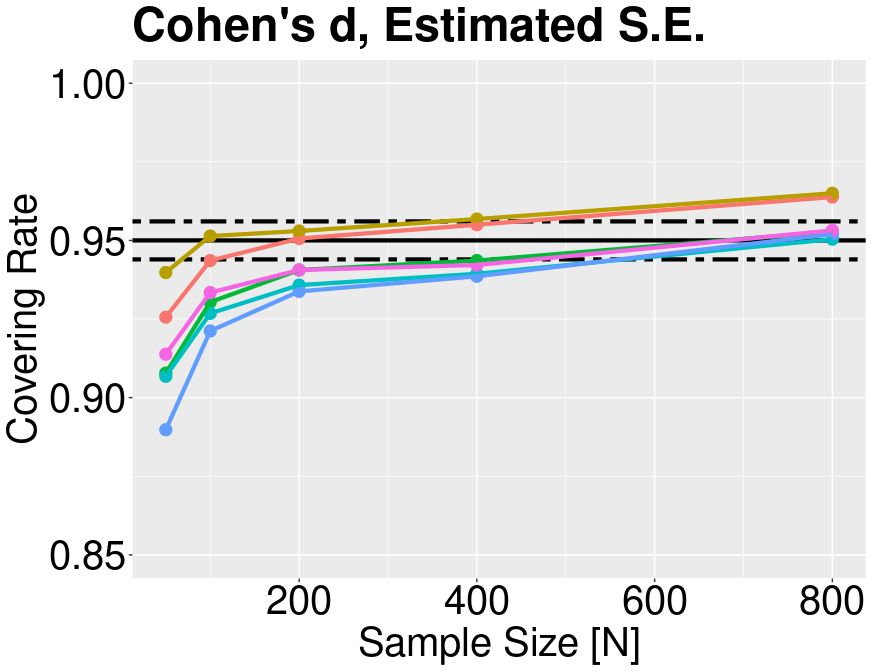}
	\caption{Simulations of coverage rates of SCBs for Model C.  The first panels are the coverage rats of SCBs without
the bias correction, while the second panel reports the coverage rate of bias-corrected SCBs. The black dashed lines are
		95\% confidence intervals for the nominal level $0.95$.\label{fig:SCBsModelC}}
\end{center}
\end{figure}

\FloatBarrier
Coverage rates of the SCBs for Model $ A $ are visualized in Figure \ref{fig:SCBsModelA}.
In this smooth Gaussian case the coverage rates are close
to the nominal coverage rate of $ 0.95 $ for all methods of quantile estimation if $ N $ is larger
than $ 200 $. The effect of estimating the bias is marginal.
If the variance is estimated then all methods apart from the tGKF, require roughly $ N = 400 $
to converge to nominal coverage. The tGKF reaches nominal coverage for
low sample sizes, which can be partly explained by the fact that the tGKF overestimates the coverage rates,
if the true variance is known. The main cause of the slow convergence of the coverage rates is that the
sample variance of the functional delta residuals substantially underestimates the true finite sample
variance of the Cohen's $d$ estimator for $ N \leq 200 $, compare Figure \ref{fig:VarianceEstimates}.
The results for Model $B$ are shown in Figure \ref{fig:SCBsModelB}. They are similar to the results
of Model $A$ with the exception that the GKF and tGKF
have over-coverage. This occurs because the GKF methods require
$ C^2 $ sample paths. 
For the non-Gaussian model $C$ the true variance of Cohen's $d$ is
not known. Hence we only simulated SCBs with estimated variances, see Figure
\ref{fig:SCBsModelC}. All methods, except those that use the GKF, have a coverage which converges asymptotically to the nominal level. The GKF approaches are slightly conservative. This holds because the EC heuristic for a Gaussian process $X$ over an interval $S\subset\mathbb{R}$, used to justify the GKF, only gives an upper bound of the probability $\Pr\big(\, \max_{s \in S} X_s \geq u \,\big)$, see for example \cite{Liebl2019}.

 \vspace{-0.7cm}

\FloatBarrier
\subsection{Coverage Rates for Skewness and Kurtosis}
\begin{figure}[h]
\begin{center}
\includegraphics[trim=0 0 0 0,clip,width=1.5in]{\figurepath 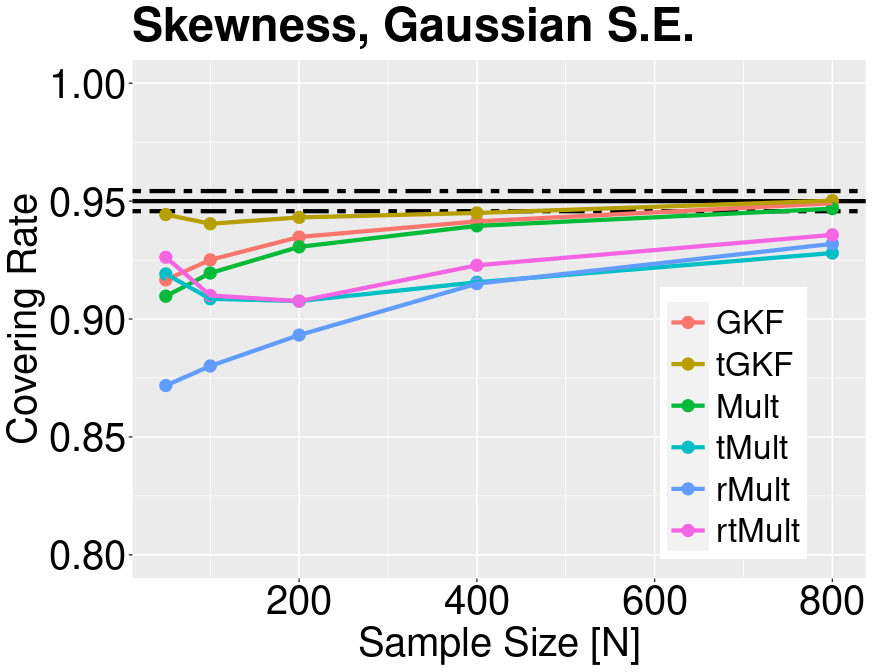}
\includegraphics[trim=0 0 0 0,clip,width=1.5in]{\figurepath 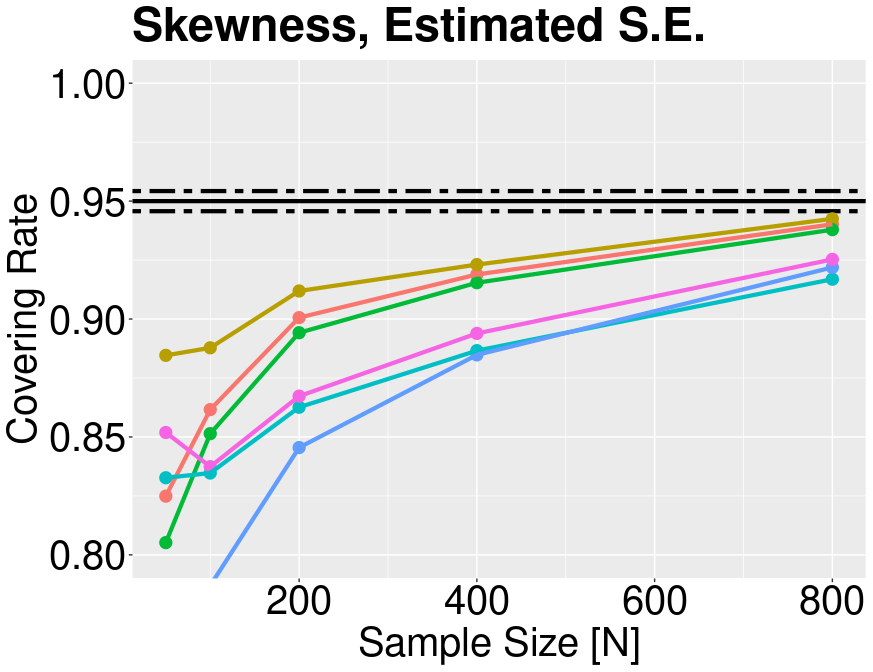}
\includegraphics[trim=0 0 0 0,clip,width=1.5in]{\figurepath 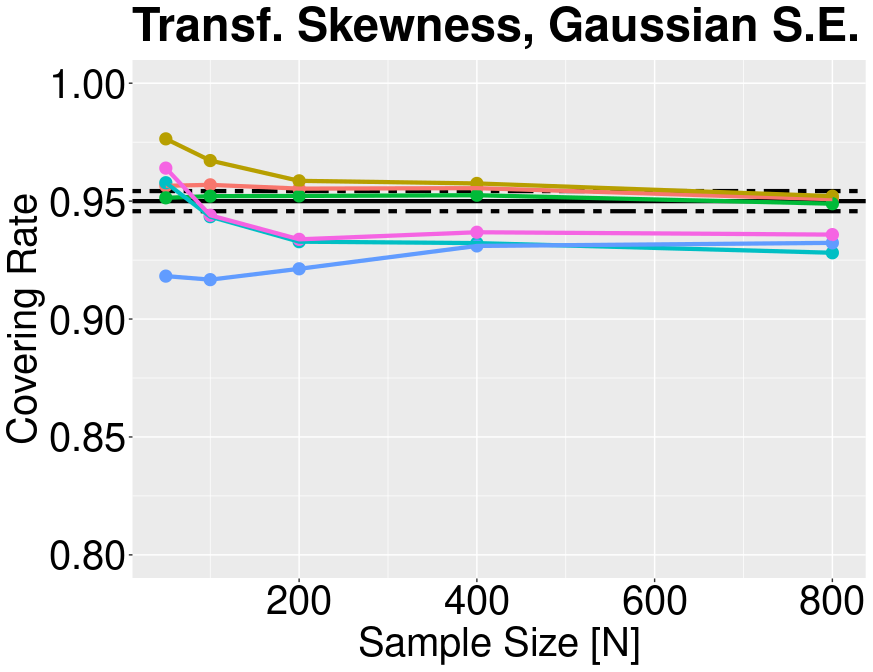}
\includegraphics[trim=0 0 0 0,clip,width=1.5in]{\figurepath 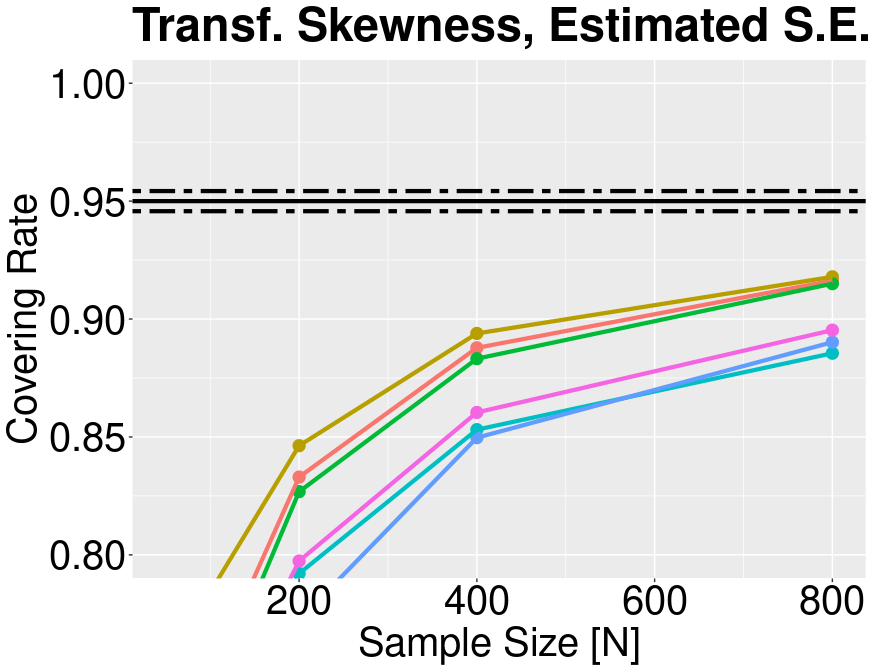}
	\caption{Simulations of coverage rates of SCBs for Model $A$. In these simulations no bias correction is applied. The black dashed lines are 95\% confidence intervals for the nominal level $0.95$. Special attention need to be placed on  the third panel. It shows that the test for Gaussianity based on the SCBs for transformed skewness given in eq. \eqref{eq:TestGaussNormalSkew} has the correct significance level even for small $N$ for this Gaussian process so long as the GKF or Mult method is used for estimation of the quantile. \label{fig:SCBsModelA-skew}}
\end{center}
\end{figure}

\FloatBarrier
For the smooth Gaussian Model $ A $ coverage rates of the SCBs for skewness and transformed skewness without including the bias estimate are shown in Figure \ref{fig:SCBsModelA-skew}. Transformed skewness means applying the additional transformation $Z_{1,N}$ to the moment-based statistic $g_1$ and constructing SCBs for the transformed parameter. Under Gaussianity the true transformed parameter is equal to zero for all $s\in S$. Convergence to nominal coverage of $0.95$ for skewness requires sample sizes larger than $ N \approx 800 $ for all methods of quantile estimation. Only when the true finite sample size variance is used, a good level of coverage is achieved  for sample sizes around $ N \approx 400 $, so long as the GKF or the Gaussian multiplier bootstrap is used for quantile estimation.
The reason behind this is twofold.

\begin{figure}[h]
\begin{center}
		\includegraphics[trim=0 0 0 0,clip,width=1.5in]{\figurepath 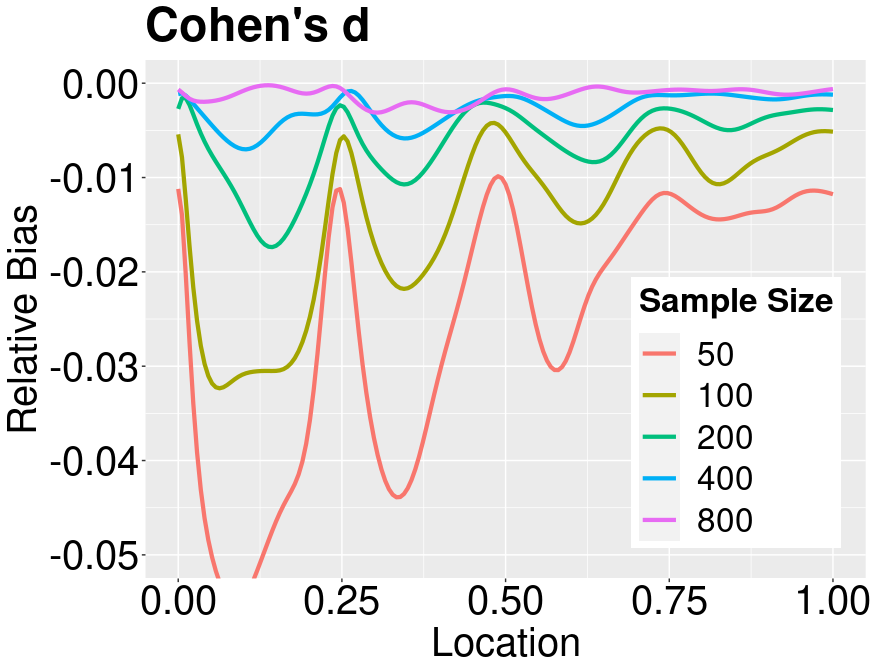}
		\includegraphics[trim=0 0 0 0,clip,width=1.5in]{\figurepath 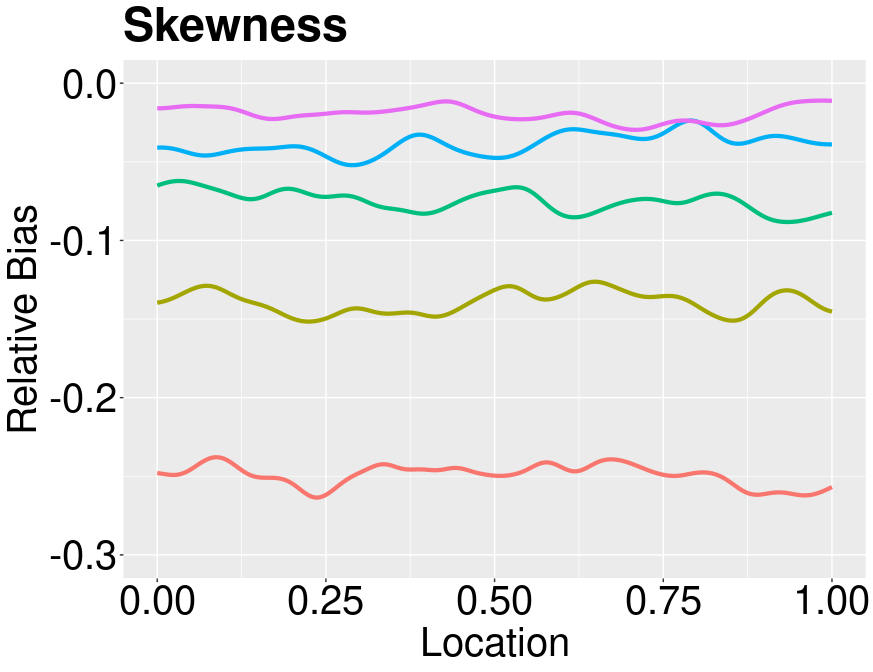}
		\includegraphics[trim=0 0 0 0,clip,width=1.5in]{\figurepath 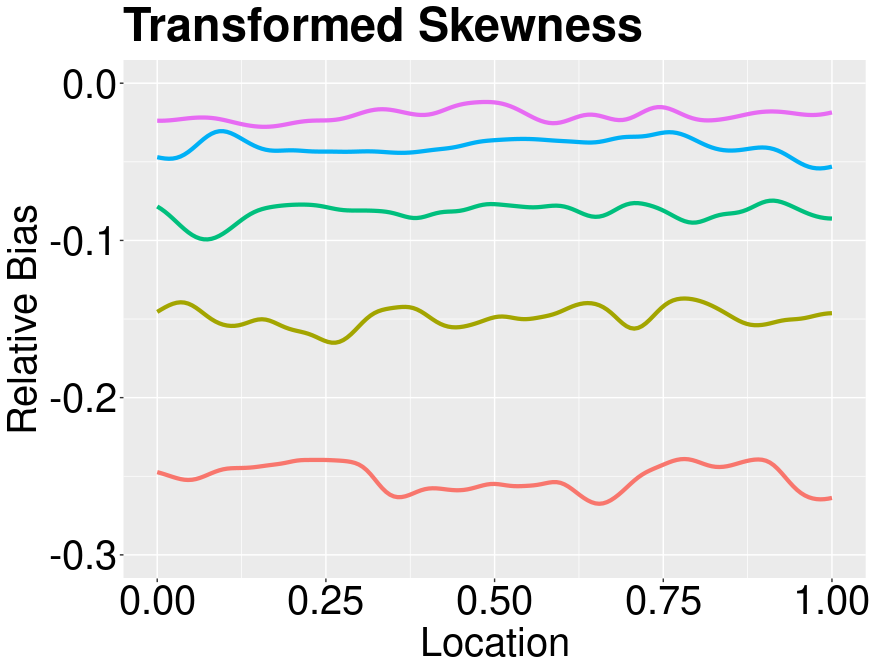}
	\caption{Dependency on sample size of the relative bias of s.e. estimates from the delta residuals for Model $A$. The colored curves are the estimates of the relative bias $\mathbb{E}\left[ \text{estimated s.e.} - \text{true s.e.} \right] / (\text{true s.e.}) $ obtained as the sample average of $5,000$ Monte Carlo simulations.\label{fig:VarianceEstimates}}
\end{center}
\end{figure}

First, Figure \ref{fig:VarianceEstimates} shows that for finite $ N $ the estimated variances from the functional delta residuals underestimate the true finite $ N $
variance massively (up to $30\%$ for $N = 50$) and only slowly converge to the correct finite sample variance of the skewness estimator. This problem is less severe for Cohen's $d$ since the finite $ N $ variance is only slightly underestimated (up to $5\%$ for $N = 50$). Hence finite sample coverage is improved, if we use the pre-knowledge of the true pointwise standard error of the skewness estimator under Gaussianity.

The second reason is that while the skewness estimator satisfies a fCLT, its convergence to a Gaussian process is slow. This can be remedied by using the transformed skewness parameter, which has a faster convergence to Gaussianity. In fact, SCBs for transformed skewness have almost exact coverage rates even for low $N$, if the true pointwise standard error of $1/\sqrt{N}$ for Gaussian samples is used. In particular, this shows that under the null hypothesis the test given in \eqref{eq:TestGaussNormalSkew} for Gaussianity using the transformed skewness SCBs has the correct significance level of $0.05$.
Similar as for skewness using the pre-knowledge of
the Gaussian variance is essential, since Figure \ref{fig:VarianceEstimates} shows that estimating the variance correctly requires large $N$.

\begin{figure}[h]
\begin{center}
\includegraphics[trim=0 0 0 0,clip,width=1.5in]{\figurepath 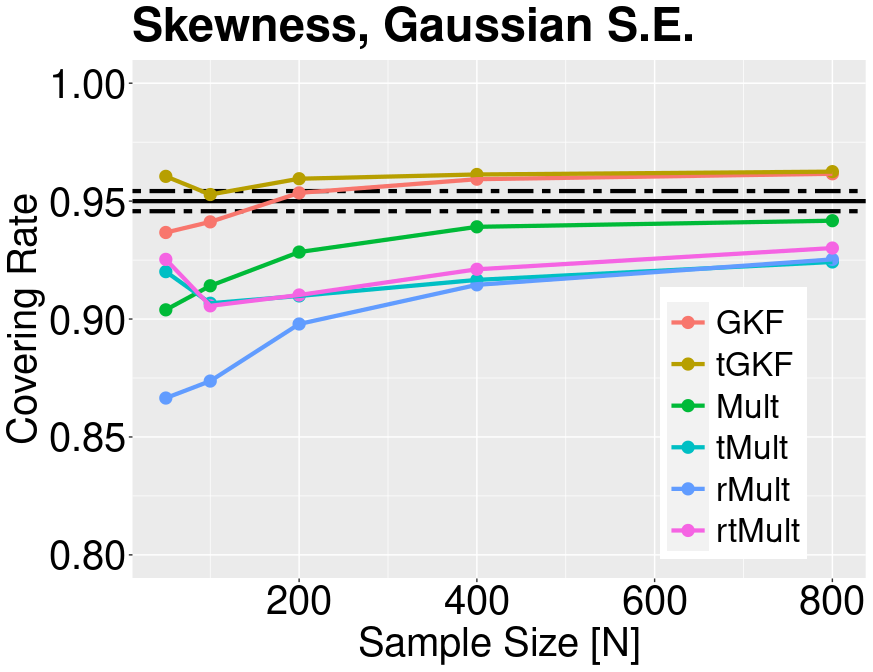}
\includegraphics[trim=0 0 0 0,clip,width=1.5in]{\figurepath 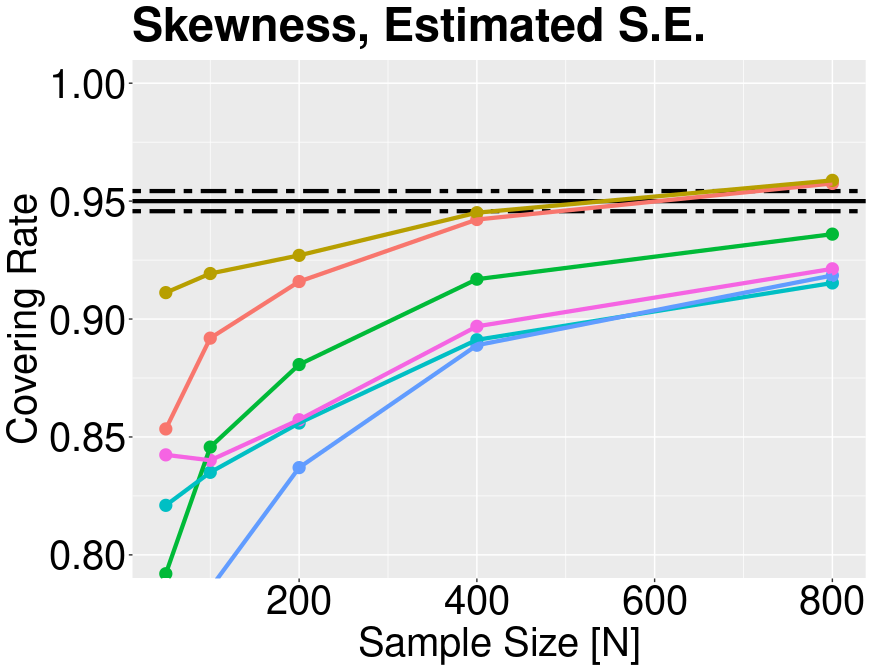}
\includegraphics[trim=0 0 0 0,clip,width=1.5in]{\figurepath 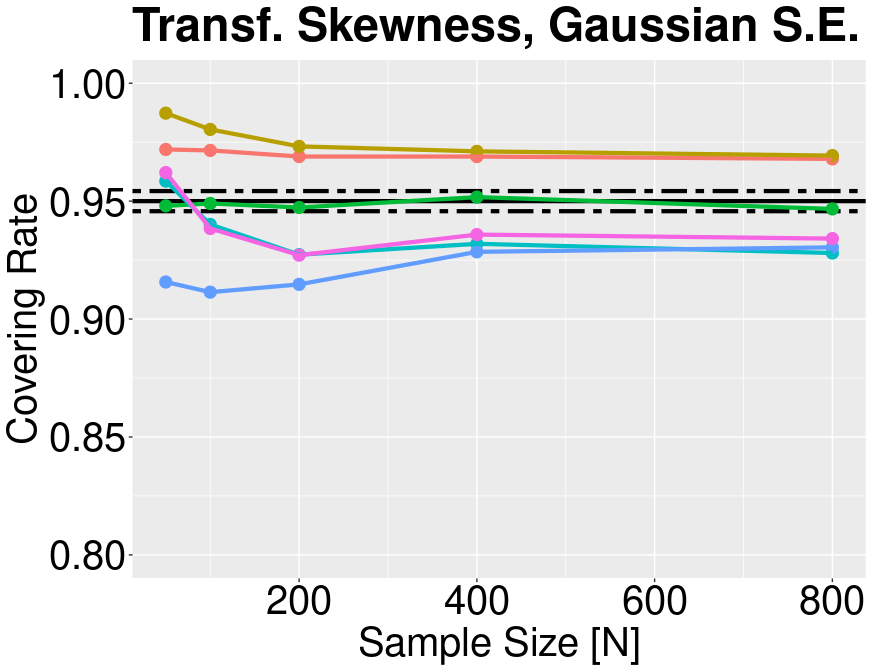}
\includegraphics[trim=0 0 0 0,clip,width=1.5in]{\figurepath 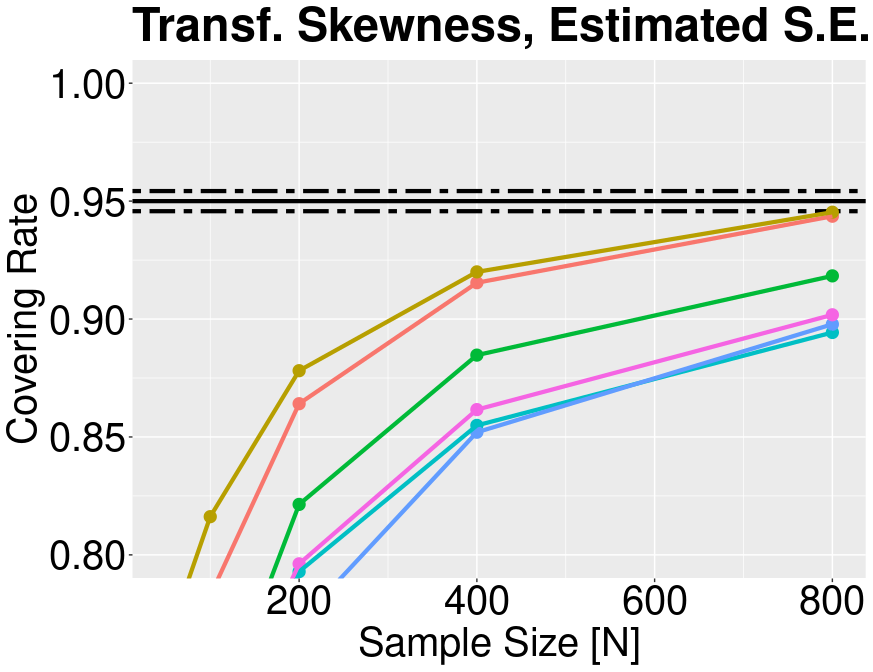}
	\caption{Simulations of coverage rates of SCBs for Model $B$. In these simulations no bias correction is applied.   The black dashed lines are
		95\% confidence intervals for the nominal level $0.95$.  Special attention need to be placed on the third panel. It shows that the test for Gaussianity based on the SCBs for transformed skewness given in eq. \eqref{eq:TestGaussNormalSkew} has the
	correct significance level even for small $N$ for this Gaussian process so long as the Mult method is used for estimation of the quantile.\label{fig:SCBsModelB-skew}}
\end{center}
\end{figure}

Furthermore, in both cases, bias estimation seems to reduce the coverage of the SCBs. Therefore these simulation results are deferred to the appendix. In general, it seems to be better to not account for the bias, which in the Gaussian case is zero for the pointwise skewness. Nevertheless, nominal coverage seems to be still reached for large $ N $ in all scenarios if the bias estimate is included into the construction of the SCBs.
For the non-differentiable Model $ B $, see Figure \ref{fig:SCBsModelB-skew}, the results are comparable to Model $ A $ except that once more the GKF methods have over-coverage for large $ N $, since the sample paths of Model $ B $ are not $ C^2 $ and therefore do not satisfy the assumptions for the GKF methods.

\begin{figure}[h]
\begin{center}
\includegraphics[trim=0 0 0 0,clip,width=1.5in]{\figurepath 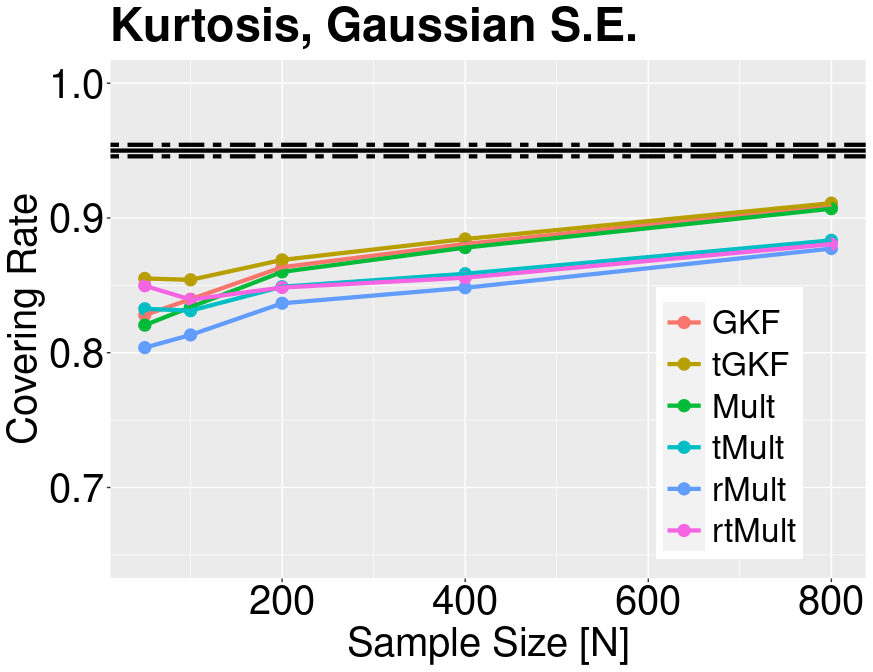}
\includegraphics[trim=0 0 0 0,clip,width=1.5in]{\figurepath 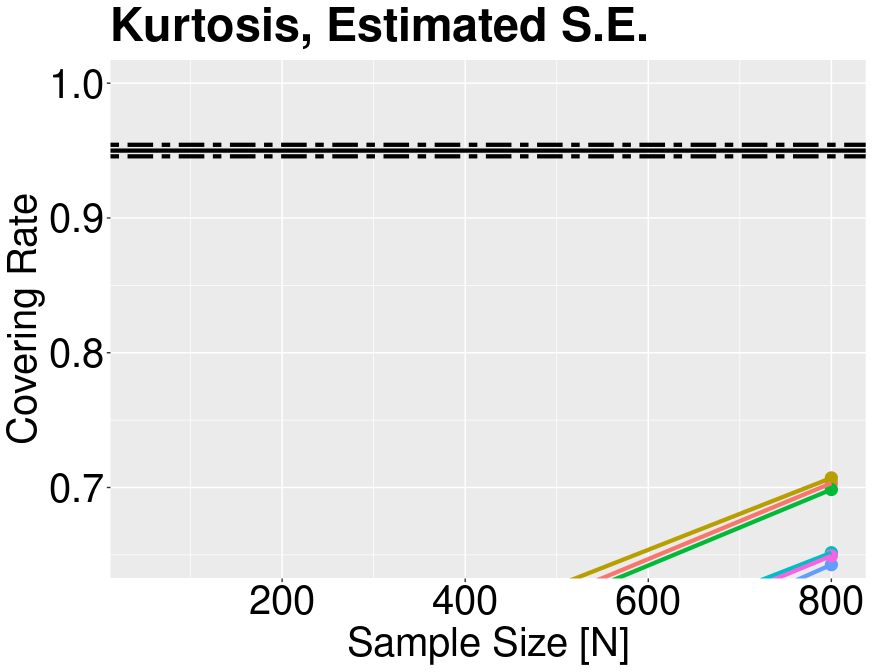}
\includegraphics[trim=0 0 0 0,clip,width=1.5in]{\figurepath 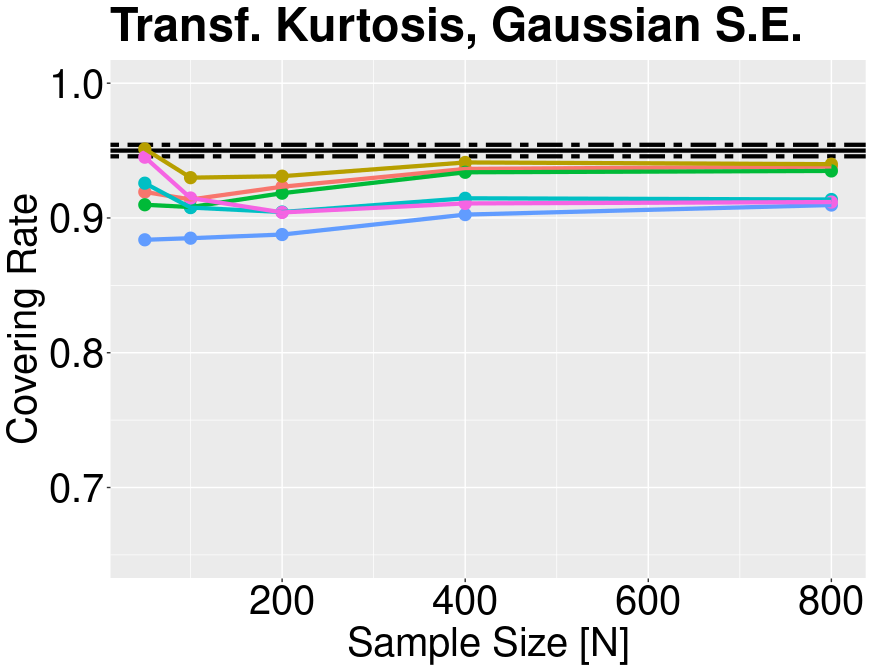}
\includegraphics[trim=0 0 0 0,clip,width=1.5in]{\figurepath 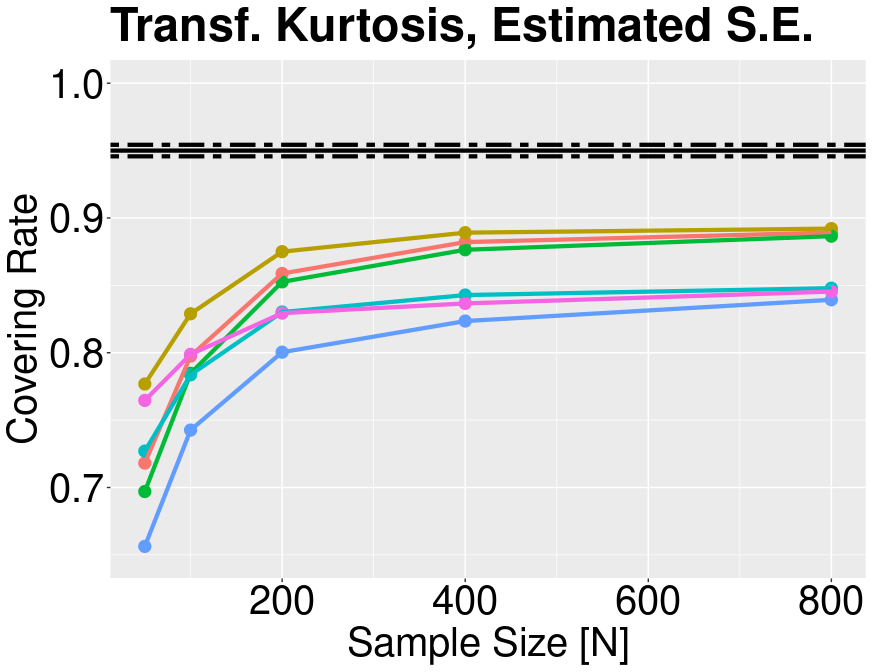}
	\caption{Simulations of coverage rates of SCBs for Model $A$.  In these simulations no bias correction is applied.  The black dashed lines are 95\% confidence intervals for the nominal level $0.95$.  Special attention need to be placed on the third panel. It shows that the test for Gaussianity based on the SCBs for transformed kurtosis given in eq. \eqref{eq:TestGaussNormalKurt} has a significance level that lies close to the nominal even for small $N$ so long as the GKF or Mult method is used for estimation of the quantile.\label{fig:SCBsModelA-kurt}}
\end{center}
\end{figure}

\begin{figure}[h]
\begin{center}
\includegraphics[trim=0 0 0 0,clip,width=1.5in]{\figurepath 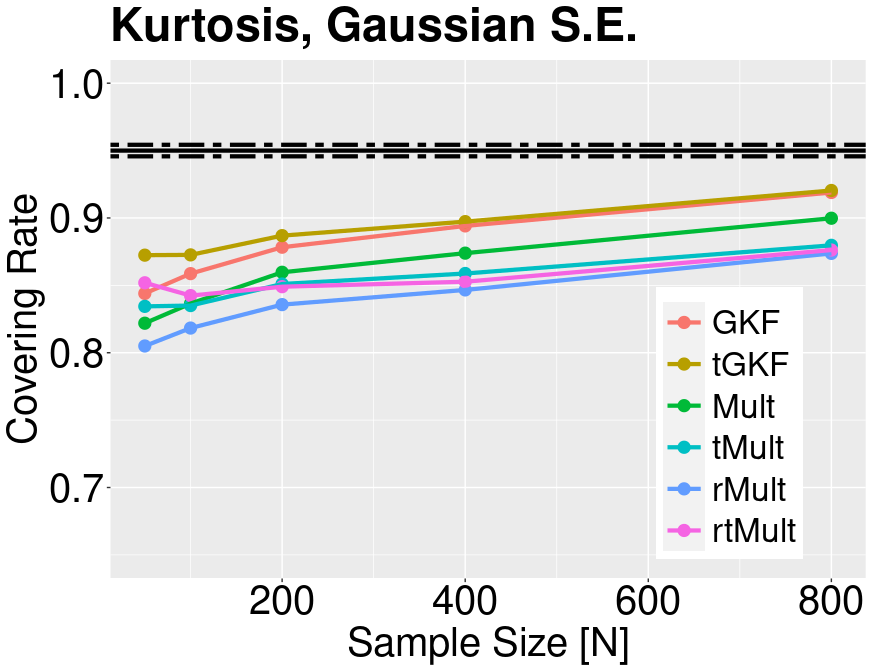}
\includegraphics[trim=0 0 0 0,clip,width=1.5in]{\figurepath 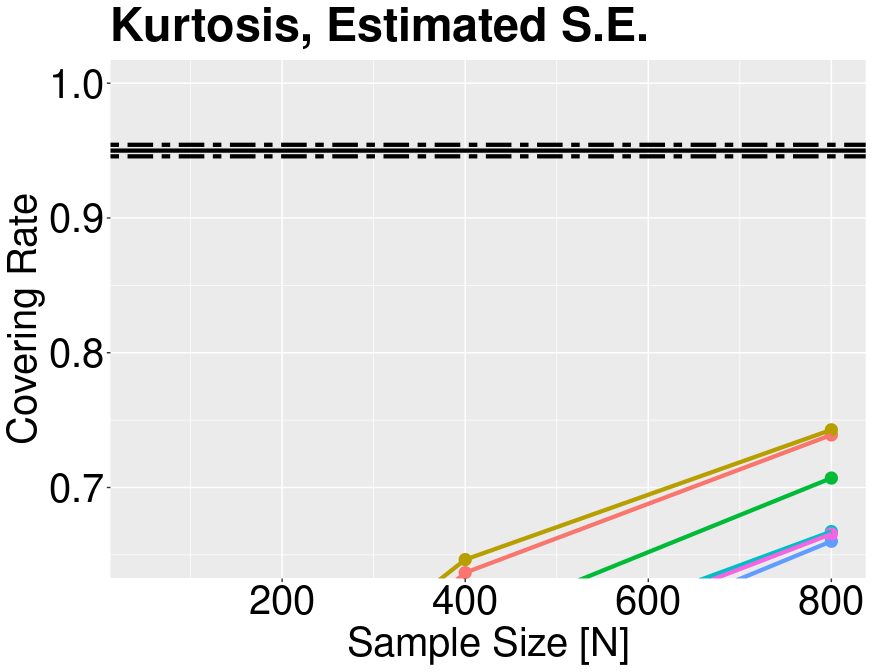}
\includegraphics[trim=0 0 0 0,clip,width=1.5in]{\figurepath 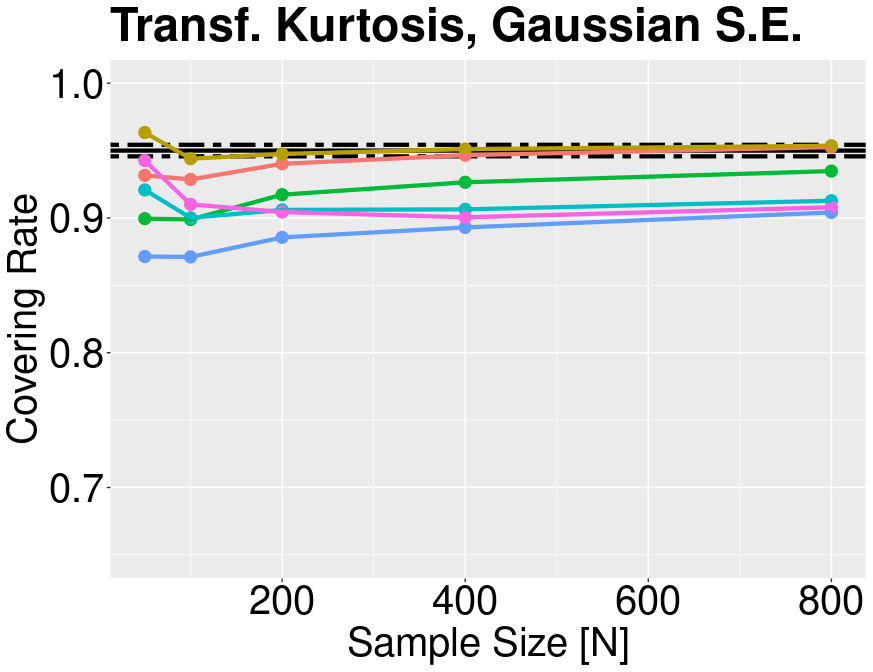}
\includegraphics[trim=0 0 0 0,clip,width=1.5in]{\figurepath 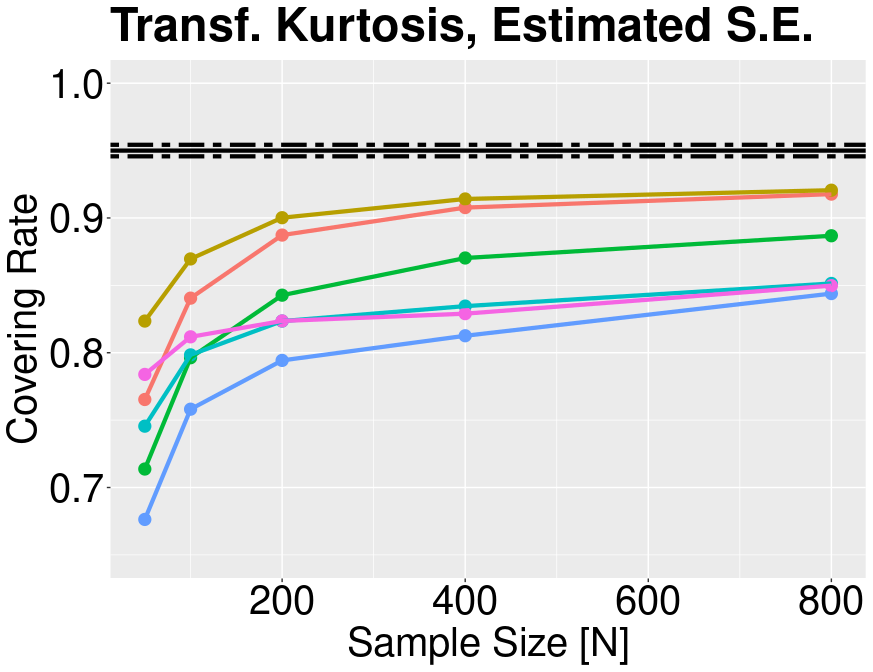}
	\caption{Simulations of coverage rates of SCBs for Model $B$.  In these simulations no bias correction is applied.  The black dashed lines are
		95\% confidence intervals for the nominal level $0.95$.
Special attention need to be placed on the third panel. It shows that the test for Gaussianity based on the SCBs for transformed kurtosis given in eq. \eqref{eq:TestGaussNormalKurt} has a significance level that lies close to the nominal level even for small $N$ as long as the GKF or Mult method is used for estimation of the quantile.
	\label{fig:SCBsModelB-kurt}}
\end{center}
\end{figure}

Simulation results of coverage rates of the SCBs for kurtosis and transformed kurtosis for the Gaussian Models $A$ and $B$ are shown in Figure \ref{fig:SCBsModelA-kurt} and \ref{fig:SCBsModelB-kurt}. They are fairly similar to the results for skewness. The main difference is that
the coverage rates are lower and convergence to nominal requires larger $N$ than in the case of skewness. Nevertheless,
transforming kurtosis is again increasing coverage rates a lot such that even for low sample sizes we always have a coverage rate above $0.9$.
Simulation results for Model $C$ can be found in \ref{app:modelC}.

\FloatBarrier
\subsection{Coverage Rates for Skewness under Sampling and Observational Noise}
\FloatBarrier

\begin{figure}[h]
\begin{center}
	\includegraphics[trim=0 0 0 0,clip,width = 1.5in]{\figurepath 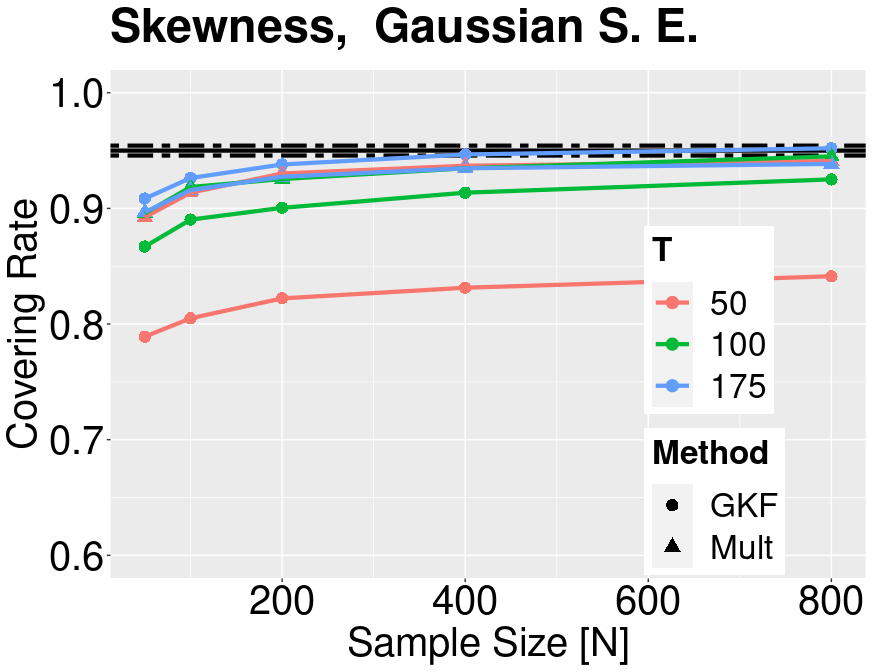}
	\includegraphics[trim=0 0 0 0,clip,width = 1.5in]{\figurepath 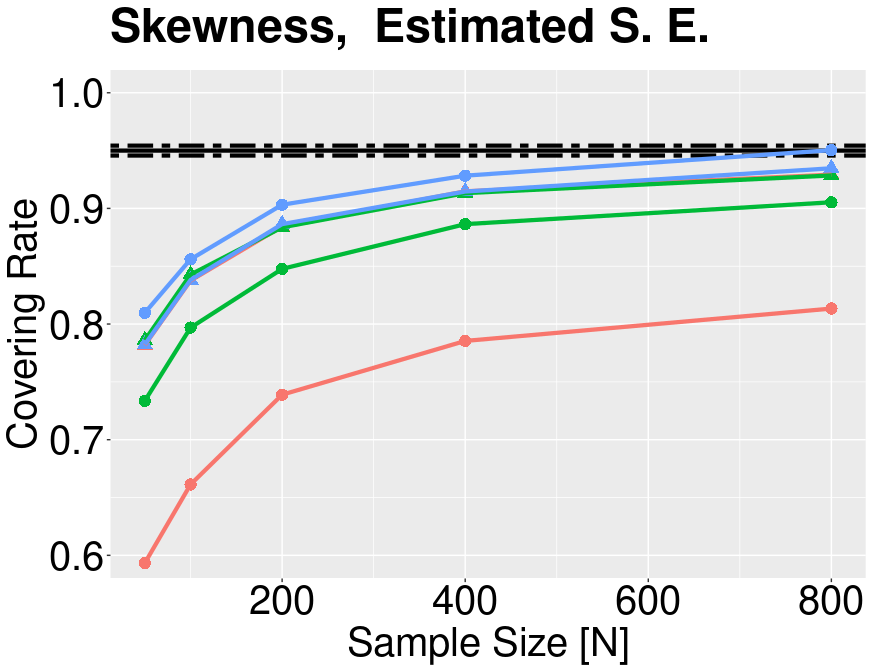}
	\includegraphics[trim=0 0 0 0,clip,width = 1.5in]{\figurepath 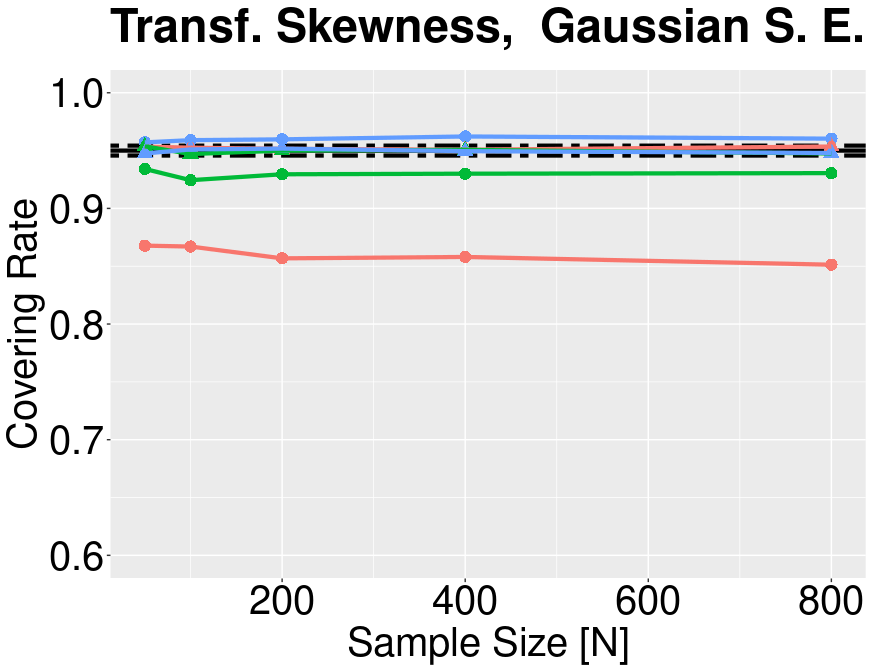}
	\includegraphics[trim=0 0 0 0,clip,width = 1.5in]{\figurepath 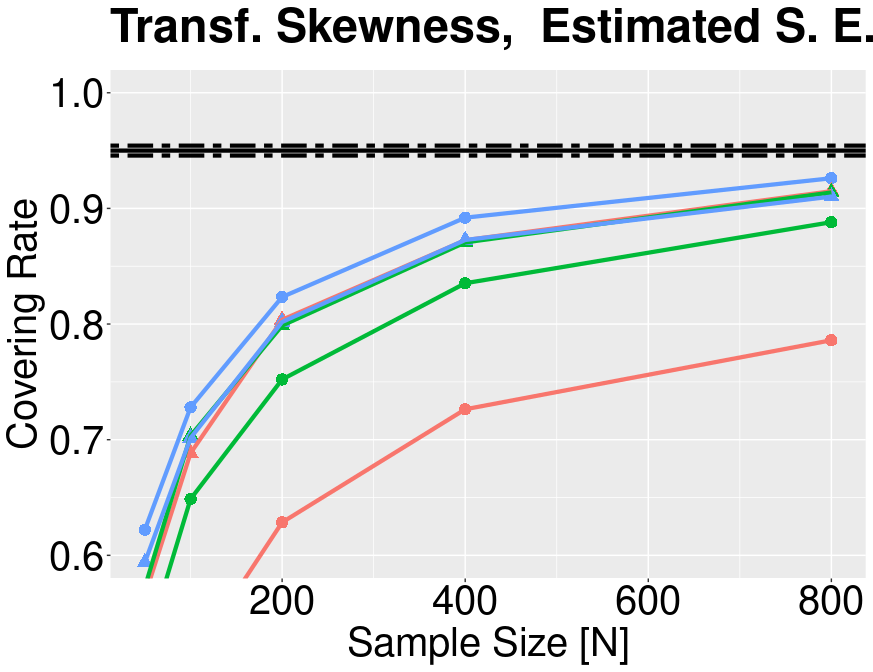}
	\caption{Dependency of coverage rate of SCBs for Model $A$ with observation noise ($\sigma = 0.05$) on the number $T$ of observation points. Panel 3 is again the most important since it shows that even in models with observation noise the significance level for a test for Gaussianity based on transformed skewness is close to the nominal so long as the sampling grid is sufficiently dense.\label{fig:ObsCovRatesA5}}
\end{center}
\end{figure}

\begin{figure}[h]
\begin{center}
	\includegraphics[trim=0 0 0 0,clip,width = 1.5in]{\figurepath 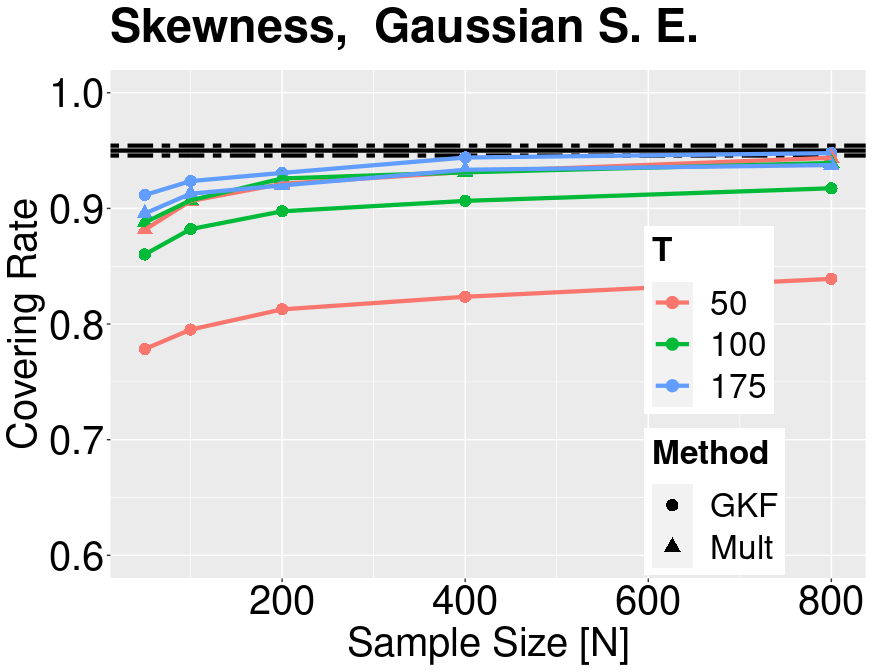}
	\includegraphics[trim=0 0 0 0,clip,width = 1.5in]{\figurepath 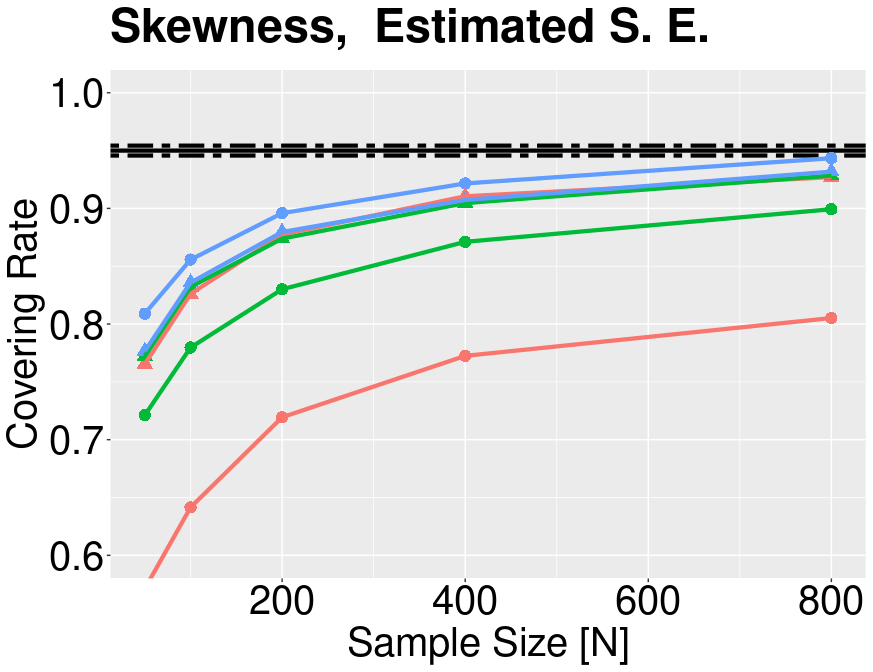}
	\includegraphics[trim=0 0 0 0,clip,width = 1.5in]{\figurepath 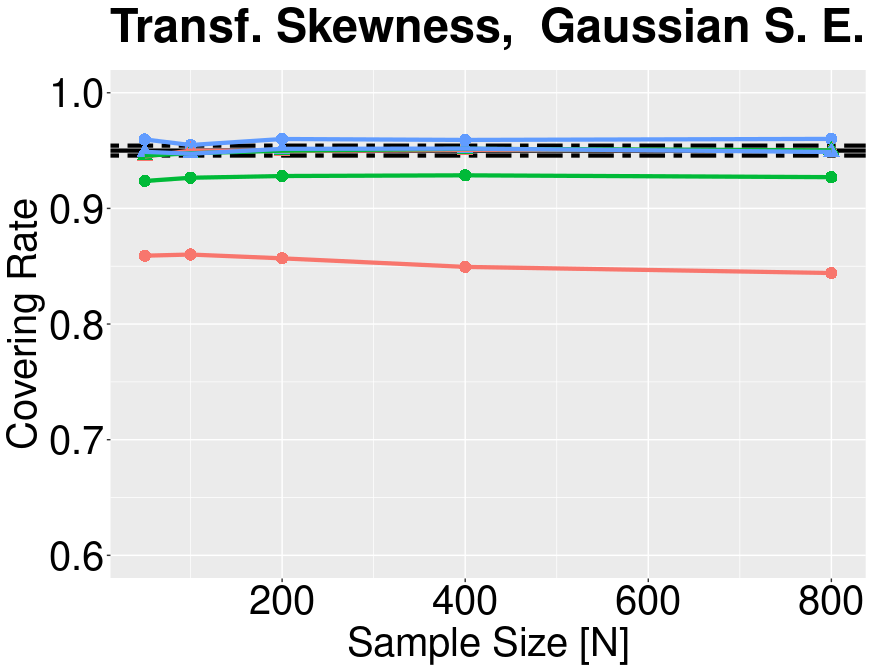}
	\includegraphics[trim=0 0 0 0,clip,width = 1.5in]{\figurepath 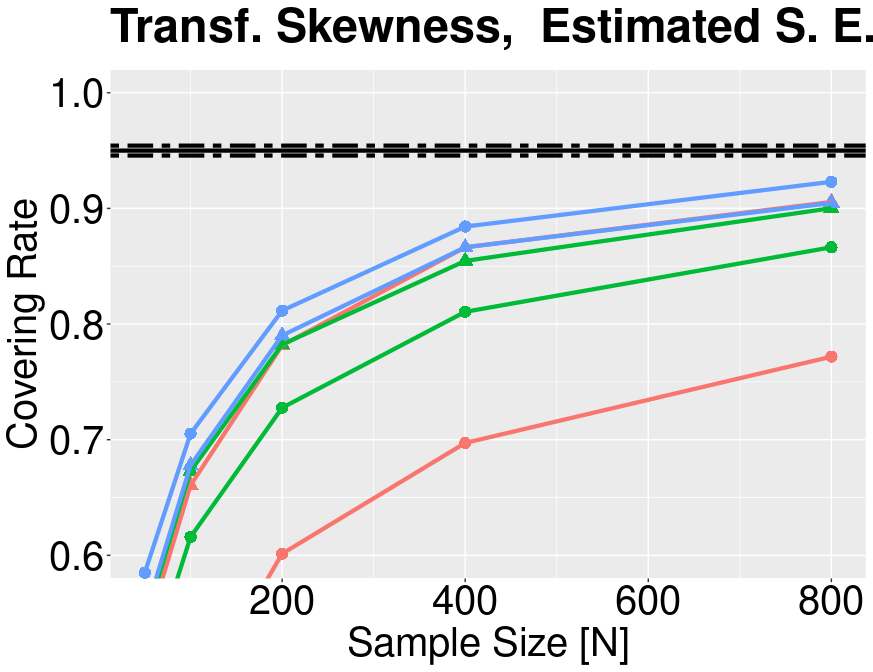}
	\caption{Dependency of coverage rate of SCBs for Model $A$ with observation noise ($\sigma = 0.1$) on the number $T$ of observation points. Panel 3 is again the most important since it shows that even in models with observation noise the significance level for a test for Gaussianity based on transformed skewness is close to the nominal so long as the sampling grid is sufficiently dense.\label{fig:ObsCovRatesA10}}
\end{center}
\end{figure}

\begin{figure}[h]
\begin{center}
	\includegraphics[trim=0 0 0 0,clip,width = 1.5in]{\figurepath 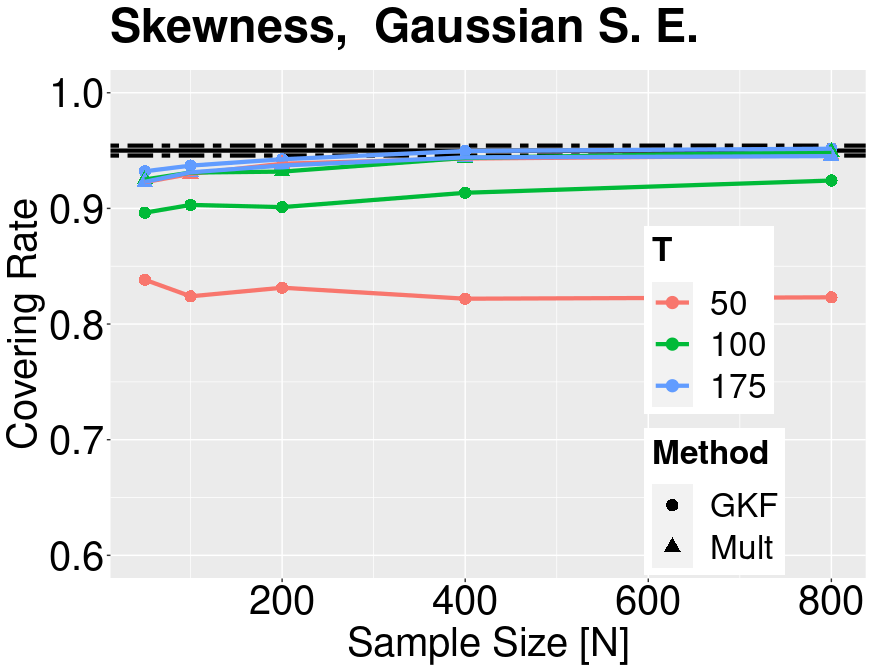}
	\includegraphics[trim=0 0 0 0,clip,width = 1.5in]{\figurepath 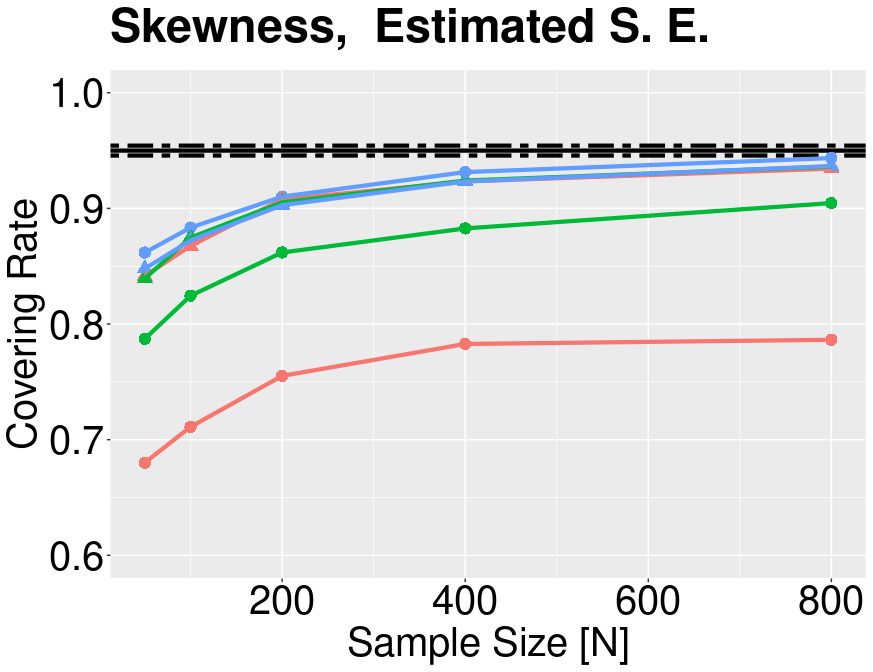}
	\includegraphics[trim=0 0 0 0,clip,width = 1.5in]{\figurepath 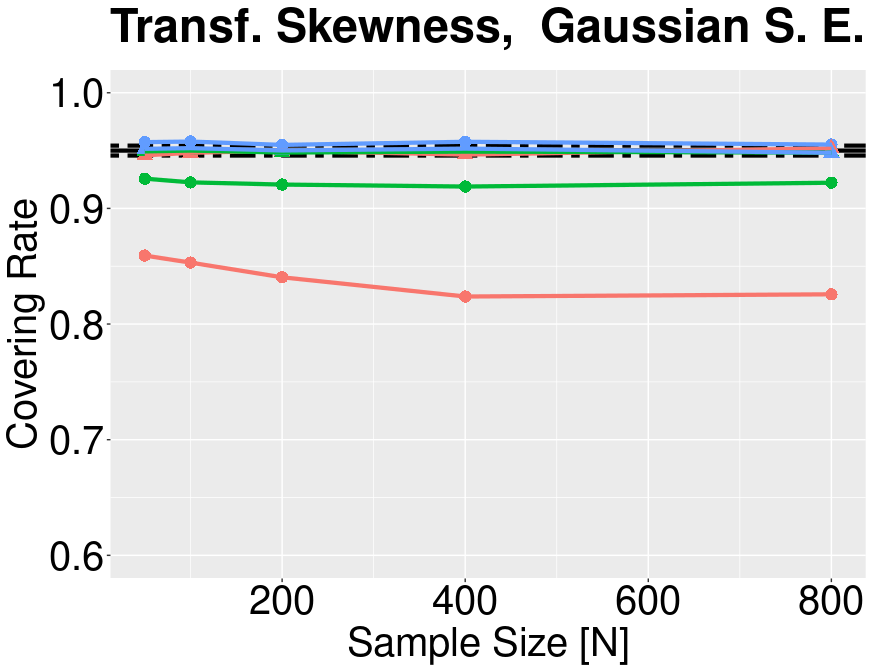}
	\includegraphics[trim=0 0 0 0,clip,width = 1.5in]{\figurepath 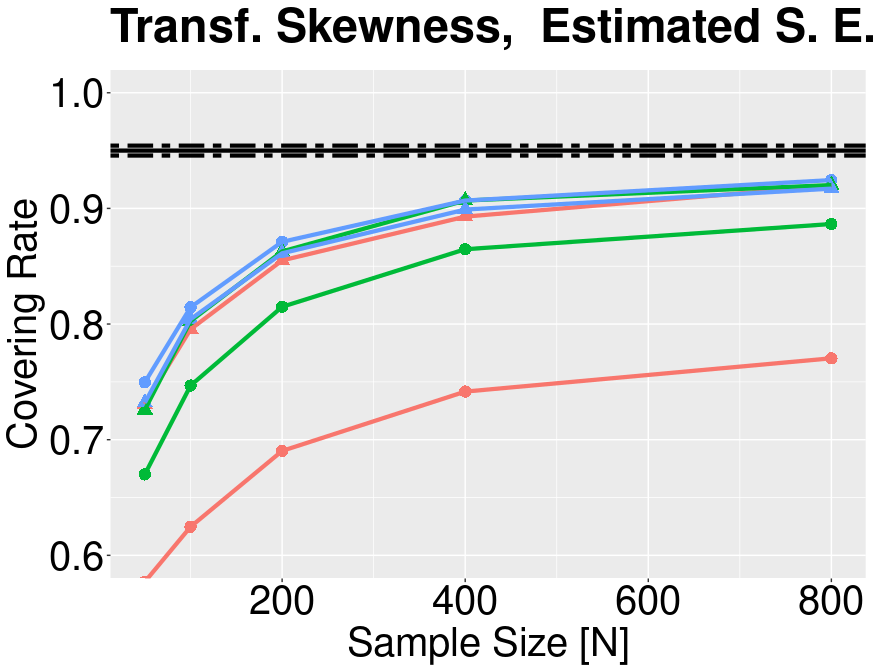}
	\caption{Dependency of coverage rate of SCBs for Model $B$ with observation noise ($\sigma = 0.05$) on the number $T$ of observation points. Panel 3 is again the most important since it shows that even in models with observation noise the significance level for a test for Gaussianity based on transformed skewness is close to the nominal so long as the sampling grid is sufficiently dense.\label{fig:ObsCovRatesB5}}
\end{center}
\end{figure}

\begin{figure}[h]
\begin{center}
	\includegraphics[trim=0 0 0 0,clip,width = 1.5in]{\figurepath 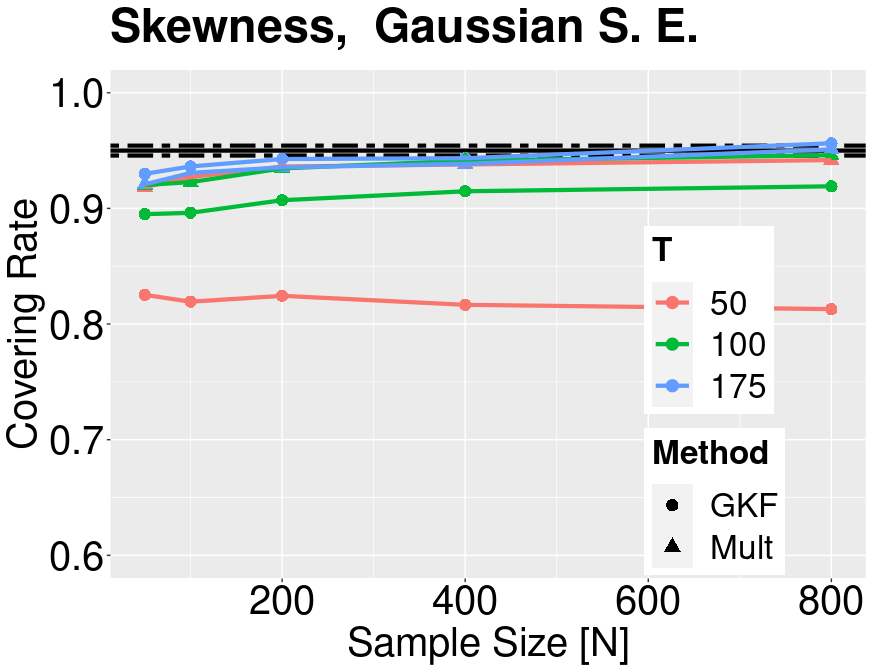}
	\includegraphics[trim=0 0 0 0,clip,width = 1.5in]{\figurepath 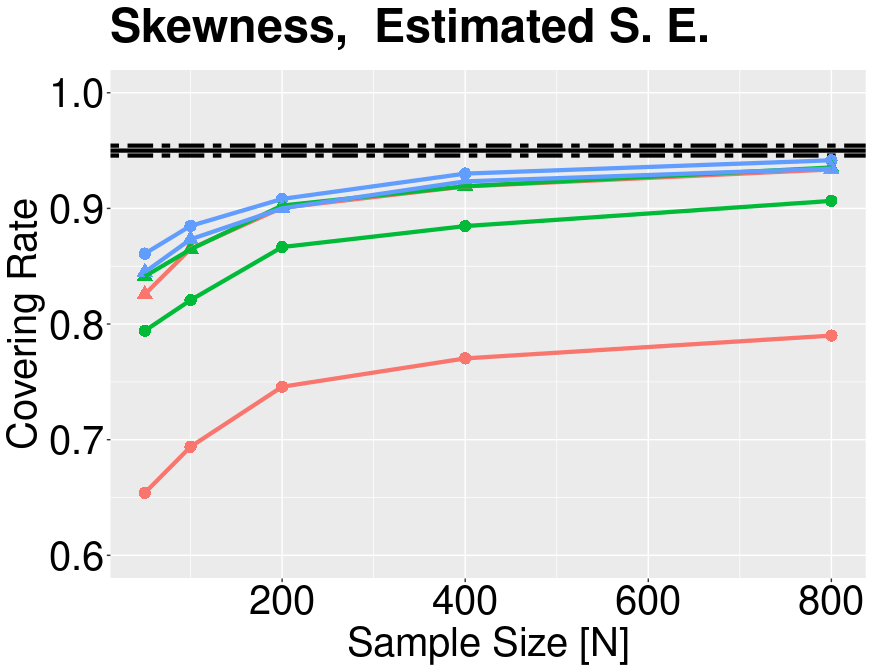}
	\includegraphics[trim=0 0 0 0,clip,width = 1.5in]{\figurepath 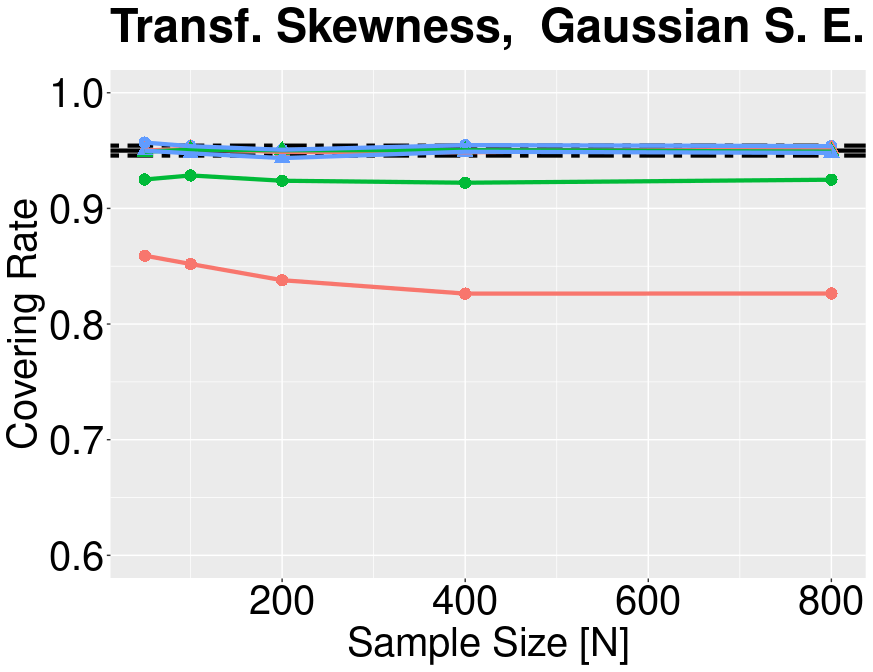}
	\includegraphics[trim=0 0 0 0,clip,width = 1.5in]{\figurepath 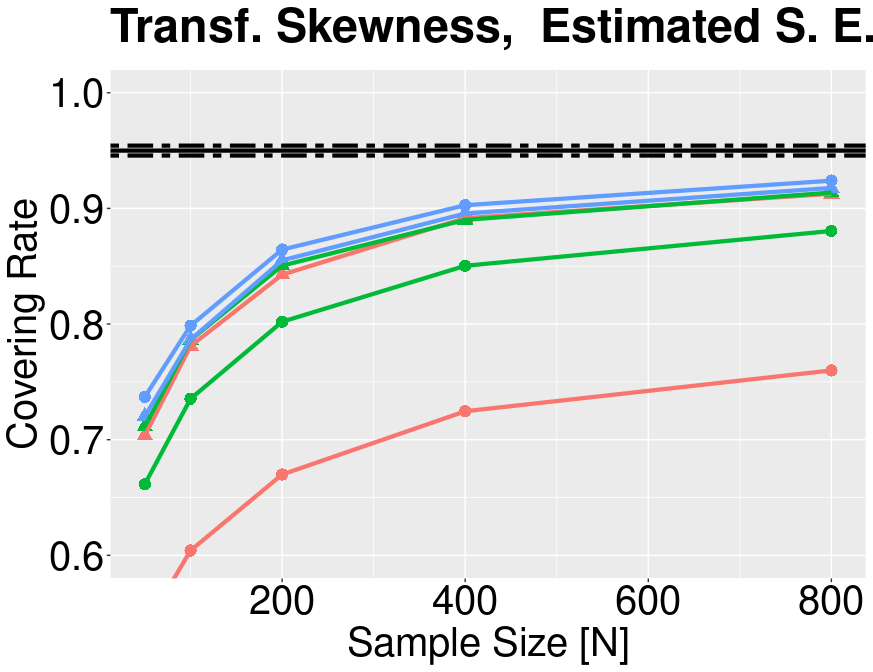}
	\caption{Dependency of coverage rate of SCBs for Model $B$ with observation noise ($\sigma = 0.1$) on the number $T$ of observation points. Panel 3 is again the most important since it shows that even in models with observation noise the significance level for a test for Gaussianity based on transformed skewness is close to the nominal so long as the sampling grid is sufficiently dense.\label{fig:ObsCovRatesB10}}
\end{center}
\end{figure}

In order to demonstrate the dependence of the coverage rates on additional observation noise and the density of the observed grid points, we simulated the coverage rate for the pointwise skewness estimator for $T \in \{ 50, 100, 175 \}$ equidistant grid points of the interval $ [0, 1] $. At each observed grid point we add iid zero-mean Gaussian noise with standard deviation $0.05$ or $0.1$. The observed noisy curves are smoothed using a local linear estimator, e.g., \cite{Degras2011}, where for each observed sample the smoothing bandwidth is obtained by cross-validation using \texttt{cv.select()} from the \texttt{R}-Package \textit{SCBmeanfd} \cite{SCBmeanfd}. SCBs for the skewness parameter using the GKF and gMult are then constructed from the smoothed observations and it is checked whether they contain zero for all $t\in[0,1]$. We decided to only report these two methods for the construction of SCBs, since the other bootstrap methods from the previous simulation perform similarly to before. The simulation results, shown in Figures \ref{fig:ObsCovRatesA5},  \ref{fig:ObsCovRatesA10}.  \ref{fig:ObsCovRatesB5} and  \ref{fig:ObsCovRatesB10}, are similar to the results without observation noise and a dense observation grid from Figure \ref{fig:SCBsModelA-skew} and \ref{fig:SCBsModelB-skew}, if the gMult approach is used. The construction based on the GKF has lower coverage  for $T \in \{50, 100\}$. This is because estimation of the quantile using the tGKF relies on estimation of the variance of the derivative of the error process which is harder when the observations are less dense, since the estimation of these derivatives from the smoothed sample curves has a larger bias.

\begin{figure}[h]
\begin{center}
	\includegraphics[trim=0 0 0 0,clip,width = 1.5in]{\figurepath 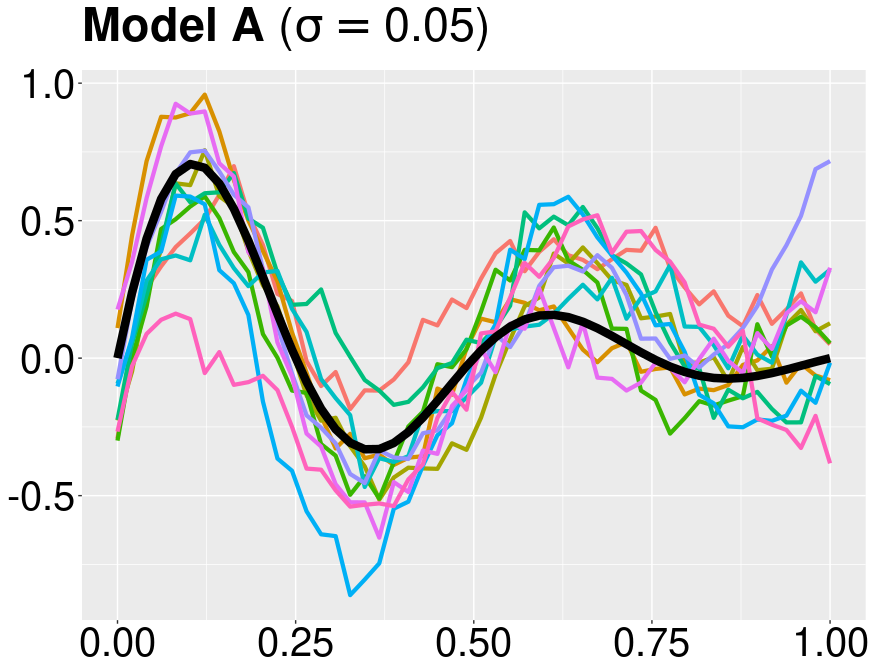}
	\includegraphics[trim=0 0 0 0,clip,width = 1.5in]{\figurepath 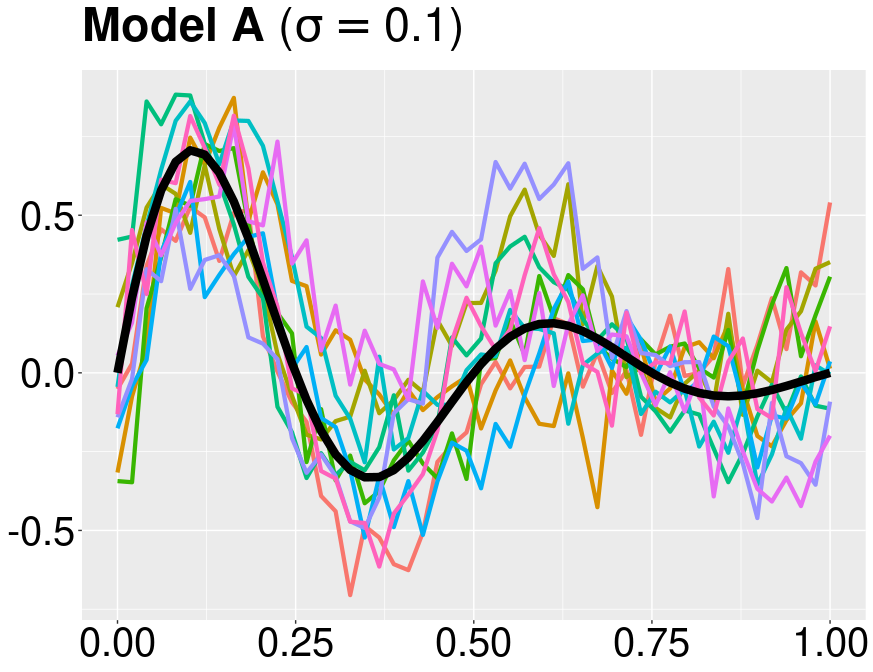}
	\includegraphics[trim=0 0 0 0,clip,width = 1.5in]{\figurepath 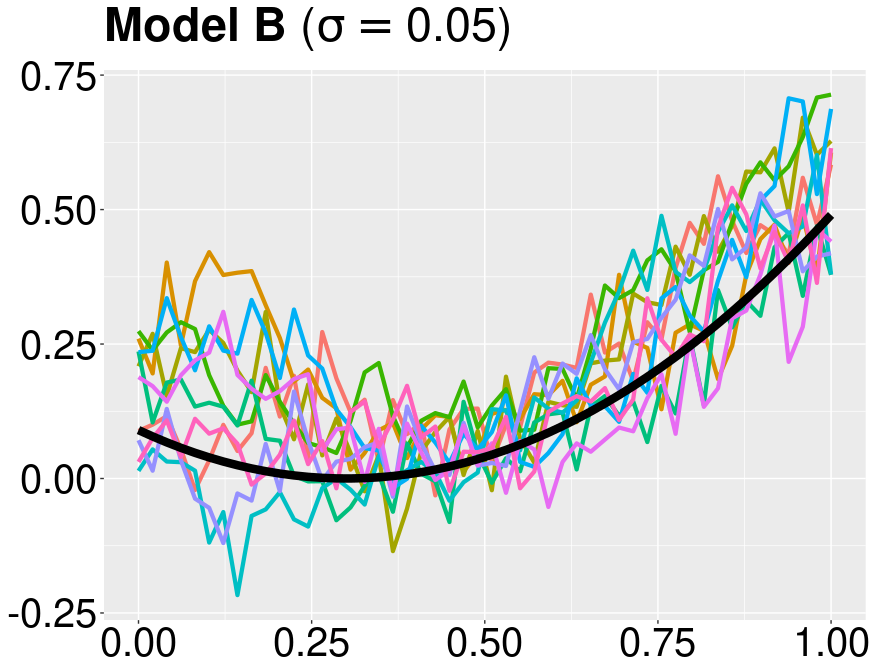}
	\includegraphics[trim=0 0 0 0,clip,width = 1.5in]{\figurepath 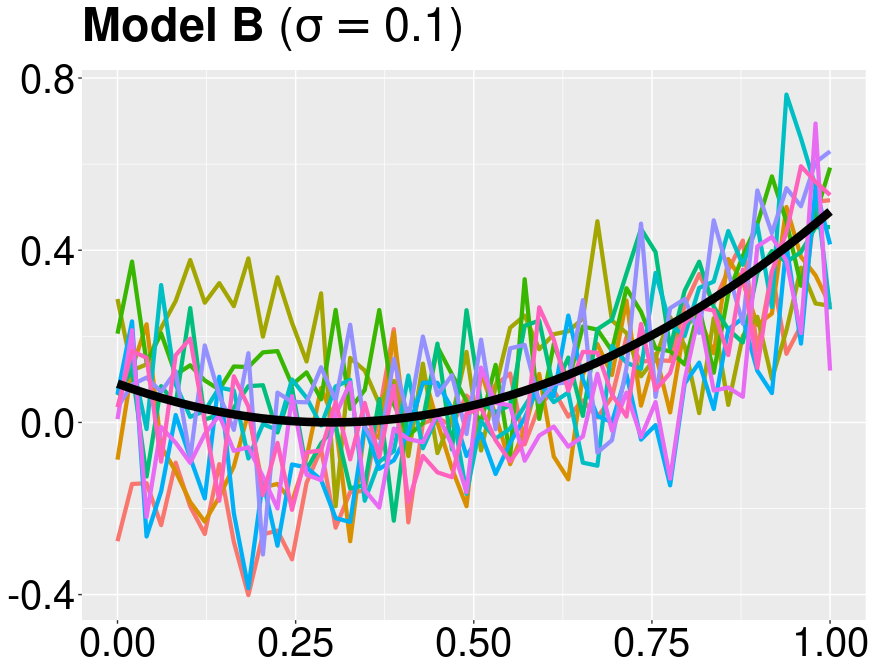}
	\caption{Samples of Model $A$ and Model $B$ with additive Gaussian observation noise. The bold black line is
		the true population mean.\label{fig:ObsCovRates}}
\end{center}
\end{figure} 

\FloatBarrier
\section{Discussion}
\FloatBarrier
Functional delta residuals are a powerful and helpful tool for performing inference on functional data.
They are useful for extending spatial coverage probability excursion sets to effect-size measures,
as demonstrated in \citet{Bowring:2020}, and for the construction of simultaneous confidence bands
for Cohen's $d$ and other moment-based statistics, as done in this article.
Moreover, as we have shown in Section \ref{sec:SCBs} and our simulations such SCBs can be used to test Gaussianity of samples of a random field based on the
transformed skewness or transformed kurtosis statistic.

The main challenge for the future is to improve upon the finite sample coverage for
moment-based statistics. This in particular requires that the estimation of the standard error of
moment-based statistics is improved.
The empirical estimate obtained from the functional delta residuals, which converges asymptotically to the standard error of the moment-based statistic, suffers from the problem that the pointwise distributions of the functional delta residuals are heavily tailed and skewed due to the non-linear transformation. As such more robust estimators for the variance, than the standard sample variance estimator, might be required and could lead to faster convergence to nominal coverage.

Future research could extend the application of the functional delta residuals to other parameter estimators satisfying a fCLT derived from the functional delta method. Potential extensions are coverage probability excursion sets or simultaneous confidence bands for $R^2$ processes from linear models such as those usually fitted in functional magnetic resonance imaging analysis or new statistical methods such as LISA \citep{Lohmann:2018}. The latter introduces spatial smoothing of the z-score field of brain activation together with a false discovery rate controlled inference. In this context, any inference based on random field theory, e.g., cluster inference or coverage probability excursion sets, that is used to detect activation of the smoothed z-score process, will require residuals having the correct asymptotic correlation structure in order to perform valid inference. The functional delta residuals discussed in this article provide a tool to solve this problem.

\section*{Acknowledgments}
F.T., S.D. and A.S. were partially supported by NIH grant R01EB026859.
F.T. is also funded by the Deutsche Forschungsgemeinschaft (DFG) under Excellence Strategy The
Berlin Mathematics Research Center MATH+ (EXC-2046/1, project ID:390685689).

\bibliographystyle{unsrtnat}
\bibliography{\auxpath paper-ref}

\begin{thebibliography}{42}
\providecommand{\natexlab}[1]{#1}
\providecommand{\url}[1]{\texttt{#1}}
\expandafter\ifx\csname urlstyle\endcsname\relax
  \providecommand{\doi}[1]{doi: #1}\else
  \providecommand{\doi}{doi: \begingroup \urlstyle{rm}\Url}\fi

\bibitem[Dette et~al.(2020)Dette, Kokot, Aue, et~al.]{Dette2020functional}
Holger Dette, Kevin Kokot, Alexander Aue, et~al.
\newblock Functional data analysis in the banach space of continuous functions.
\newblock \emph{Annals of Statistics}, 48\penalty0 (2):\penalty0 1168--1192,
  2020.

\bibitem[Dette and Kokot(2020)]{Dette2020bio}
Holger Dette and Kevin Kokot.
\newblock Bio-equivalence tests in functional data by maximum deviation.
\newblock \emph{Biometrika}, 2020.

\bibitem[Adler(1981)]{AdlerBook}
Robert~J Adler.
\newblock \emph{The geometry of random fields}.
\newblock Wiley, 1981.

\bibitem[Worsley et~al.(2004)Worsley, Taylor, Tomaiuolo, and
  Lerch]{Worsley:2004}
Keith~J Worsley, Jonathan~E Taylor, Francesco Tomaiuolo, and Jason Lerch.
\newblock Unified univariate and multivariate random field theory.
\newblock \emph{Neuroimage}, 23:\penalty0 S189--S195, 2004.

\bibitem[Taylor and Worsley(2007)]{Taylor:2007}
Jonathan~E Taylor and Keith~J Worsley.
\newblock Detecting sparse signals in random fields, with an application to
  brain mapping.
\newblock \emph{Journal of the American Statistical Association}, 102\penalty0
  (479):\penalty0 913--928, 2007.

\bibitem[Degras(2011)]{Degras2011}
David~A Degras.
\newblock Simultaneous confidence bands for nonparametric regression with
  functional data.
\newblock \emph{Statistica Sinica}, 21\penalty0 (4):\penalty0 1735--1765, 2011.

\bibitem[Cao et~al.(2012)Cao, Yang, and Todem]{Cao2012}
Guanqun Cao, Lijian Yang, and David Todem.
\newblock Simultaneous inference for the mean function based on dense
  functional data.
\newblock \emph{Journal of Nonparametric Statistics}, 24\penalty0 (2):\penalty0
  359--377, 2012.

\bibitem[Cao(2014)]{Cao2014}
Guanqun Cao.
\newblock Simultaneous confidence bands for derivatives of dependent functional
  data.
\newblock \emph{Electronic Journal of Statistics}, 8\penalty0 (2):\penalty0
  2639--2663, 2014.

\bibitem[Chang et~al.(2017)Chang, Lin, and Ogden]{Chang2017}
Chung Chang, Xuejing Lin, and R~Todd Ogden.
\newblock Simultaneous confidence bands for functional regression models.
\newblock \emph{Journal of Statistical Planning and Inference}, 188:\penalty0
  67--81, 2017.

\bibitem[Wang et~al.(2019)Wang, Wang, Wang, and Ogden]{Wang2019}
Yueying Wang, Guannan Wang, Li~Wang, and R~Todd Ogden.
\newblock Simultaneous confidence corridors for mean functions in functional
  data analysis of imaging data.
\newblock \emph{Biometrics}, 2019.

\bibitem[Telschow and Schwartzman(2022)]{Telschow:2019}
Fabian~JE Telschow and Armin Schwartzman.
\newblock Simultaneous confidence bands for functional data using the gaussian
  kinematic formula.
\newblock \emph{Journal of Statistical Planning and Inference}, 216:\penalty0
  70--94, 2022.

\bibitem[Liebl and Reimherr(2019)]{Liebl2019}
Dominik Liebl and Matthew Reimherr.
\newblock Fast and fair simultaneous confidence bands for functional
  parameters.
\newblock \emph{arXiv preprint arXiv:1910.00131}, 2019.

\bibitem[Cao et~al.(2016)Cao, Wang, Li, and Yang]{Cao2016}
Guanqun Cao, Li~Wang, Yehua Li, and Lijian Yang.
\newblock Oracle-efficient confidence envelopes for covariance functions in
  dense functional data.
\newblock \emph{Statistica Sinica}, pages 359--383, 2016.

\bibitem[Wang et~al.(2020)Wang, Cao, Wang, and Yang]{Wang2020}
Jiangyan Wang, Guanqun Cao, Li~Wang, and Lijian Yang.
\newblock Simultaneous confidence band for stationary covariance function of
  dense functional data.
\newblock \emph{Journal of Multivariate Analysis}, 176:\penalty0 104584, 2020.

\bibitem[Guo et~al.(2018)Guo, Zhou, and Zhang]{Guo2018}
Jia Guo, Bu~Zhou, and Jin-Ting Zhang.
\newblock Testing the equality of several covariance functions for functional
  data: A supremum-norm based test.
\newblock \emph{Computational Statistics \& Data Analysis}, 124:\penalty0
  15--26, 2018.

\bibitem[Sommerfeld et~al.(2018)Sommerfeld, Sain, and
  Schwartzman]{Sommerfeld:2018}
Max Sommerfeld, Stephan Sain, and Armin Schwartzman.
\newblock Confidence regions for spatial excursion sets from repeated random
  field observations, with an application to climate.
\newblock \emph{Journal of the American Statistical Association}, 113\penalty0
  (523):\penalty0 1327--1340, 2018.

\bibitem[Dette and Kokot(2021)]{Dette2020cov}
Holger Dette and Kevin Kokot.
\newblock Detecting relevant differences in the covariance operators of
  functional time series: a sup-norm approach.
\newblock \emph{Annals of the Institute of Statistical Mathematics}, pages
  1--37, 2021.

\bibitem[Bowring et~al.(2019)Bowring, Telschow, Schwartzman, and
  Nichols]{Bowring:2019a}
Alexander Bowring, Fabian Telschow, Armin Schwartzman, and Thomas~E Nichols.
\newblock Spatial confidence sets for raw effect size images.
\newblock \emph{NeuroImage}, 203:\penalty0 116187, 2019.

\bibitem[Chang and Ogden(2009)]{Chang2009}
Chung Chang and R~Todd Ogden.
\newblock Bootstrapping sums of independent but not identically distributed
  continuous processes with applications to functional data.
\newblock \emph{Journal of multivariate analysis}, 100\penalty0 (6):\penalty0
  1291--1303, 2009.

\bibitem[Adler and Taylor(2009)]{Adler:2009}
Robert~J Adler and Jonathan~E Taylor.
\newblock \emph{Random fields and geometry}.
\newblock Springer Science \& Business Media, 2009.

\bibitem[Kosorok(2003)]{Kosorok2003}
Michael~R Kosorok.
\newblock Bootstraps of sums of independent but not identically distributed
  stochastic processes.
\newblock \emph{Journal of Multivariate Analysis}, 84\penalty0 (2):\penalty0
  299--318, 2003.

\bibitem[Bowring et~al.(2021)Bowring, Telschow, Schwartzman, and
  Nichols]{Bowring:2020}
Alexander Bowring, Fabian~JE Telschow, Armin Schwartzman, and Thomas~E Nichols.
\newblock Confidence sets for cohen’sd effect size images.
\newblock \emph{NeuroImage}, 226:\penalty0 117477, 2021.

\bibitem[Davenport and Nichols(2020)]{Davenport:2020}
Samuel Davenport and Thomas~E Nichols.
\newblock Selective peak inference: Unbiased estimation of raw and standardized
  effect size at local maxima.
\newblock \emph{Neuroimage}, 209:\penalty0 116375, 2020.

\bibitem[D'agostino et~al.(1990)D'agostino, Belanger, and
  D'Agostino~Jr]{DAgostino:1990}
Ralph~B D'agostino, Albert Belanger, and Ralph~B D'Agostino~Jr.
\newblock A suggestion for using powerful and informative tests of normality.
\newblock \emph{The American Statistician}, 44\penalty0 (4):\penalty0 316--321,
  1990.

\bibitem[Kosorok(2008)]{Kosorok:2008}
Michael~R Kosorok.
\newblock \emph{Introduction to empirical processes and semiparametric
  inference.}
\newblock Springer, 2008.

\bibitem[Jain and Marcus(1975)]{Jain:1975}
Naresh~C Jain and Michael~B Marcus.
\newblock Central limit theorems for $c(s)$-valued random variables.
\newblock \emph{Journal of Functional Analysis}, 19\penalty0 (3):\penalty0
  216--231, 1975.

\bibitem[Landau and Shepp(1970)]{Landau:1970}
Henry~J Landau and Lawrence~A Shepp.
\newblock On the supremum of a gaussian process.
\newblock \emph{Sankhy{\=a}: The Indian Journal of Statistics, Series A}, pages
  369--378, 1970.

\bibitem[Davidson(1994)]{Davidson1994}
James Davidson.
\newblock \emph{Stochastic Limit Theory: An introduction for econometricians}.
\newblock OUP Oxford, 1994.

\bibitem[Van Der~Vaart and Wellner(1996)]{Vaart1996}
Aad~W Van Der~Vaart and Jon~A Wellner.
\newblock \emph{Weak convergence and empirical processes}.
\newblock Springer, 1996.

\bibitem[Billingsley(1999)]{Billingsley:1999}
Patrick Billingsley.
\newblock \emph{Convergence of probability measures}.
\newblock John Wiley \& Sons, 1999.

\bibitem[Vandekar and Stephens(2021)]{Vandekar2021improving}
Simon~N Vandekar and Jeremy Stephens.
\newblock Improving the replicability of neuroimaging findings by thresholding
  effect sizes instead of p-values.
\newblock \emph{Human Brain Mapping}, 2021.

\bibitem[Vignat(2012)]{Vignat:2012}
Christophe Vignat.
\newblock A generalized isserlis theorem for location mixtures of gaussian
  random vectors.
\newblock \emph{Statistics \& Probability Letters}, 82\penalty0 (1):\penalty0
  67--71, 2012.

\bibitem[Joanes and Gill(1998)]{Joanes:1998}
Derrick~N Joanes and Christine~A Gill.
\newblock Comparing measures of sample skewness and kurtosis.
\newblock \emph{Journal of the Royal Statistical Society: Series D (The
  Statistician)}, 47\penalty0 (1):\penalty0 183--189, 1998.

\bibitem[Fisher(1930)]{Fisher:1930}
Ronald~A Fisher.
\newblock The moments of the distribution for normal samples of measures of
  departure from normality.
\newblock \emph{Proceedings of the Royal Society of London. Series A,
  Containing Papers of a Mathematical and Physical Character}, 130\penalty0
  (812):\penalty0 16--28, 1930.

\bibitem[Taylor et~al.(2005)Taylor, Takemura, and Adler]{Taylor2005}
Jonathan Taylor, Akimichi Takemura, and Robert~J Adler.
\newblock Validity of the expected euler characteristic heuristic.
\newblock \emph{The Annals of Probability}, 33\penalty0 (4):\penalty0
  1362--1396, 2005.

\bibitem[Telschow et~al.(2020)Telschow, Schwartzman, Cheng, and
  Pranav]{Schwartzman:2019}
Fabian Telschow, Armin Schwartzman, Dan Cheng, and Pratyush Pranav.
\newblock Estimation of expected euler characteristic curves of nonstationary
  smooth gaussian random fields.
\newblock \emph{arXiv preprint arXiv:1908.02493}, 2020.

\bibitem[Lehmann et~al.(2005)Lehmann, Romano, and Casella]{Lehmann2005}
Erich~Leo Lehmann, Joseph~P Romano, and George Casella.
\newblock \emph{Testing statistical hypotheses}, volume~3.
\newblock Springer, 2005.

\bibitem[D'Agostino(1970)]{DAgostino:1970}
Ralph~B D'Agostino.
\newblock Transformation to normality of the null distribution of g1.
\newblock \emph{Biometrika}, pages 679--681, 1970.

\bibitem[Anscombe and Glynn(1983)]{Anscombe:1983}
Francis~J Anscombe and William~J Glynn.
\newblock Distribution of the kurtosis statistic b 2 for normal samples.
\newblock \emph{Biometrika}, 70\penalty0 (1):\penalty0 227--234, 1983.

\bibitem[Degras(2016)]{SCBmeanfd}
David Degras.
\newblock \emph{SCBmeanfd: Simultaneous Confidence Bands for the Mean of
  Functional Data}, 2016.
\newblock URL \url{https://CRAN.R-project.org/package=SCBmeanfd}.
\newblock R package version 1.2.2.

\bibitem[Lohmann et~al.(2018)Lohmann, Stelzer, Lacosse, Kumar, Mueller, Kuehn,
  Grodd, and Scheffler]{Lohmann:2018}
Gabriele Lohmann, Johannes Stelzer, Eric Lacosse, Vinod~J Kumar, Karsten
  Mueller, Esther Kuehn, Wolfgang Grodd, and Klaus Scheffler.
\newblock Lisa improves statistical analysis for fmri.
\newblock \emph{Nature communications}, 9\penalty0 (1):\penalty0 4014, 2018.

\bibitem[Lang(1993)]{Lang:1993}
Serge Lang.
\newblock \emph{Real and functional analysis}, volume 142.
\newblock Springer Science \& Business Media, 1993.

\end{thebibliography}

\section{Auxiliary Lemmas}\label{app:auxlemmas}
We now prove a small extension of the functional delta method (for example, \cite[Theorem 2.8, p.235]{Kosorok:2008})
such that the transformation $H$ can be a sequence of transformations $H_N$. Afterwards we will show that the transformations
$Z_{1,N}$ and $Z_{2,N}$ satisfy the assumptions of the extended functional delta method.
\begin{lemma}
	Assume that $ (H_N)_{N\in \mathbb{N}} \subset C^1\big( \mathbb{R}^D, \mathbb{R} \big) $
	be a sequence and $H^1\in C\big( \mathbb{R}^D, \mathbb{R} \big)$. Let $\nabla H\vert_x$, $\nabla H_N\vert_x$
	denote the gradients at $x\in \mathbb{R}^D$. For $\theta \in C\big(S, \mathbb{R}^D\big)$ define the set
	\begin{equation}
		D^{\epsilon}_\theta = \Big\{\, x\in \mathbb{R}^D:~x = \theta(s) \pm \eta ~\text{ for some }s\in S\,,~\vert \eta\vert \leq \epsilon \,\Big\}
	\end{equation}
	and assume that
	\begin{equation}\label{UnifCond}
		\lim_{N\rightarrow\infty} \max_{x \in D^{\epsilon}_\theta} \big\vert\, \nabla H_N\vert_x - \nabla H\vert_x \,\big\vert = 0\,.
	\end{equation}
	Let $r_N\rightarrow \infty$ and $(X_N)_{N\in \mathbb{N}} \subset C\big( S, \mathbb{R}^D \big)$ be a sequence of random
	fields such that
	\begin{equation}
		r_N(\, X_N - \theta \,) \rightsquigarrow X\,,~ ~\text{ in } C\big( S, \mathbb{R}^D \big)
	\end{equation}
	with $X \in C\big( S, \mathbb{R}^D \big)$ being a limiting random field. Then
	\begin{equation}
		r_N\big(\, H_N(X_N) - H_N(\theta) \,\big) \rightsquigarrow \nabla H\vert_\theta X\,,~ ~\text{ in } C\big( S, \mathbb{R}^D \big)
	\end{equation}
\end{lemma}
\begin{proof}
	We follow the idea of the proof of the functional delta method from \cite[Theorem 2.8, p.235]{Kosorok:2008}.
	The key step is establishing
	\begin{equation}
		r_N\big(\, H_N(\theta + r_N^{-1}h_N) - H_N(\theta) \,\big) \rightarrow \nabla H\vert_\theta h
	\end{equation}
	for all $h_N \in C\big( S, \mathbb{R}^D \big)$ satisfying $h_N \rightarrow h$ in $C\big( S, \mathbb{R}^D \big)$.
	Using Taylor's theorem we obtain
	\begin{equation}
	\begin{split}
		r_N\Big(\, H_N(\theta + r_N^{-1}h_N) - H_N(\theta) \,\Big)
		= \nabla H_N\vert_\theta h_N + \int_0^1 \left(\nabla H_N\vert_{\theta + yr_N^{-1}h_N} - \nabla H_N\vert_\theta\right) h_N dy\,,
	\end{split}
	\end{equation}
	see for example \cite[p. 349f]{Lang:1993}. First note that
	\begin{equation}
	\begin{split}
		\Big\Vert\, \nabla H_N\vert_\theta h_N - \nabla H\vert_\theta h \,\Big\Vert_\infty
		\leq \Big\Vert\, \nabla H_N\vert_\theta \Big\Vert_\infty \Big\Vert  h_N -  h  \,\Big\Vert_\infty
		+ \Big\Vert\, \nabla H_N\vert_\theta - \nabla H\vert_\theta \,\Big\Vert_\infty \Vert  h\Vert_\infty
	\end{split}
	\end{equation}	
	converges to zero by \eqref{UnifCond}. To see this note that $\big\Vert \nabla H_N\vert_\theta \big\Vert_\infty \rightarrow \big\Vert \nabla H\vert_\theta \big\Vert_\infty $ by inverse triangle inequality.
	For the second term we have for $N$ large enough that
	\begin{equation}
	\begin{split}
		 \frac{1}{\Vert h_N \Vert_\infty}\max_{s\in S}&\left\vert \int_0^1 \left(\nabla H_N\vert_{\theta(s) + yr_N^{-1}h_N(s)} - \nabla H_N\vert_{\theta(s)}\right) h_N(s) dy
		\right\vert\\
		&\hspace{1cm}\leq
		\max_{s\in S}\max_{y\in[0,1]}\left\vert \nabla H_N\vert_{\theta(s) + yr_N^{-1}h_N(s)} - \nabla H_N\vert_{\theta(s)}
		\right\vert\\
		&\hspace{1cm}\leq
		\max_{s\in S}\max_{y\in[0,1]}\left\vert \nabla H_N\vert_{\theta(s) + yr_N^{-1}h_N(s)} - \nabla H\vert_{\theta(s) + yr_N^{-1}h_N(s)}
		\right\vert\\
		&\hspace{1.4cm}+
		\max_{s\in S}\max_{y\in[0,1]}\left\vert \nabla H\vert_{\theta(s) + yr_N^{-1}h_N(s)} - \nabla H\vert_{\theta(s)}
		\right\vert\\
		&\hspace{1.4cm}+
		\max_{s\in S}\max_{y\in[0,1]}\Big\vert \nabla H\vert_{\theta(s)} - \nabla H_N\vert_{\theta(s)}
		\Big\vert\\
		&\hspace{1cm}\leq
		2\max_{x \in D^{\epsilon}_\theta}\Big\vert \nabla H_N\vert_{x} - \nabla H\vert_{x}
		\Big\vert
		+ \max_{x\in  D^{\epsilon}_\theta}\max_{\vert y \vert < \epsilon}\Big\vert \nabla H\vert_{x} - \nabla H\vert_{x + y}
		\Big\vert
	\end{split}
	\end{equation}
	The first summand converges to zero by condition \eqref{UnifCond}, the second can be made arbitrarily small, since
	as a continuous function $\nabla H\vert_{x}$ is uniform continuous on the compact set $ D^{\epsilon}_\theta$.
	The rest of the proof follows Kosorok's proof.
\end{proof}
\begin{remark}
	It can be easily seen from the proof that \eqref{UnifCond} can be replaced with $ (H_N)_{N\in \mathbb{N}} \subset C^2\big( \mathbb{R}^D, \mathbb{R} \big) $ and
	\begin{equation}
	\begin{split}
		\lim_{N\rightarrow\infty} \max_{s \in S} \big\vert\, \nabla H_N\vert_\theta - \nabla H\vert_\theta \,\big\vert = 0\,~ ~ ~ ~ ~ ~\text{ and },~ ~ ~ ~ ~ ~
		\sup_{N \in \mathbb{N}} \max_{x \in D^{\epsilon}_\theta} \big\vert\, \nabla^2 H_{x} \,\big\vert < \infty\,.
	\end{split}
	\end{equation}
\end{remark}

It remains to show that $Z_{1,N}$ and $Z_{2,N}$ satisfy the uniform convergence condition. Note that by the chain rule it is enough to show the assumptions only for these functions and not for the
whole transformation which maps the moments to $Z_{1,N}(g_1)$, $Z_{2,N}(g_2)$ respectively.
\begin{lemma}
	Let $(a_N)_{N\in\mathbb{N}}$, $(b_N)_{N\in\mathbb{N}}$, $(c_N)_{N\in\mathbb{N}} \subset\mathbb{R}$ be sequences converging to $a,b,c\in \mathbb{R}$. Let $H\in C^1(\mathbb{R})$. Define the sequence of transformations to be $H_N(x) = a_NH(b_Nx + c_N)$. Then for any compact set $D\subset\mathbb{R}$ it follows that
	\begin{equation}
		\lim_{N\rightarrow\infty} \max_{x \in D} \big\vert\,H_N'\vert_x - aH'\vert_{bx+c} \,\big\vert = 0\,.
	\end{equation}
\end{lemma}
\begin{proof}
	We compute:
	\begin{equation}
	\begin{split}
		\big\vert\,H_N'\vert_x - aH'\vert_{bx+c} \,\big\vert
		&= \big\vert\,a_NH'\vert_{b_Nx+c_N} - aH'\vert_{bx+c} \,\big\vert\\
		&\leq \vert\,a_N - a \,\vert \big\vert H'\vert_{b_Nx+c_N}\,\big\vert + \vert a\vert \big\vert\,H'\vert_{b_Nx+c_N} - H'\vert_{bx+c} \,\big\vert\,.
	\end{split}
	\end{equation}
	Both summands converge to zero uniformly for all $x\in D$, since $\big\vert H'\vert_{b_Nx+c_N}\,\big\vert$ is bounded
	on the compact set $b_N D c_N$ and since
	$\vert b_Nx+c_N - bx+c\vert \leq \vert b_N -b \vert \max_{x\in D} \vert x\vert + \vert c_N - c \vert $ goes to zero
	 uniformly and $H'\vert_x$ is uniform continuous on compact sets.
\end{proof}
 
\section{Bias corrected SCBs for Skewness and Kurtosis in Model A and Model B}

In this section we report additional simulation results showing that bias correction in the construction of the SCBs for (transformed) skewness and (transformed) kurtosis in general does not improve the coverage rate of the SCBs, but rather decreases it. For the pointwise skewness estimator this is clear a priori for Gaussian data since the skewness estimator is unbiased. Therefore estimation of the bias only increases the variance. For transformed skewness, kurtosis and transformed
kurtosis it is not immediately clear whether inclusion of the bias estimate is preferable. However, our simulations
suggest that not including the bias estimate seems to be preferable. The other qualitative aspects of these simulations are
very similar to the discussion of the SCBs from the simulation section of the main article.

\begin{figure}[h]
\begin{center}
\includegraphics[trim=0 0 0 0,clip,width=1.5in]{\figurepath 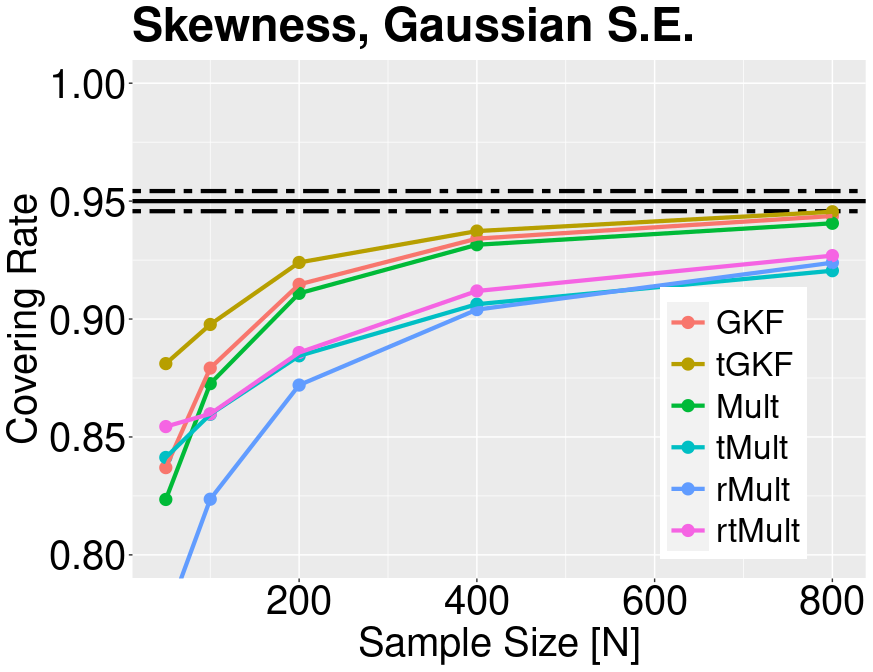}
\includegraphics[trim=0 0 0 0,clip,width=1.5in]{\figurepath 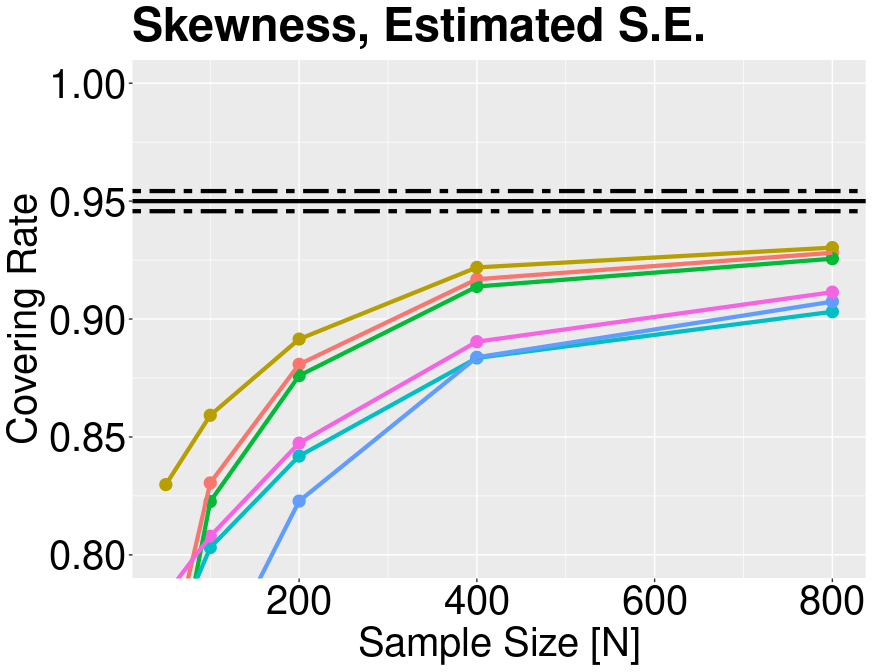}
\includegraphics[trim=0 0 0 0,clip,width=1.5in]{\figurepath 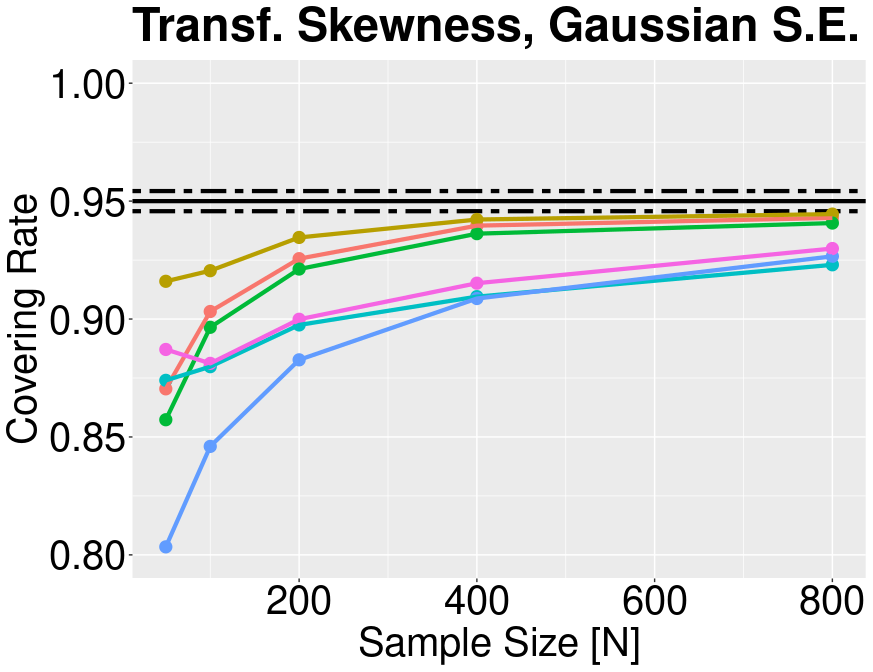}
\includegraphics[trim=0 0 0 0,clip,width=1.5in]{\figurepath 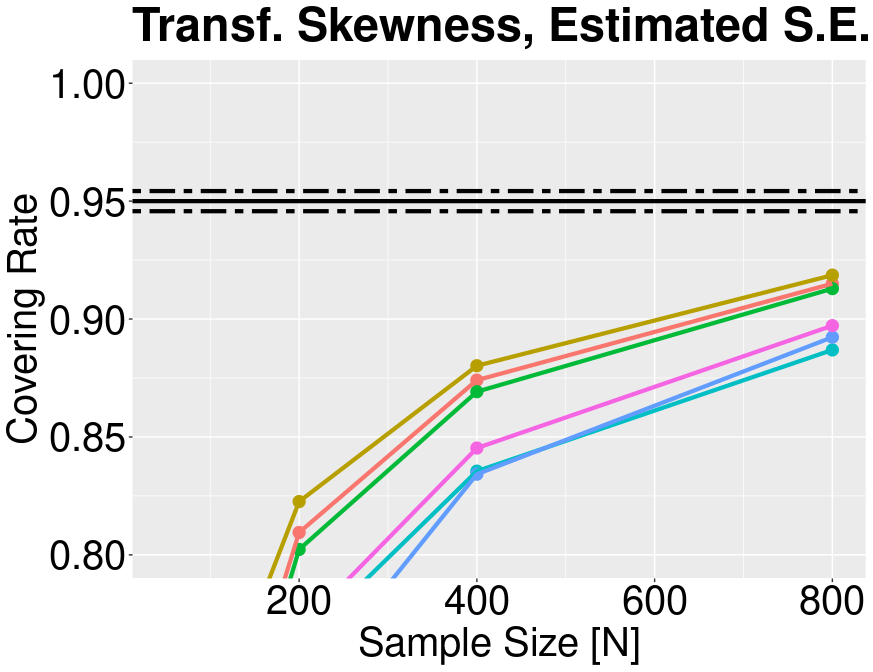}
	\caption{Simulations of coverage of bias corrected SCBs for skewness  and transformed skewness for Model A.  The black dashed lines are a 95\% confidence intervals for the nominal level $0.95$. \label{fig:SCBsModelA-skewBias}}
\end{center}
\end{figure}

\begin{figure}[h]
\begin{center}
\includegraphics[trim=0 0 0 0,clip,width=1.5in]{\figurepath 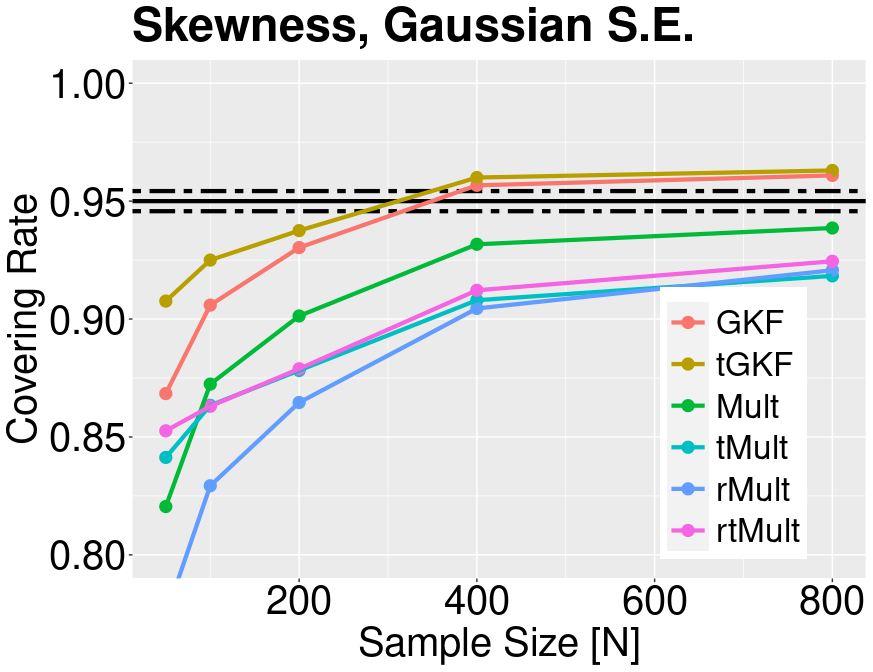}
\includegraphics[trim=0 0 0 0,clip,width=1.5in]{\figurepath 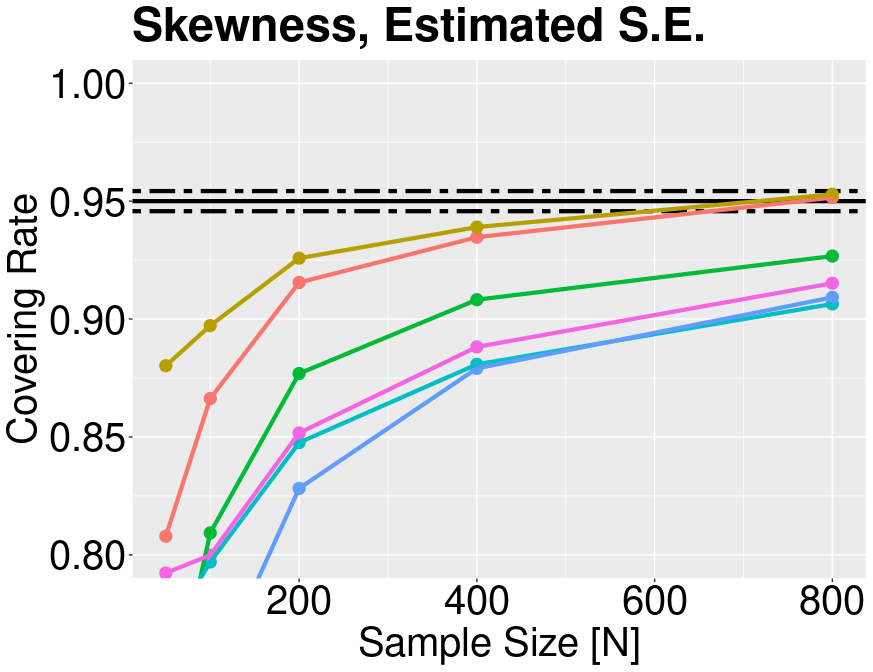}
\includegraphics[trim=0 0 0 0,clip,width=1.5in]{\figurepath 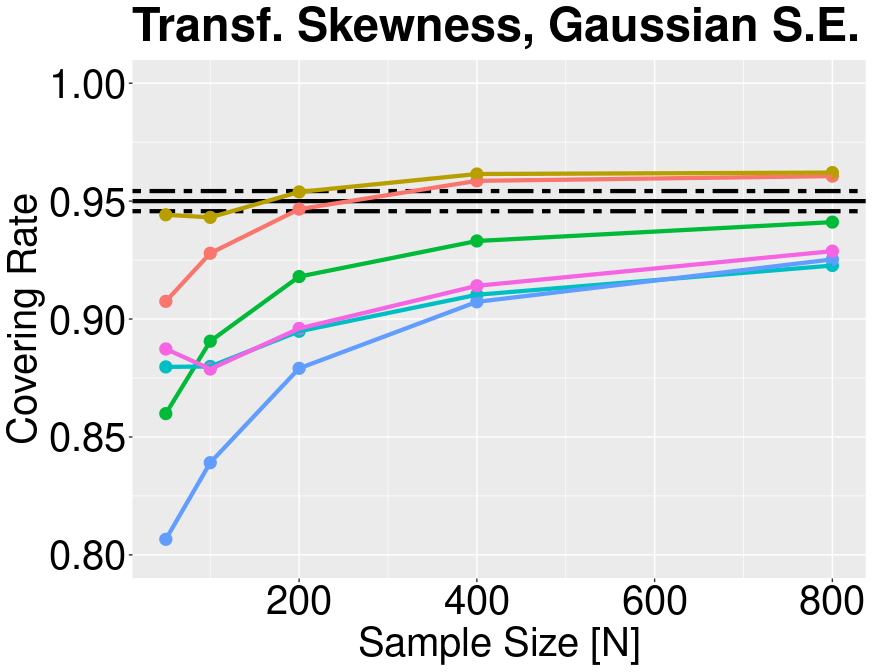}
\includegraphics[trim=0 0 0 0,clip,width=1.5in]{\figurepath 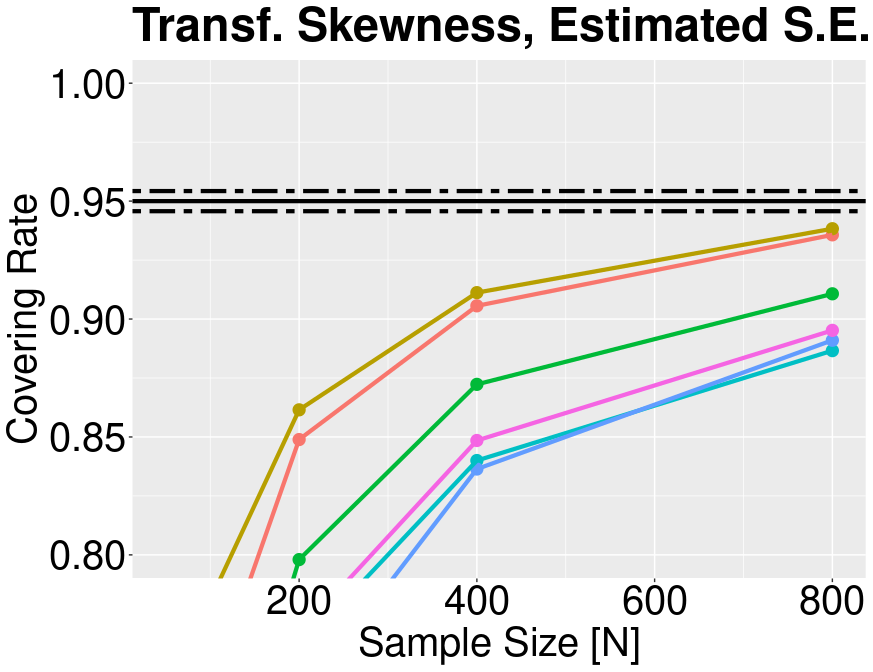}
	\caption{Simulations of coverage of bias corrected SCBs for skewness and transformed skewness for Model B.   The black dashed lines are a 95\% confidence intervals for the nominal level $0.95$. \label{fig:SCBsModelB-skewBias}}
\end{center}
\end{figure}

\begin{figure}[h]
\begin{center}
\includegraphics[trim=0 0 0 0,clip,width=1.5in]{\figurepath 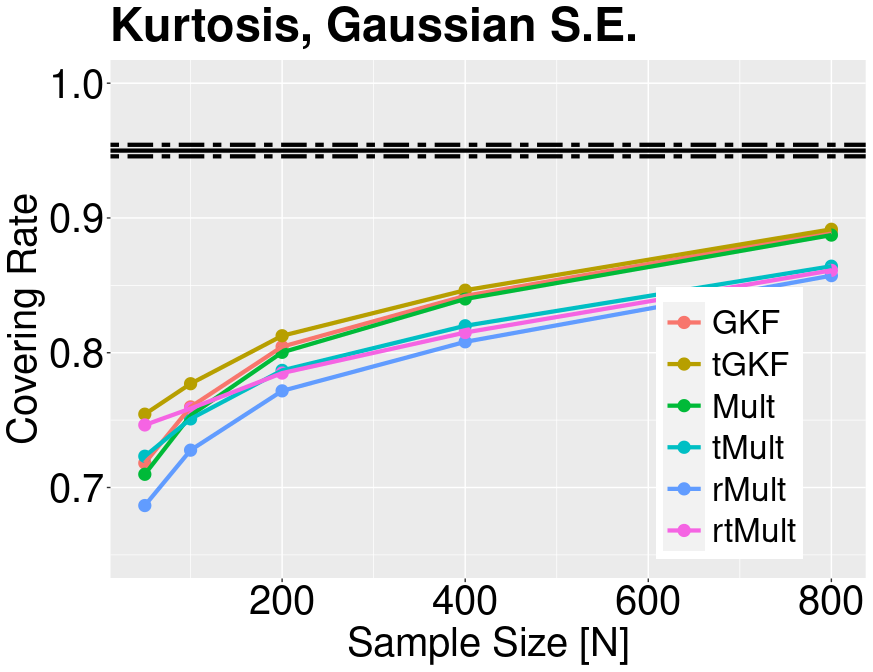}
\includegraphics[trim=0 0 0 0,clip,width=1.5in]{\figurepath 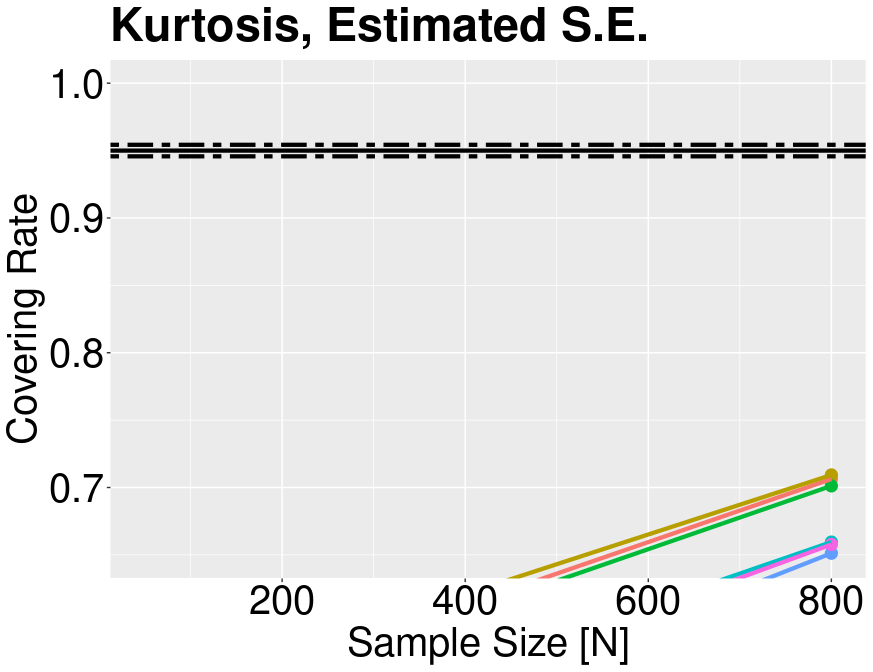}
\includegraphics[trim=0 0 0 0,clip,width=1.5in]{\figurepath 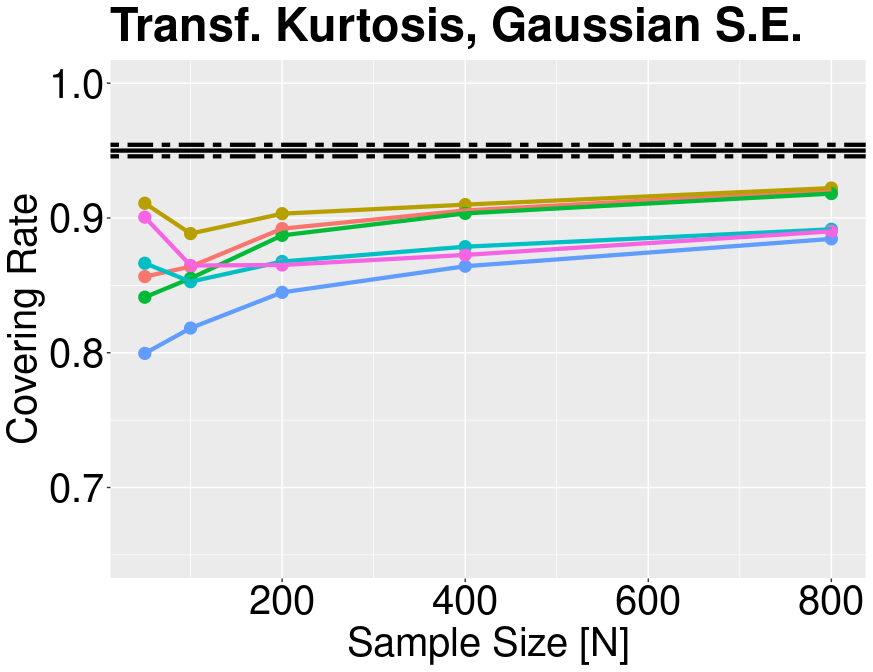}
\includegraphics[trim=0 0 0 0,clip,width=1.5in]{\figurepath 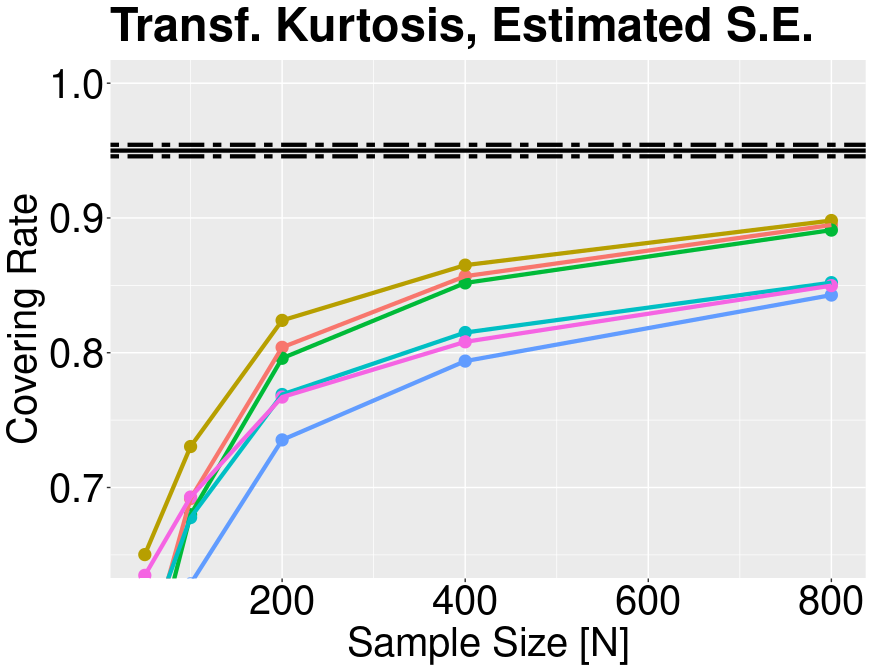}
	\caption{Simulations of coverage of bias corrected  SCBs for kurtosis and transformed kurtosis for Model A.   The black dashed lines are a 95\% confidence intervals for the nominal level $0.95$. \label{fig:SCBsModelA-kurtBias}}
\end{center}
\end{figure}

\begin{figure}[h]
\begin{center}
\includegraphics[trim=0 0 0 0,clip,width=1.5in]{\figurepath 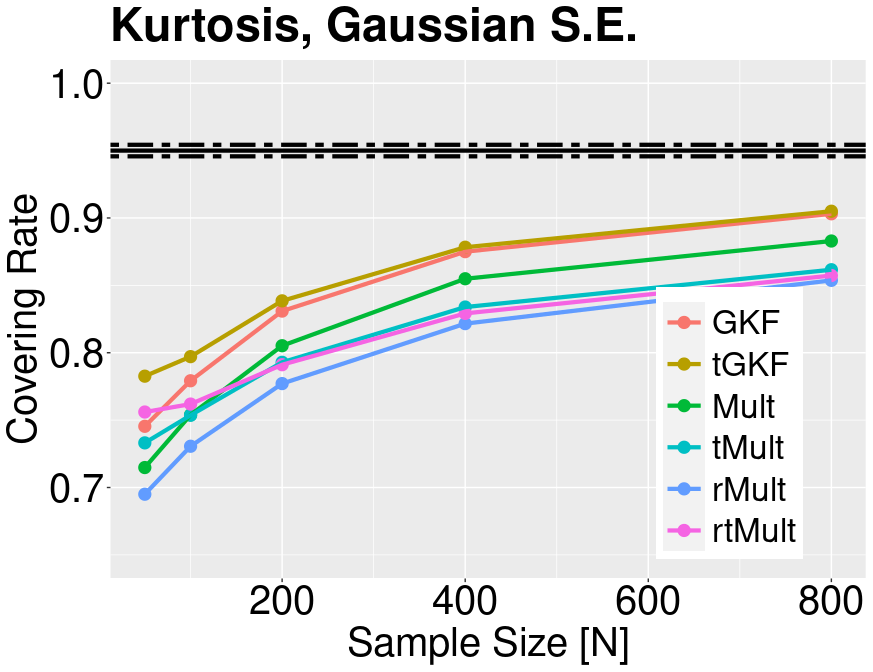}
\includegraphics[trim=0 0 0 0,clip,width=1.5in]{\figurepath 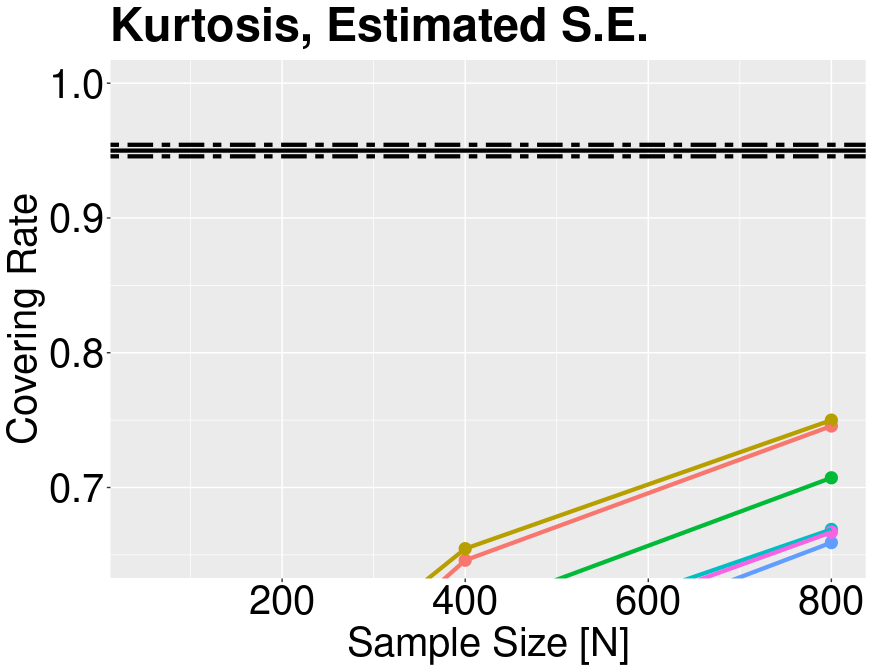}
\includegraphics[trim=0 0 0 0,clip,width=1.5in]{\figurepath 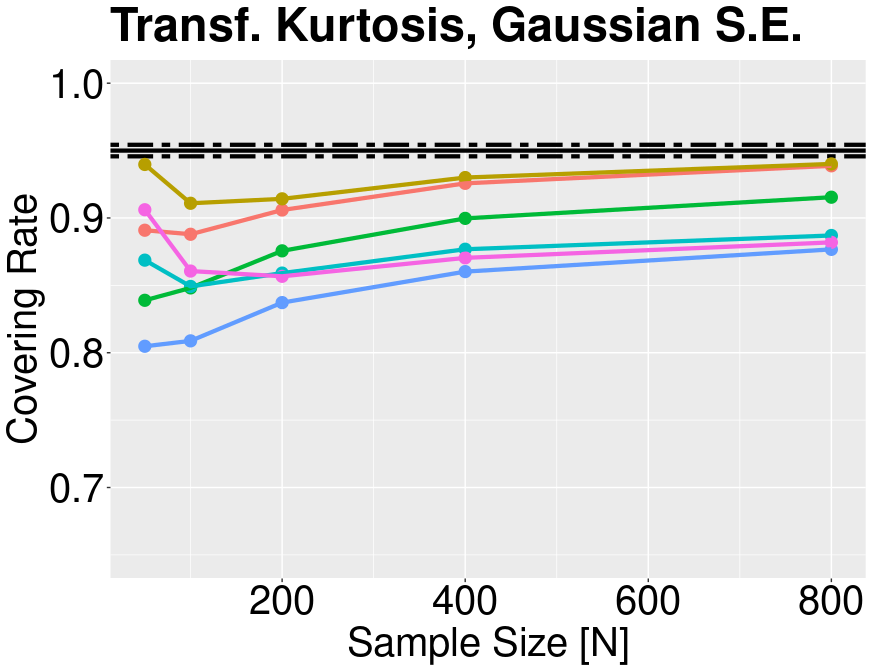}
\includegraphics[trim=0 0 0 0,clip,width=1.5in]{\figurepath 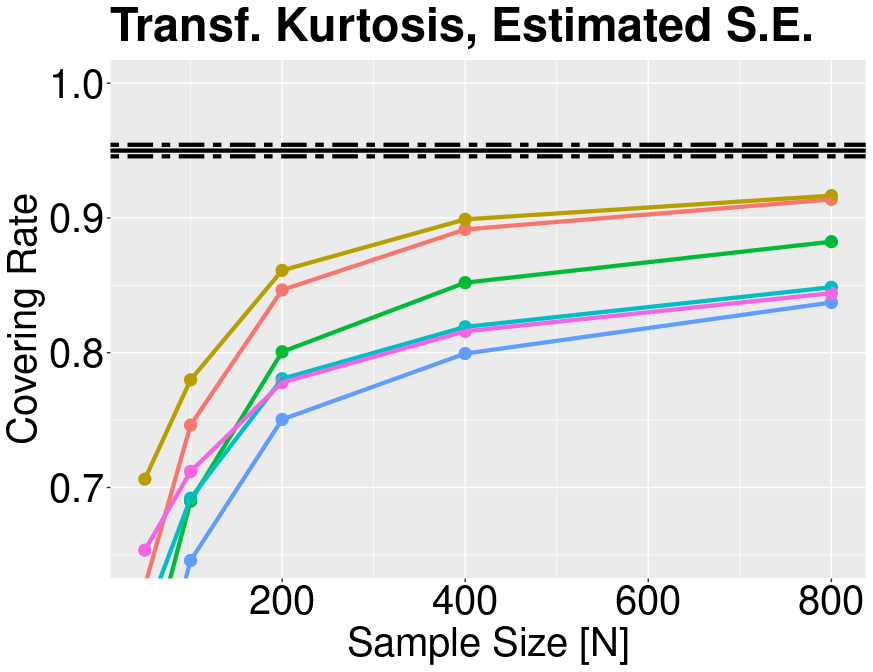}
	\caption{Simulations of coverage of bias corrected SCBs for kurtosis and transformed kurtosis for Model B. The black dashed lines are
		a 95\% confidence intervals for the nominal level $0.95$.	\label{fig:SCBsModelB-kurtBias}}
\end{center}
\end{figure}

\section{Simulations for SCBs for Skewness and Kurtosis in Model C}\label{app:modelC}

\begin{figure}[h]
\begin{center}
\includegraphics[trim=0 0 0 0,clip,width=1.5in]{\figurepath 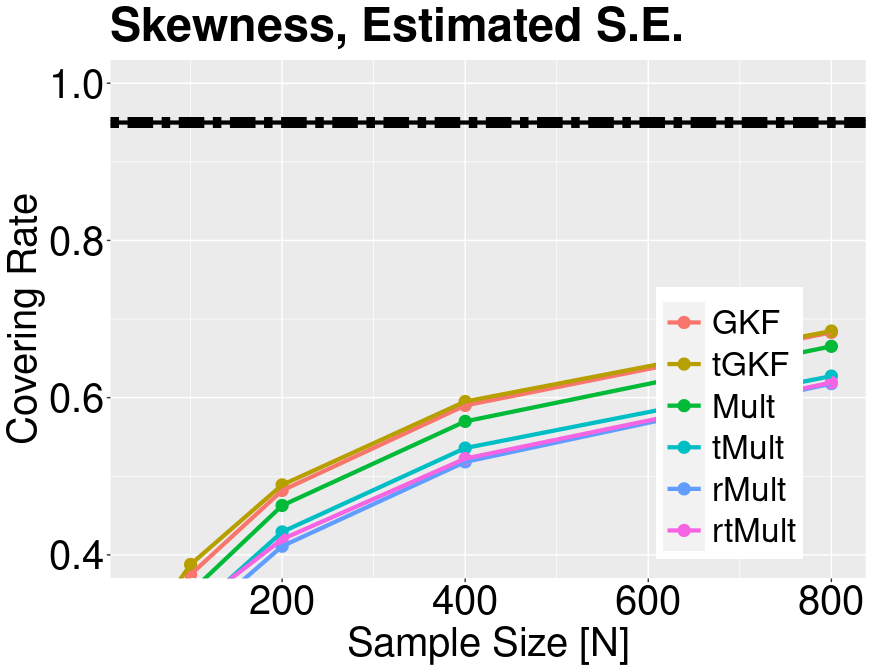}
\includegraphics[trim=0 0 0 0,clip,width=1.5in]{\figurepath 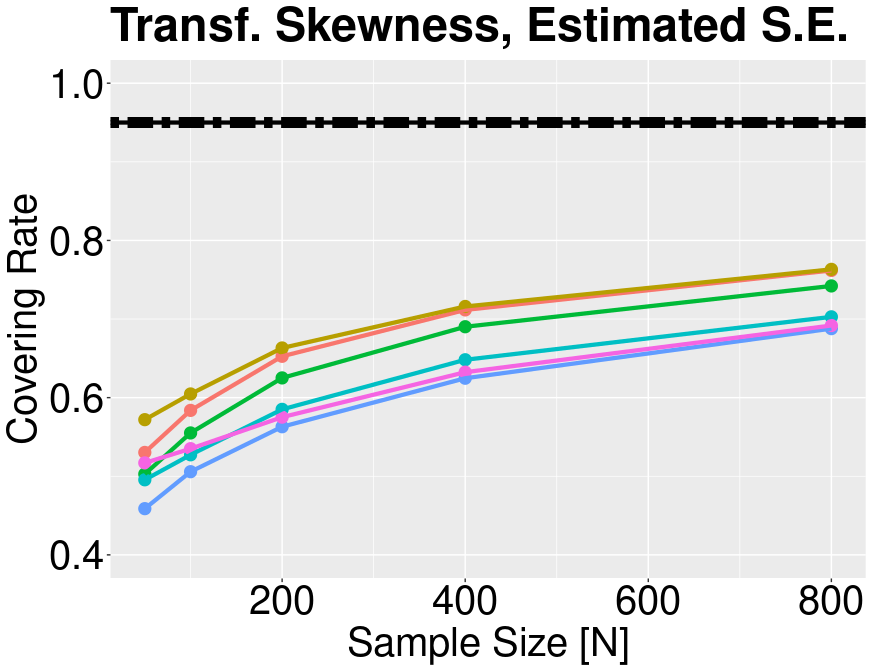}
\includegraphics[trim=0 0 0 0,clip,width=1.5in]{\figurepath 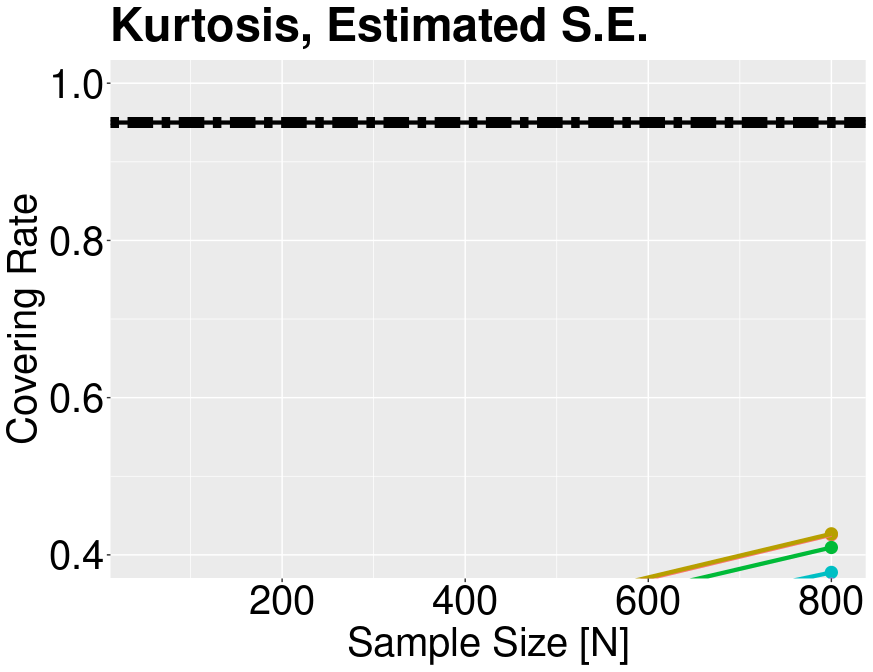}
\includegraphics[trim=0 0 0 0,clip,width=1.5in]{\figurepath 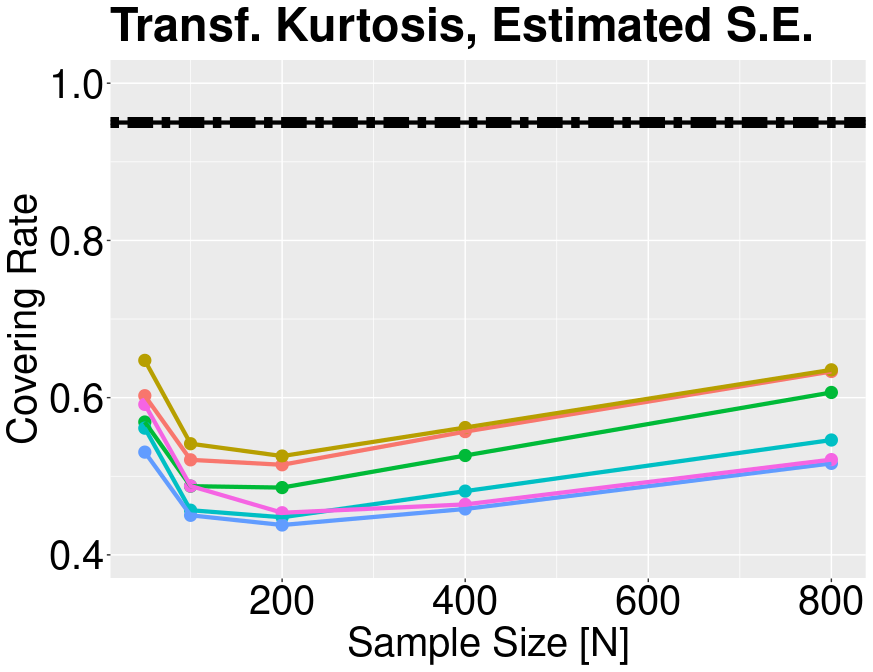}
	\caption{Simulations of coverage rates of SCBs for Model $C$.  In these simulations no bias correction is applied.  The black dashed lines are
		95\% confidence intervals for the nominal level $0.95$.\label{fig:SCBsModelC-skew}}
\end{center}
\end{figure}

\begin{figure}[h]
\begin{center}
\includegraphics[trim=0 0 0 0,clip,width=1.5in]{\figurepath 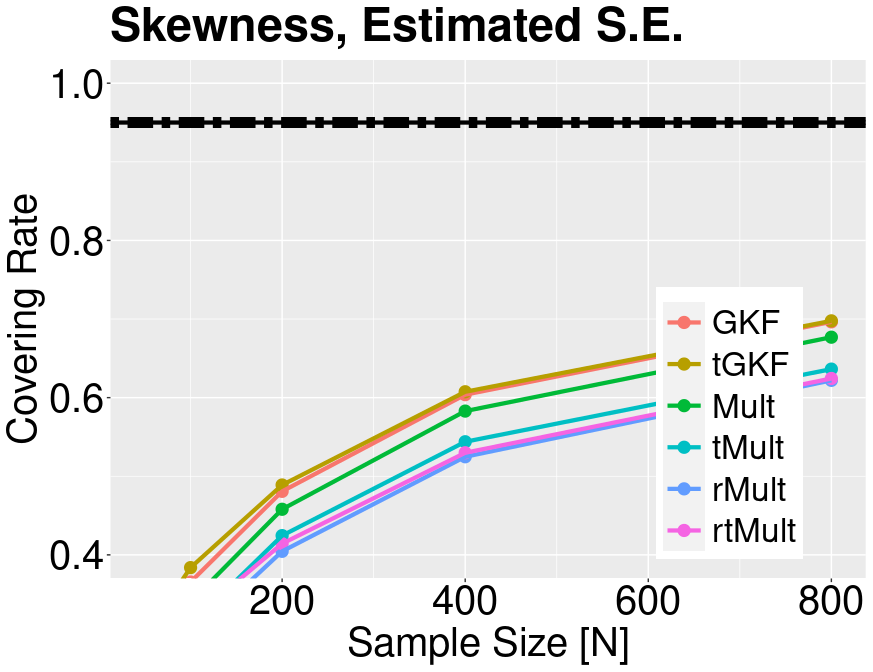}
\includegraphics[trim=0 0 0 0,clip,width=1.5in]{\figurepath 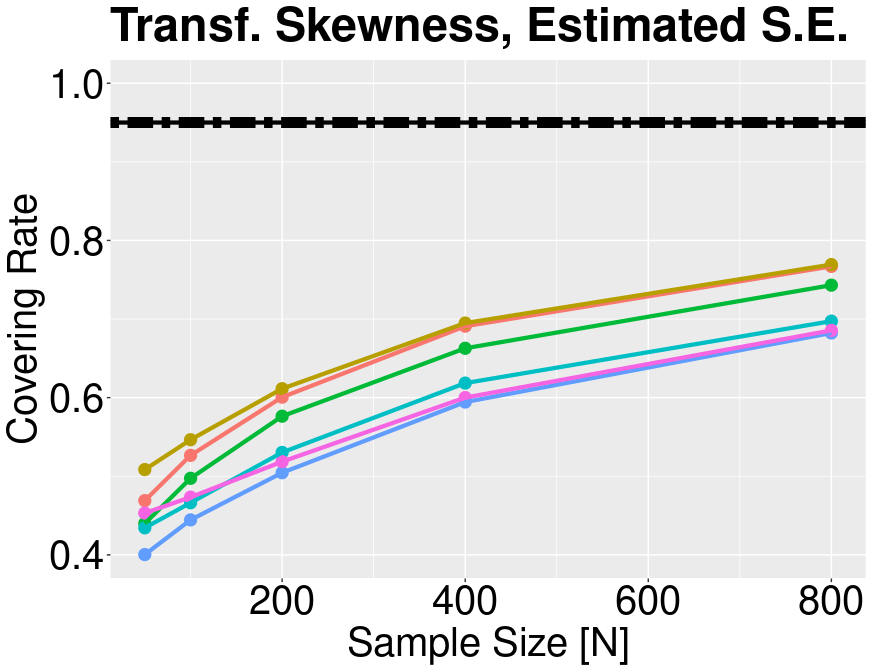}
\includegraphics[trim=0 0 0 0,clip,width=1.5in]{\figurepath 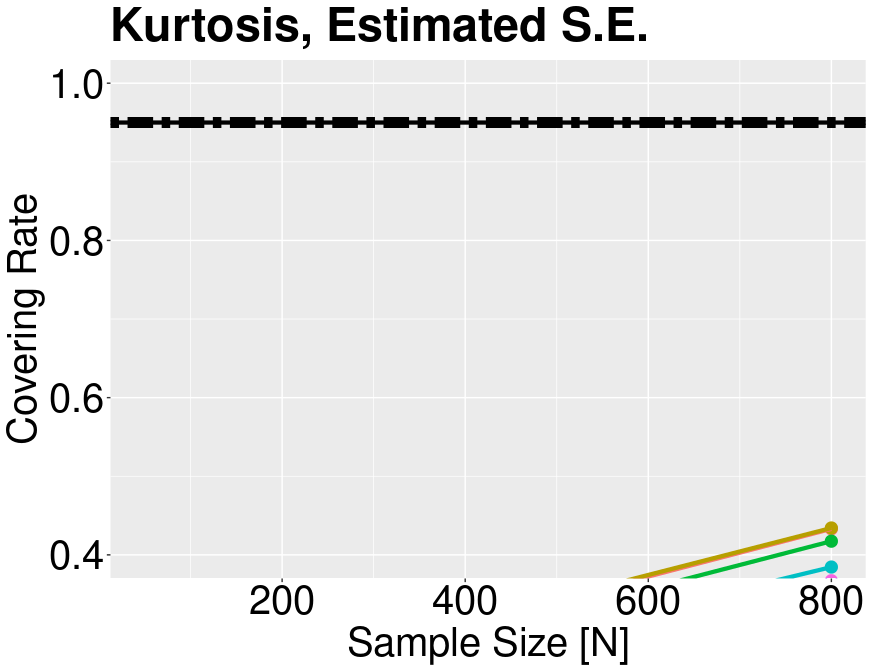}
\includegraphics[trim=0 0 0 0,clip,width=1.5in]{\figurepath 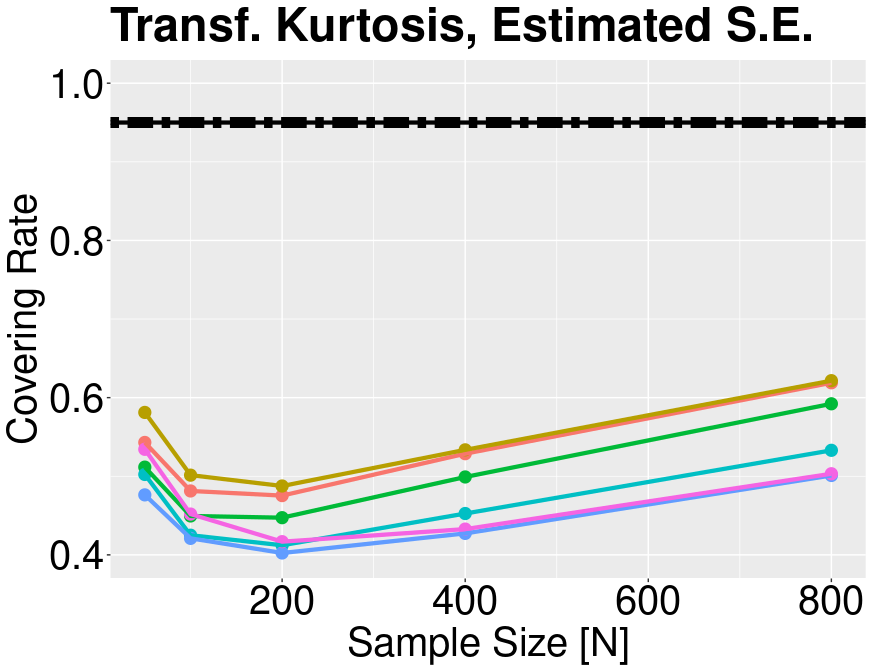}
	\caption{Simulations of coverage of SCBs for skewness, transformed skewness, kurtosis and transformed kurtosis (from left to right) for Model C.   The black dashed lines are
		a 95\% confidence intervals for the nominal level $0.95$.\label{fig:SCBsModelC-skewBias}}
\end{center}
\end{figure}

The simulations of coverage rates of SCBs for Model $C$ for (transformed) skewness and (transformed) kurtosis 
are reported in Figures \ref{fig:SCBsModelC-skew} and \ref{fig:SCBsModelC-skewBias}.
Here we do not report simulations using the Gaussian standard error since it results in very low coverage rates. The pointwise standard errors of (transformed) skewness and (transformed) kurtosis for Model $ C $ are very different from the exact Gaussian quantities due to the high non-Gaussianity of Model $C$.
As shown in Figures \ref{fig:SCBsModelC-skew} and \ref{fig:SCBsModelC-skewBias} the SCBs using the estiamted standard error from the delta residuals require $ N $ to be much larger than $ 800 $ to converge to the nominal value.
This is partially because the variance estimate of the skewness estimator requires even longer to converge under the non-Gaussian noise model than for Gaussian models. A solution for this extreme case to construct SCBs is left for future work.

\end{document}